\documentclass[reqno,twoside]{amsart}

\usepackage{amsfonts}
\usepackage{amsmath,amssymb,amsthm}
\usepackage{mathtools}
\usepackage{etoolbox}
\usepackage{color}
\usepackage[english]{babel}
\usepackage{hyperref}
\usepackage{microtype}
\usepackage{fullpage}
\usepackage{floatrow}
\usepackage[utf8]{inputenc}
\usepackage[T1]{fontenc}
\usepackage{lmodern}
\usepackage[leftcaption]{sidecap}
\sidecaptionvpos{figure}{m}

\DisableLigatures{encoding = *, family = * }

\usepackage{float}

\newtheorem{theorem}{Theorem}[section]
\newtheorem{definition}[theorem]{Definition}

\newtheorem{lemma}[theorem]{Lemma}
\newtheorem{proposition}[theorem]{Proposition}
\newtheorem{remark}[theorem]{Remark}

\numberwithin{equation}{section}
\numberwithin{figure}{section}

\newcommand{\RR}{\mathbb{R}}
\newcommand{\NN}{\mathbb{N}}

\newcommand{\norm}[2]{{\left\|#1\right\|}_{#2}}
\newcommand{\fl}[2]{(-\Delta)^#1#2}

\newcommand{\Om}{\Omega}

\title[Control and numerics of fractional diffusion]{Control and numerical approximation of fractional diffusion equations}

\author{Umberto Biccari\textsuperscript{\,$\ast$}}  
\address{\textsuperscript{$\ast$}\, [1] Chair of Computational Mathematics, Fundaci\'on Deusto, Avenida de las Universidades 24, 48007 Bilbao, Basque Country, Spain} 
\address{[2]\,Facultad de Ingenier\'ia, Universidad de Deusto, Avenida de las Universidades 24, 48007 Bilbao, Basque Country, Spain.}
\email{umberto.biccari@deusto.es, u.biccari@gmail.com}

\thanks{This project has received funding from the European Research Council (ERC) under the European Union’s Horizon 2020 research and innovation programme (grant agreement NO: 694126-DyCon). The work of the three authors is supported by the Air Force Office of Scientific Research (AFOSR) under Award NO: FA9550-18-1-0242. The work of MW is supported by the US Army Research Office (ARO) under Award NO: W911NF-20-1-0115. The work of U.B. and E.Z. is supported by Grant MTM2017-92996-C2-1-R COSNET of MINECO (Spain) and by the Elkartek grant KK-2020/00091 CONVADP of the Basque government. The work of E.Z. is funded by the Alexander von Humboldt-Professorship program, the European Unions Horizon 2020 research and innovation programme under the Marie Sklodowska-Curie grant agreement No.765579-ConFlex and the Transregio 154 Project ‘‘Mathematical Modelling, Simulation and Optimization Using the Example of Gas Networks’’, project C08, of the German DFG} 

\author{Mahamadi Warma\textsuperscript{\,$\dagger$}}  
\address{\textsuperscript{$\dagger$}\,Department of Mathematical Sciences and the Center for Mathematics and Artificial Intelligence, George Mason University.  Fairfax VA 22030 (USA)}
\email{	mwarma@gmu.edu}

\author{Enrique Zuazua\textsuperscript{\,$\ddagger$}}
\address{\textsuperscript{$\ddagger$}\, [1] Chair for Dynamics, Control and Numerics, Alexander von Humboldt-Professorship, Department of Data Science, Friedrich-Alexander-Universit\"at Erlangen-N\"urnberg, 91058 Erlangen, Germany.}
\address{[2] Chair of Computational Mathematics, Fundaci\'on Deusto, Avenida de las Universidades 24, 48007 Bilbao, Basque Country, Spain.} 
\address{[3] Departamento de Matem\'aticas, Universidad Aut\'onoma de Madrid, 28049 Madrid, Spain.}
\email{enrique.zuazua@fau.de}


\begin{document}
	
\begin{abstract}
The aim of this work is to give a broad panorama of the control properties of fractional diffusive models from a numerical analysis and simulation perspective. We do this by surveying several research results we obtained in the last years, focusing in particular on the numerical computation of controls, though not forgetting to recall other relevant contributions which can be currently found in the literature of this prolific field. Our reference model will be a non-local diffusive dynamics driven by the fractional Laplacian on a bounded domain $\Omega$. The starting point of our analysis will be a Finite Element approximation for the associated elliptic model in one and two space-dimensions, for which we also present error estimates and convergence rates in the $L^2$ and energy norm. Secondly, we will address two specific control scenarios: firstly, we consider the standard interior control problem, in which the control is acting from a small subset $\omega\subset\Omega$. Secondly, we move our attention to the exterior control problem, in which the control region $\mathcal O\subset\Omega^c$ is located outside $\Omega$. This exterior control notion extends boundary control to the fractional framework, in which the non-local nature of the models does not allow for controls supported on $\partial\Omega$. We will conclude by discussing the interesting problem of simultaneous control, in which we consider families of parameter-dependent fractional heat equations and we aim at designing a unique control function capable of steering all the different realizations of the model to the same target configuration. In this framework, we will see how the employment of stochastic optimization techniques may help in alleviating the computational burden for the approximation of simultaneous controls.  Our discussion is complemented by several open problems related with fractional models which are currently unsolved and may be of interest for future investigation.
\end{abstract}
	
\maketitle

\section{Introduction}

A non-local PDE is a particular type of differential equation involving integral or pseudo-differential terms. For this reason, these models are often referred as integro-differential equations. 

In recent years, there has been a growing interest on non-local models because of their relevance in a wide spectrum of practical applications. Indeed, there is a plethora of situations in which a non-local equation gives a significantly better description than a local PDE of the problem one wants to analyze. A widely studied class of non-local models involves fractional order operators, which have nowadays emerged as a modeling alternative in various branches of science. They usually describe anomalous diffusion. Typical examples in which non-local or fractional equations appear are models in turbulence (\cite{bakunin2008turbulence}), population dynamics (\cite{de2002size}), image processing (\cite{gilboa2008nonlocal}), laser design (\cite{longhi2015fractional}), and porous media flow (\cite{vazquez2012nonlinear}). Besides, a number of stochastic models associated with fractional operators have been introduced in the literature for explaining anomalous diffusion. Among them we quote the fractional Brownian motion, the continuous time random walk, the L\'evy flights, the Schneider gray Brownian motion, and more generally, random walk models based on evolution equations of single and distributed fractional order in  space (see e.g. \cite{dubkov2008levy,gorenflo2007continuous,mandelbrot1968fractional,schneider1990stochastic}). In general, a fractional diffusion operator corresponds to a diverging jump length variance in the random walk. Finally, we can refer to \cite{antil2019sobolev,weiss2018fractional} for the relevance of fractional operators in geophysics and imaging science.

This surge of interest towards non-local and fractional models has opened a very challenging field in the applied mathematical research, since most of the existing techniques in PDE analysis were not adapted to treat non-local effects. 

In the wide spectrum of non-local and fractional models, of particular interest are the ones involving the fractional Laplacian. From a mathematical perspective, there is nowadays a well established and rich literature on this operator and its employment in PDE models. Among many others contributions, we remind here the works \cite{biccari2018poisson,biccari2017local,biccari2018local,caffarelli2007extension,caffarelli2009regularity,leonori2015basic,ros2014dirichlet,servadei2012mountain,servadei2013variational,warma2015the}. Besides, it is well known that the fractional Laplacian is the generator of s-stable L\'evy processes, and it is often used in stochastic models with applications, for instance, in mathematical finance (\cite{levendorskii2004pricing,pham1997optimal}).

In addition to that, control problems in the non-local setting have been largely considered in recent years. An incomplete bibliography at this respect includes \cite{antil2020controllability,biccari2020internal,biccari2019controllability,biccari2019null,biccari2020null,biccari2020positive,chaves2017controllability,fernandez2016null,lu2017null,micu2006controllability,miller2006controllability,warma2019approximate,warma2020analysis,warma2018null}.

The aim of this chapter is to give a broad panorama of the control properties of fractional diffusive models from a numerical analysis and simulation perspective. We will do this by surveying several research results we obtained in the last years, focusing in particular on aspects related with the numerical computation of the controls, though not without recalling other relevant contributions which can be currently found in the literature of this prolific field. In more details, we will consider the following specific situations.
\begin{itemize}
	\item[1.] {\bf Elliptic problems}: we will start in Section \ref{sec:2} by considering the elliptic problem associated with the fractional Laplace operator, for which we will present a Finite Element (FE) approximation in one and two space dimensions. This FE scheme will be at the basis of all our numerical simulations in the sections afterward. Our discussion in this section is based on the contributions \cite{acosta2017short,acosta2017fractional,biccari2019controllability}. 
	
	\item[2.] {\bf Interior control of fractional heat equations}: secondly, we will consider in Section \ref{sec:3} the interior control of fractional heat equations. We will start by recalling the theoretical controllability results we have obtained in \cite{biccari2019controllability,biccari2020positive,biccari2021null}. In a second moment, we will deal with the numerical computation of controls in one and two space dimensions by means of the penalized Hilbert Uniqueness Method. 
		
	\item[3.] {\bf Exterior control of fractional heat equations}: thirdly, we will deal in Section \ref{sec:4} with the exterior control properties of fractional heat equations. As we shall see, this is the equivalent of the boundary controllability in the framework of the fractional Laplacian. Also in this case, we will first remind the theoretical controllability results we have obtained in \cite{antil2020controllability,warma2018null}. After that, we will focus on the numerical computation of controls and discuss some relevant differences with respect to the interior control problem of Section \ref{sec:3}.
	
	\item[4.] {\bf Simultaneous control of parameter-dependent fractional heat equations}: finally, we will address in Section \ref{sec:5} the simultaneous control of parameter-dependent fractional heat equations, in which we aim at designing a unique parameter-independent control capable of steering different realizations of the same fractional dynamics to zero in finite time. Our discussion here will focus on the associated optimal control problem and its computational aspects. In particular, we will propose the employment of stochastic optimization algorithms for the computation of the simultaneous controls, and discuss the advantages and disadvantages of these techniques with respect to some more classical ones (namely, Gradient Descent and Conjugate Gradient) typically applied in PDE control. 	
\end{itemize}

Finally, this chapter will be complemented with Appendices \ref{appendixA} and \ref{appendixB}, in which we gather some results on fractional Sobolev spaces and the fractional Laplacian which will be used in our analysis.

\section{Finite Element approximation of the fractional Laplace operator}\label{sec:2}

In this section, we give an abridged presentation of the fundamental aspects required to convey a complete FE analysis of the following elliptic problem
\begin{align}\label{eq:Sec2_DP}
	\begin{cases}
		\fl{s}{u} = f & \mbox{ in }\;\Omega,
		\\
		u = 0 & \mbox{ in }\;\Omega^c,
	\end{cases}
\end{align}
where $\Omega\subset\RR^d$ ($d=1,2$) is a bounded and $C^{1,1}$ domain, $\Omega^c:=\RR^d\setminus\Omega$, $f$ is a smooth enough function whose regularity will be specified later, and where for all $s\in (0,1)$ we denote with $\fl{s}{}$ the fractional Laplace operator defined as 
\begin{align}\label{eq:Sec2_fl}
	\fl{s}{u}(x) := C_{d,s} P.V. \int_{\mathbb R^d}\frac{u(x)-u(y)}{|x-y|^{d+2s}}\,dy,
\end{align} 
$C_{d,s}$ being an explicit normalization constant (see \eqref{eq:App1_FLO}).

Numerical approximation for the fractional Laplacian on bounded domains has been extensively addressed in the last years. In \cite{huang2014numerical}, the authors proposed a method combining finite differences with numerical quadrature to obtain a discrete convolution operator which accurately approximates the singular integral defining the fractional Laplacian. Nevertheless, it is noteworthy that the convergence of the proposed algorithm is proved assuming that solutions to \eqref{eq:Sec2_DP} are of class $C^4$, which is generally not the case.

These regularity requirements can be avoided if we discretize the fractional Laplacian \eqref{eq:Sec2_fl} via a FE approach, which is based on the weak variational formulation associated with \eqref{eq:Sec2_DP} (see \ref{def:App1_weakSol}). This issue was firstly addressed in \cite{acosta2017short,acosta2017fractional,borthagaray2017laplaciano}, in the 2-D case, by combining techniques borrowed from the theory of the Boundary Element Method (\cite{sauter2010boundary}) together with an appropriate numerical treatment of the integrals involving the unbounded domain $\Omega^c$. An alternative yet complementary approach, allowing to treat also 3-D problems, has been proposed in \cite{bonito2019numerical}, where the authors based their FE approximation on a different (but equivalent) variational formulation of the Dirichlet problem. Moreover, in \cite{biccari2019controllability}, motivated by control purposes, we specifically focused on the one-dimensional fractional Laplacian. In this setting, we showed that the stiffness matrix discretizing the operator can be computed entirely offline, without need of any numerical integration, with entries depending only on their position, fractional exponent $s$ and the mesh-size $h$. 

For completeness, we shall also mention some works dealing with the so-called \textit{spectral} Dirichlet fractional Laplacian, that is, the fractional $s$-power of the realization in $L^2(\Omega)$ of the Laplace operator $-\Delta$ with null Dirichlet boundary conditions. To discretize this operator, in \cite{nochetto2015pde} the authors proposed a FE approach combined with the famous \textit{Caffarelli-Silvestre extension} (see \cite{caffarelli2007extension}) which allows to work in a local framework. Nevertheless, we have to stress that the spectral fractional Laplacian is different from the integral one that we consider in this work (see e.g. \cite{servadei2014spectrum} and the references therein).

In what follows, we will summarize the main steps for the FE approximation of the fractional Laplacian in space dimension $d=1,2$, as presented in \cite{acosta2017short,acosta2017fractional,biccari2019controllability,borthagaray2017laplaciano}. We will start by discussing in Subsection \ref{subs:Sec2_stiffness} the construction of the stiffness matrix. In the one-dimensional case, we will follow the procedure of \cite{biccari2019controllability}, while the presentation of the two-dimensional problem is based on \cite{acosta2017short,acosta2017fractional,borthagaray2017laplaciano}. Secondly, we will give in Subsection \ref{subs:Sec2_error} an overview of the error analysis and the convergence rate of the FE scheme in $H^s$ and $L^2$. Finally, Subsection \ref{subs:Sec2_simul} will be devoted to some numerical experiments.

\subsection{Computation of the stiffness matrix}\label{subs:Sec2_stiffness}

\begin{figure}[h]
	\centering
	\includegraphics[scale=0.7]{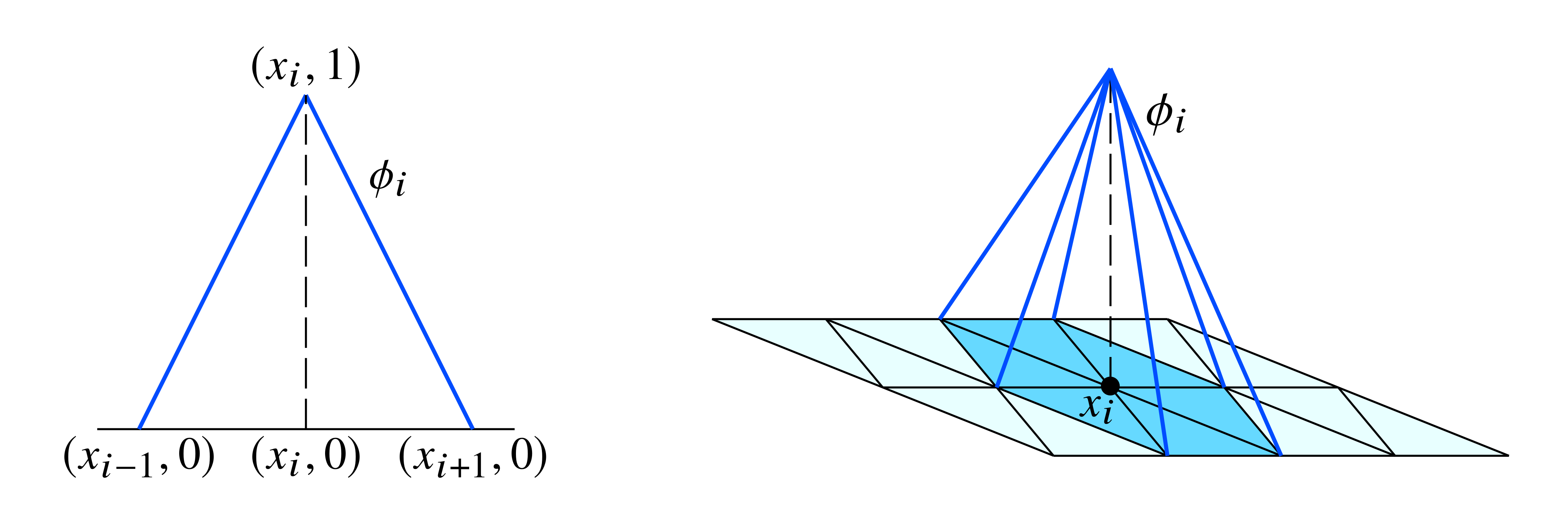}
	\caption{Basis functions $\phi_i(x)$ in space dimension $d=1$ (left) and $d=2$ (right).}\label{fig:Sec2_basis_int}
\end{figure}

As in any FE scheme, the first step is to introduce the variational formulation associated with the problem we aim to solve. In our case, 
this variational formulation reads as follows (see Definition \ref{def:App1_weakSol}): find $u\in H_0^s(\Omega)$ such that, for all $v\in H_0^s(\Omega)$,	
\begin{align}\label{eq:Sec2_WF}
	a(u,v):= \frac{C_{d,s}}{2}\int_{\RR^d}\int_{\RR^d}\frac{(u(x)-u(y))(v(x)-v(y))}{|x-y|^{d+2s}}\;dxdy = \int_\Omega fv\,dx.
\end{align}
Let $\mathfrak M = \{T_i\}_{i=1}^N$ be a partition of the domain $\Omega$ (i.e. $\overline{\Omega}=\bigcup_{i=1}^N T_i$) composed by $N$ elements as follows:
\begin{itemize}
	\item In dimension $d=1$, we consider a uniform partition of $\Omega=(-1,1)$ 
	\begin{align*}
		-1 = x_0<x_1<\ldots <x_i<x_{i+1}<\ldots<x_{N+1}=1\,,
	\end{align*}
	with $x_{i+1}=x_i+h$, $i=0,\ldots N$, and we denote $T_i:=[x_i,x_{i+1}]$. 

	\item In dimension $d=2$, we consider triangular elements $T_i$, $i = 1,\ldots,n$. We indicate with $h_i$ the diameter of the element $T_i$ and with $\rho_i$ its inner radius, i.e. the diameter of the largest ball contained in $T_i$. We define 
	\begin{align*}
		h = \max_{i\in\{1,\ldots,N\}} h_i.
	\end{align*} 
	Moreover, we require that the triangulation satisfies the usual regularity and local uniformity conditions:
	\begin{equation}\label{eq:Sec2_mesh2D}
		\begin{array}{l} 
			\mbox{ there exists } \sigma > 0 \mbox{ s.t. } h_i\leq \sigma\rho_i \mbox{ for all } i=1,\ldots,N,
			\\
			\mbox{ there exists } \kappa > 0 \mbox{ s.t. } h_i\leq \kappa h_j \mbox{ for all } i,j = 1,\ldots,N,\; T_i\cap T_j = \emptyset. 
		\end{array} 
	\end{equation}
	Naturally the second condition is a consequence of the first one. In this way $\kappa$ can be expressed in terms of $\sigma$.
\end{itemize}
Consider the discrete space 
\begin{align*}
	\mathbb V :=\Big\{v\in H_0^s(\Omega)\,:\, \left. v\,\right|_{T_i}\in \mathcal{P}^1\Big\},
\end{align*} 
$\mathcal{P}^1$ being the space of the continuous and piece-wise linear functions. We approximate \eqref{eq:Sec2_WF} with the following discrete problem: find $u_h\in \mathbb V$ such that
\begin{align}\label{eq:Sec2_FEid}
	\frac{C_{d,s}}{2} \int_{\RR}\int_{\RR}\frac{(u_h(x)-u_h(y))(v(x)-v(y))}{|x-y|^{1+2s}}\,dxdy = \int_\Omega fv\,dx,	
\end{align}
for all $v\in \mathbb V$. Let $\big\{\phi_i\big\}_{i=1}^N$ be the standard nodal basis of $\mathbb V$ corresponding to the internal nodes $\{x_1,\ldots,x_N\}$, that is $\phi_i(x_j)=\delta_{i,j}$ (see Figure \ref{fig:Sec2_basis_int}). By decomposing
\begin{align*}
	u(x) = \sum_{j=1}^N u_j\phi_j(x)\quad\mbox{ and }\quad f(x) = \sum_{j=1}^N f_j\phi_j(x),
\end{align*} 
and taking $v=\phi_i$, \eqref{eq:Sec2_FEid} becomes the linear system $\mathcal A_h {\bf u} = \mathcal M_h {\bf f}$, where ${\bf u} = (u_1,\ldots,u_N)\in\RR^N$, ${\bf f} = (f_1,\ldots,f_N)\in\RR^N$, and the stiffness and mass matrices $\mathcal A_h,\mathcal M_h\in \RR^{N\times N}$ have components
\begin{equation}\label{eq:Sec2_stiffness_nc}
	\begin{array}{ll}
		\displaystyle a_{i,j}=\frac{C_{d,s}}{2} \int_{\RR^d}\int_{\RR^d}\frac{(\phi_i(x)-\phi_i(y))(\phi_j(x)-\phi_j(y))}{|x-y|^{1+2s}}\,dxdy \notag 
		\\
		\displaystyle m_{i,j} = \int_\Omega \phi_i(x)\phi_j(x)\,dx &i,j = 1,\ldots,N.
	\end{array} 
\end{equation}

Notice that $a(\phi_i,\phi_j)$ in \eqref{eq:Sec2_stiffness_nc} provides non-local interactions between the basis functions on the entire real line. Consequently, the stiffness matrix $\mathcal A_h$ will be full.

\begin{remark}
Although in our presentation we focused on uniform meshes, adapted grids could have also been employed. This would allow to treat problems with corner singularities or to cope with boundary layers arising for instance in convection-dominated problems (see \cite{acosta2017fractional}). Nevertheless, for the problems that we are considering here, non-uniform meshes will not be necessary. Besides, we stress that, in the 1-D case, the employment of a uniform mesh actually allows to compute explicitly the entries of the stiffness matrix $\mathcal A_h$, without need of using any quadrature formula (see \cite{biccari2019controllability}).
\end{remark}

We now give an abridged presentation of how to compute the elements of $\mathcal A_h$. To this end, we shall consider the one and two-dimensional cases separately. 

\subsubsection{1-D case}
\begin{SCfigure}[1.2][h]
	\centering
	\includegraphics[scale=0.6]{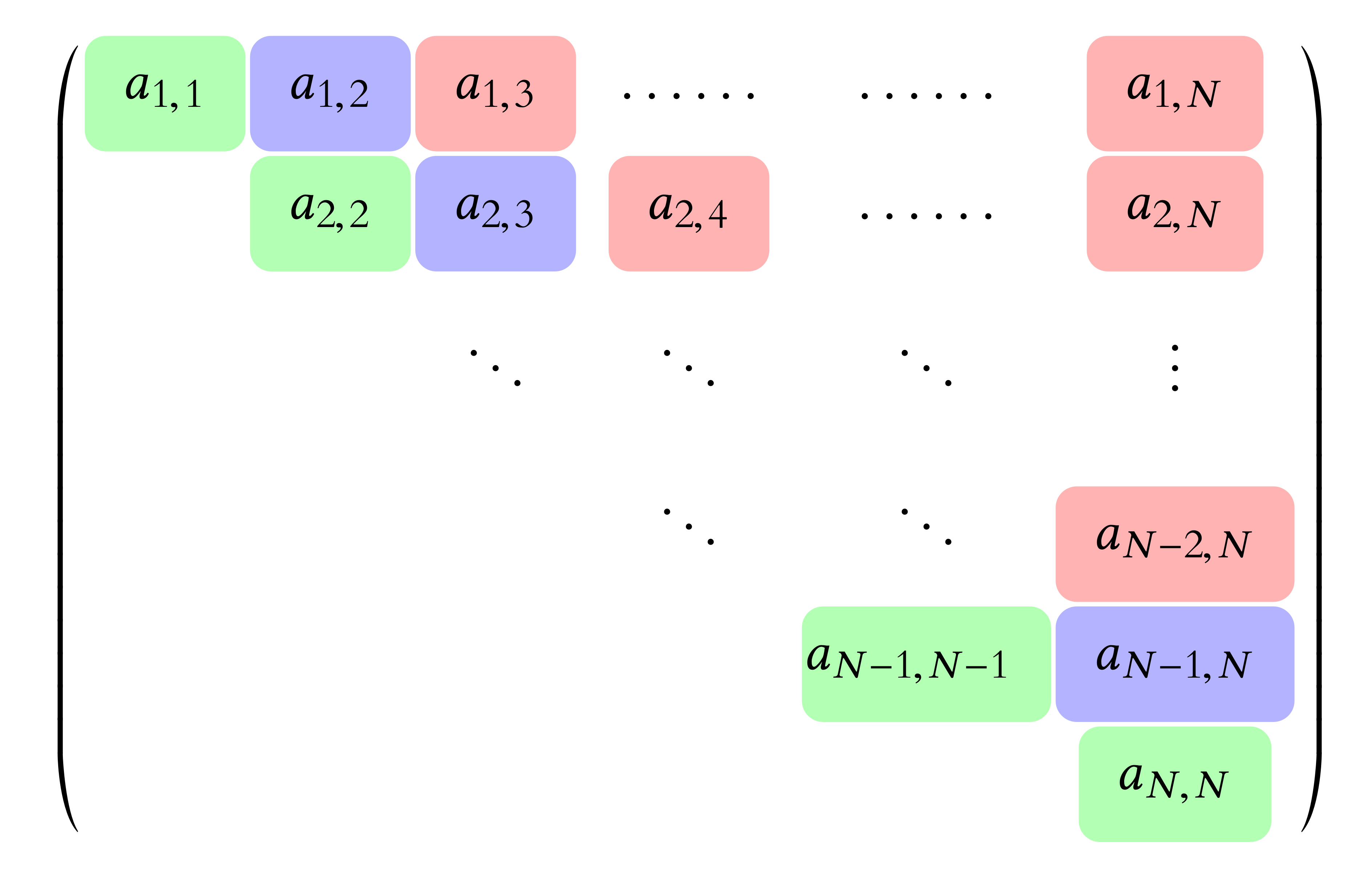}
	\caption{Structure of the stiffness matrix $\mathcal{A}_h$. The red, purple, and green elements are computed through the integrals \eqref{eq:Sec2_Atriangle}, \eqref{eq:Sec2_AdiagonalUpp}, and \eqref{eq:Sec2_Adiagonal}, respectively.}\label{fig:Sec2_matrix_fig}
\end{SCfigure}

In space dimension $d=1$, it is possible to identify three well defined regions in which we have different intersections among the support of the basis functions, thus generating different values of the coefficients $a_{i,j}$. 
\begin{itemize}
	\item[1.] In the upper triangle in red in Figure \ref{fig:Sec2_matrix_fig}, corresponding to $j\geq i+2$, we have $supp(\phi_i)\cap supp(\phi_j) =\emptyset$ and we get
	\begin{align}\label{eq:Sec2_Atriangle}
		a_{i,j}=-C_{1,s} \int_{x_{j-1}}^{x_{j+1}}\int_{x_{i-1}}^{x_{i+1}}\frac{\phi_i(x)\phi_j(y)}{|x-y|^{1+2s}}\,dxdy.
	\end{align}
	\item[2.] In the upper diagonal in purple in Figure \ref{fig:Sec2_matrix_fig}, corresponding to $j=i+1$, we have
	\begin{align}\label{eq:Sec2_AdiagonalUpp}
		a_{i,i+1}= & \frac{C_{1,s}}{2}\int_{\RR}\int_{\RR}\frac{(\phi_i(x)-\phi_i(y))(\phi_{i+1}(x)-\phi_{i+1}(y))}{|x-y|^{1+2s}}\,dxdy.
	\end{align}
	\item[3.] In the diagonal in green in Figure \ref{fig:Sec2_matrix_fig}, corresponding to $j=i$, we have
	\begin{align}\label{eq:Sec2_Adiagonal}
		a_{i,i}= & \frac{C_{1,s}}{2}\int_{\RR}\int_{\RR}\frac{(\phi_i(x)-\phi_i(y))^2}{|x-y|^{1+2s}}\,dxdy.
	\end{align}
\end{itemize}

Let us remark that the integrals in \eqref{eq:Sec2_AdiagonalUpp} and \eqref{eq:Sec2_Adiagonal} are defined on the entire real line and, therefore, they need to be treated with special attention. Those integrals, as well as \eqref{eq:Sec2_Adiagonal}, have been computed exactly in \cite{biccari2019controllability} yielding to explicit values for the elements $a_{i,j}$, which only depend on $i$, $j$, $s$ and $h$. For the sake of brevity, the complete computations are omitted here. 

\begin{remark}\label{rem:Sec2_rem_prel}
It is worth noticing that for our implementation, we used $\mathcal P^1$ elements, meaning that we are considering piece-wise linear and continuous functions for the FE approximation. Actually, for $s<1/2$, we could have used $\mathcal P^0$ (that is, piece-wise constant) elements too. Indeed, in this case, the spaces $H^s_0(-1,1)$ and $H^s(-1,1)$ coincide (see \cite[Chapter 1]{grisvard2011elliptic} and \cite[Chapter 11]{lions1968problemes}), therefore it is not necessary that the discrete functions $v_h$ vanish on the boundary in order to have a conforming method. On the other hand, a $\mathcal P^1$ basis is indeed necessary when $s\geq 1/2$. The case of $\mathcal P^0$ for $s<1/2$ is not addressed here in order to present a unified approach for the whole range $0<s<1$.
\end{remark}

\subsubsection{2-D case}

We now focus on the Dirichlet problem \eqref{eq:Sec2_DP} in dimension $d=2$, for which the FE procedure has been described in \cite{acosta2017short,acosta2017fractional,borthagaray2017laplaciano}. In this case, we have 
\begin{align*}
	& a_{i,j} = \frac{C_{2,s}}{2} \int_{\RR^2}\int_{\RR^2}\frac{(\phi_i(x)-\phi_i(y))(\phi_j(x)-\phi_j(y))}{|x-y|^{2+2s}}\,dxdy.
\end{align*}

Moreover, since $\phi_i=0$ in $\Omega^c$, the integral on $\RR^2\times\RR^2$ is reduced to integrals on the set $(\Omega\times\Omega)\cup(\Omega\times\Omega^c)\cup(\Omega^c\times\Omega)$ and, taking into account that the interactions in $\Omega\times\Omega^c$ and $\Omega^c\times\Omega$ are symmetric with respect to $x$ and $y$, we get
\begin{align}\label{eq:Sec2_Kij}
	a_{i,j} =& \, \frac{C_{2,s}}{2} \int_{\Omega}\int_{\Omega}\frac{(\phi_i(x)-\phi_i(y))(\phi_j(x)-\phi_j(y))}{|x-y|^{2+2s}}\,dxdy \notag
	\\
	&+ C_{2,s} \int_{\Omega}\int_{\Omega^c}\frac{\phi_i(x)\phi_j(x)}{|x-y|^{2+2s}}\,dxdy.
\end{align}

\begin{SCfigure}[1.2][h]
	\centering 
	\includegraphics[width=4.5cm,height=4cm]{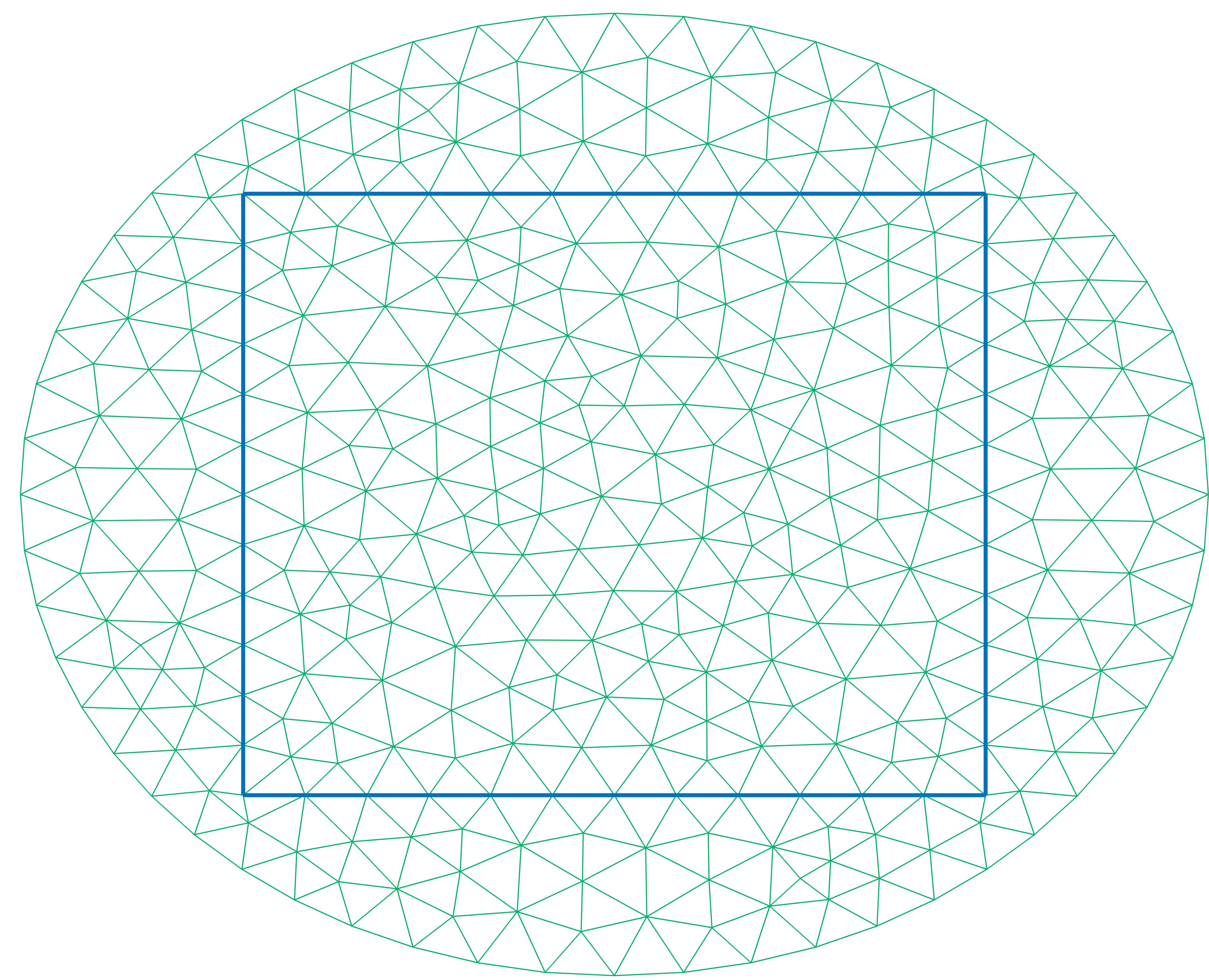}
	\caption{A square domain $\Omega$ (delimited by the blue lines) and an auxiliary ball containing it. Regular triangulations $\mathcal T$ and $\mathcal T_A$ for $\Omega$ and $B\setminus\Omega$ are shown.}\label{fig:Sec2_mesh}
\end{SCfigure}

Notice that the second integral in \eqref{eq:Sec2_Kij} has to be computed over the unbounded domain $\Omega^c$. To do that, it is convenient to introduce a ball $B$ containing $\Omega$ as in Figure \ref{fig:Sec2_mesh}, since this allows to employ polar coordinates and exploit symmetry properties. 

We indicate with $N_B$ the number of elements on the triangulation of $B$. Then, recalling \eqref{eq:Sec2_Kij}, the coefficients $a_{i,j}$ are given by the expression 
\begin{align*}
	a_{i,j} = \frac{C_{2,s}}{2}\sum_{\ell = 1}^{N_B}\left(\sum_{m = 1}^{N_B}  \mathcal I_{\ell,m}^{i,j} + 2\mathcal J_\ell^{i,j}\right),
\end{align*}
where, for $1\leq \ell,m\leq N_B$, we defined 
\begin{subequations}
	\begin{align}
		& \mathcal I_{\ell,m}^{i,j} := \int_{T_\ell}\int_{T_m} \frac{(\phi_i(x)-\phi_i(y))(\phi_j(x)-\phi_j(y))}{|x-y|^{2+2s}}\,dxdy, \label{eq:Sec2_FEint1}
		\\[5pt]
		& \mathcal J_\ell^{i,j}:= \int_{T_\ell}\int_{B^c} \frac{\phi_i(x)\phi_j(x)}{|x-y|^{2+2s}}\,dxdy. \label{eq:Sec2_FEint2}
	\end{align}
\end{subequations}

The computations for each of the above integrals are challenging for different reasons: \eqref{eq:Sec2_FEint1} involves a singular integrand if $\overline{T}_\ell\cap\overline{T}_m\neq\emptyset$, while \eqref{eq:Sec2_FEint2} needs to be calculated on an unbounded domain. A complete discussion on the actual computation of \eqref{eq:Sec2_FEint1} and \eqref{eq:Sec2_FEint2}, together with the corresponding algorithm, can be found in \cite{acosta2017short,borthagaray2017laplaciano}. For the sake of brevity, we omit this discussion here.

\subsubsection{The fractional Laplacian with non-homogeneous exterior condition}\label{subs:Sec2_FLext}

We conclude this section with an abridged discussion on the FE approximation of the fractional Laplacian with non-homogeneous exterior condition. That is, for any $f\in L^2(\Omega)$ and $g\in H^s(\Omega^c)$, we consider the following elliptic problem
\begin{align}\label{eq:Sec2_DPext}
	\begin{cases}
		\fl{s}{u} = f & \mbox{ in }\;\Omega
		\\
		u = g & \mbox{ in }\;\Omega^c.
	\end{cases}
\end{align}
There are several possibilities to discretize \eqref{eq:Sec2_DPext} in FE:

\begin{itemize}
	\item[1.] A first possibility is to start from the variational formulation associated with \eqref{eq:Sec2_DPext}, as we did before for the problem with homogeneous exterior condition. Employing the integration by parts formula in Proposition \ref{prop:App2_prop}, this variational formulation is given by: find $u\in H^s(\mathbb{R}^d)$ such that $u=g$ in $\Omega^c$ and, for all $v\in H^s_0(\Omega)$,
	\begin{align}\label{eq:Sec2_WFext}
		\int_\Omega fv\,dx = a(u,v),
	\end{align}
	where $a(\cdot,\cdot)$ is the bilinear form given in \eqref{eq:Sec2_WF}. Moreover, notice that because of the non-homogeneous exterior datum in \eqref{eq:Sec2_DPext}, this bilinear form reads as follows
	\begin{align}\label{eq:Sec2_bilinear}
		a(u,v) =&\; \frac{C_{d,s}}{2}\int_{\Omega}\int_{\Omega}\frac{(u(x)-u(y))(v(x)-v(y))}{|x-y|^{d+2s}}\;dxdy
		\\
		& + C_{d,s}\int_{\Omega}\int_{\Omega^c}\frac{(u(x)-g(y))(v(x)-v(y))}{|x-y|^{d+2s}}\;dxdy. \notag 
	\end{align}
	One then needs to introduce some appropriate FE space, with the corresponding basis on which to decompose the functions $u$, $f$ and $g$, so that \eqref{eq:Sec2_bilinear} is reduced to solve a linear system. We stress that this may be a delicate issue, due to the fact that these functions are supported on different domains. As a matter of fact, to our knowledge this approach has never been used in the literature.

	\item[2.] A second possibility would be to convert the non-homogeneous exterior problem \eqref{eq:Sec2_DPext} into a homogeneous one. This can be done in the following way: let $G\in H^s(\RR^d)$ such that $G|_{\Omega^c} = g$ and denote $u_G$ the solution to the Dirichlet problem
	\begin{align*}
	\begin{cases}
		\fl{s}{u}_G = f-\fl{s}{G} & \mbox{ in }\;\Omega
		\\
		u_G = 0 & \mbox{ in }\;\Omega^c.
	\end{cases}
	\end{align*}
	Notice that, since $f-\fl{s}{G}\in H^{-s}(\Omega)$, then we have from Proposition \ref{prop:App1_WPprop} that $u_G\in H_0^s(\Omega)$. Moreover, $u_G$ can be approximated using the FE scheme for the homogeneous exterior problem that we have introduced in the previous section. Then, finally, the solution of \eqref{eq:Sec2_DPext} would simply be given by $u=G+u_G$.

	\item[3.] A third possibility is to approximate the exterior Dirichlet problem \eqref{eq:Sec2_DPext} with an exterior Robin one, as it has been done in our recent contributions \cite{antil2019external,antil2020external}. This is the approach that we will use later for our simulations. A more detailed presentation will be given in Section \ref{subs:Sec4_remarks}. 

\end{itemize}

\subsection{Error analysis}\label{subs:Sec2_error}

We devote this section to a brief discussion on the approximation error of the FE discretization we just introduced. In particular, we will present the convergence rates of the FE scheme in the $H^s_0(\Omega)$ and in the $L^2(\Omega)$-norms. The former is the energy norm associated with the elliptic problem \eqref{eq:Sec2_DP}. In this case, we have the following result.

\begin{theorem}[{\cite[Theorem 4.6]{acosta2017fractional}}]\label{thm:Sec2_errorHs}
Let $f$ satisfy the following regularity assumptions:
\begin{displaymath}
	\begin{array}{ll}
		f\in C^{\frac 12 -s}(\Omega), & \textnormal{if } 0<s<1/2,
		\\
		f\in L^\infty(\Omega), & \textnormal{if } s=1/2,
		\\
		f\in C^\beta(\Omega) \textnormal{ for some } \beta>0, & \textnormal{if } 1/2<s<1,
	\end{array}
\end{displaymath}
where $C^{s-\frac 12}(\Omega)$ and $C^{\beta}(\Omega)$ denote the standard H\"older spaces of order $s-1/2$ and $\beta$, respectively. Then, for the solution $u$ of \eqref{eq:Sec2_WF} and its FE approximation $u_h$ on a uniform mesh with size $h$, we have the following a priori estimates 
\begin{align*}
	&\norm{u-u_h}{H^s_0(\Omega)}\leq \frac{C(s,\sigma)}{\varepsilon} h^{\frac 12-\varepsilon}\,\norm{f}{C^{\frac 12-s}(\Omega)},& &\forall\varepsilon>0,\;\;\textnormal{if}\quad s<1/2, 
	\\
	&\norm{u-u_h}{H^s_0(\Omega)}\leq \frac{C(\sigma)}{\varepsilon} h^{\frac 12-\varepsilon}\, \norm{f}{L^\infty(\Omega)},& &\forall\varepsilon>0,\;\;\textnormal{if}\quad  s=1/2 
	\\
	&\norm{u-u_h}{H^s_0(\Omega)}\leq \tfrac{C(s,\beta,\sigma)}{\sqrt{\varepsilon}(2s-1)} h^{\frac 12-\varepsilon}\, \norm{f}{C^\beta(\Omega)},& &\forall\varepsilon>0,\;\;\textnormal{if} \quad s>1/2,\;\;\beta>0
\end{align*}
where $C$ is a positive constant not depending on $h$ and $\sigma$ is given by \eqref{eq:Sec2_mesh2D}. 
\end{theorem}

The $L^2$ convergence rate, instead, will come into play in Section \ref{sec:3} when addressing the implementation of the null controllability problem for the fractional heat equation. 

\begin{theorem}[{\cite[Proposition 3.3.2]{borthagaray2017laplaciano}}]\label{thm:Sec2_errorL2}
Let $\alpha:=\min\{s, 1/2 -\delta\}$, with $\delta>0$ arbitrary small. If $f\in L^2(\Omega)$ and $u$ is the solution to \eqref{eq:Sec2_DP}, for its FE approximation on a uniform mesh with size $h$ it holds that
\begin{align*}
	\norm{u-u_h}{L^2(\Omega)}\leq C(s,\alpha)h^{2\alpha}\norm{f}{L^2(\Omega)}.
\end{align*}
\end{theorem}

\subsection{Numerical experiments}\label{subs:Sec2_simul}

We conclude this section with some numerical experiments. Here we focus on the one-dimensional case on the interval $\Omega = (-1,1)$. Analogous simulations in dimension $d=2$ can be found in \cite{acosta2017fractional}. To test the efficiency of our FE scheme, we consider the Dirichlet problem \eqref{eq:Sec2_DP} with $f=1$, whose explicit solution is given by (see \cite{getoor1961first}) 
\begin{align}\label{eq:Sec2_real_sol}
	u(x)=\gamma_s\Big(1-x^{\,2}\Big)^s\cdot\chi_{(-1,1)}, \quad\quad \gamma_s = \frac{2^{-2s}\sqrt{\pi}}{\Gamma\left(\frac{1+2s}{2}\right)\Gamma(1+s)}.
\end{align}

\begin{figure}[h]
	\centering
	\includegraphics[scale=0.28]{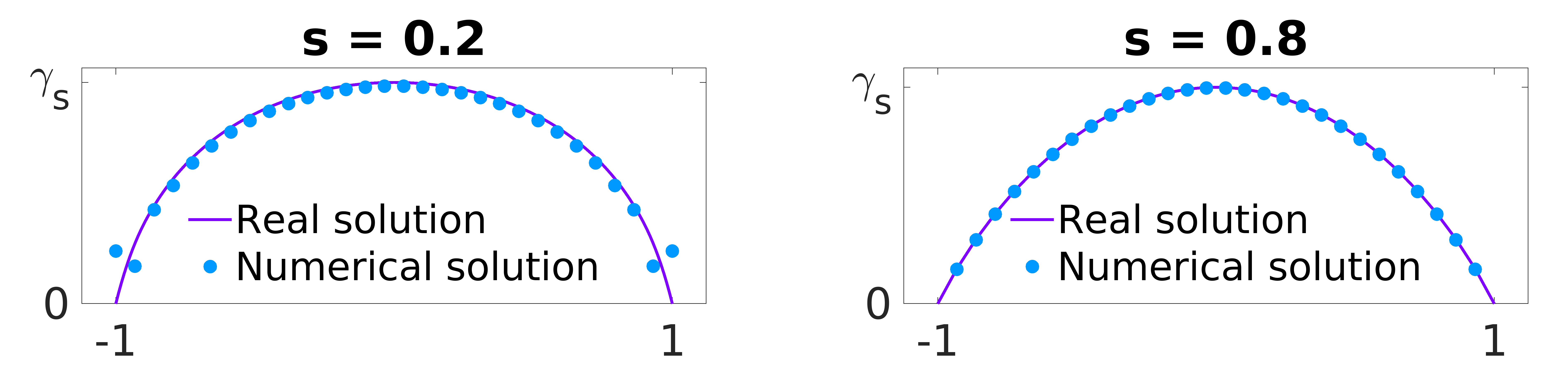}
	\caption{Exact and numerical solution of the Dirichlet problem \eqref{eq:Sec2_DP} with $f=1$ for different values of $s\in (0,1)$.}\label{fig:Sec2_comparison}
\end{figure}

In Figure \ref{fig:Sec2_comparison}, we show a comparison for different values of $s$ between the exact solution \eqref{eq:Sec2_real_sol} and the computed numerical approximation. One can notice that for large $s\geq 1/2$ the FE scheme provides a good approximation. On the other hand, when $s<1/2$, the computed solution is to a certain extent different from the exact one, as there is a discrepancy approaching the boundary. 

\begin{SCfigure}[1.1][h]
	\centering 
	\includegraphics[scale=0.19]{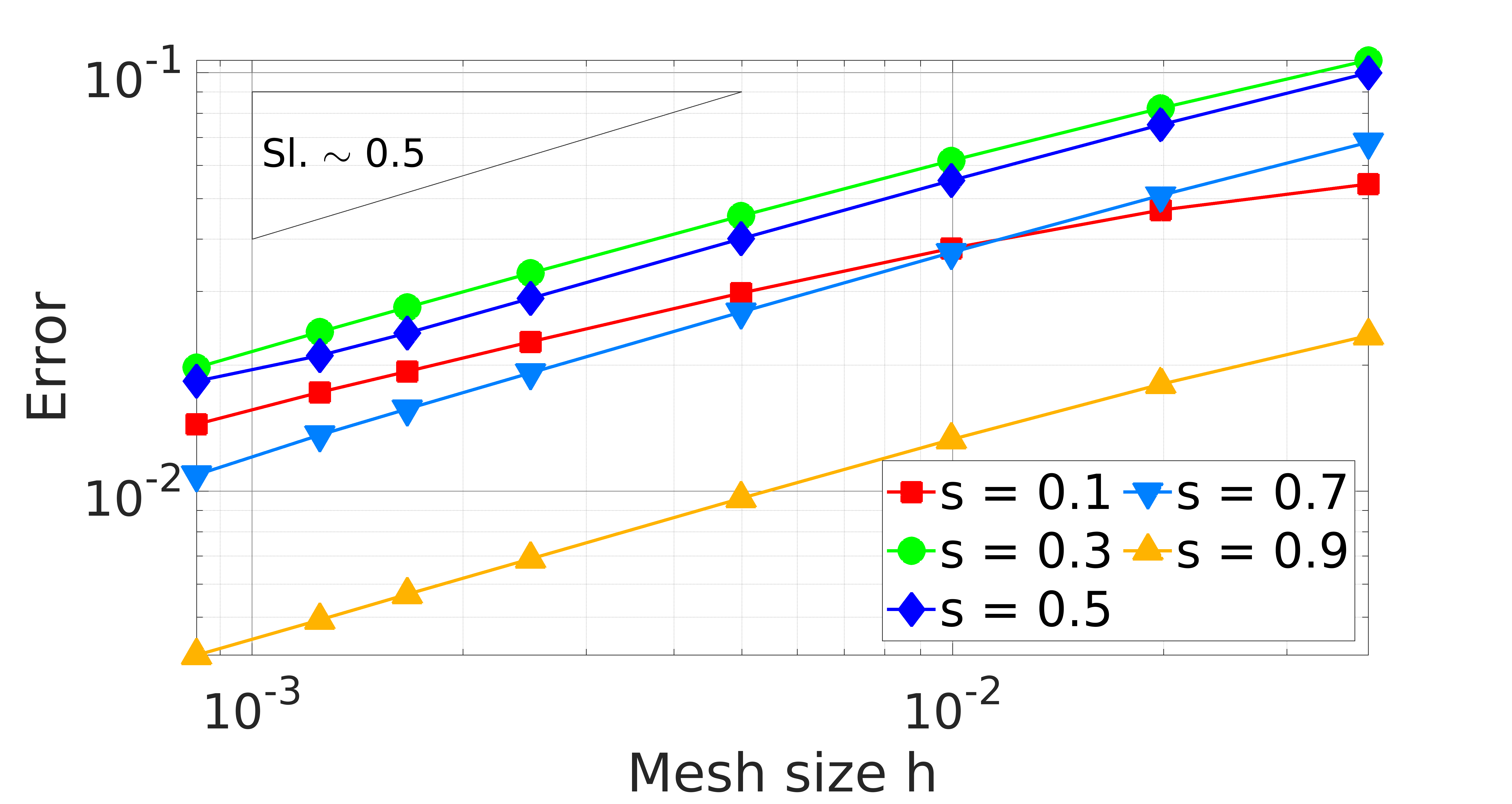}
	\caption{Convergence in $H^s(-1,1)$ of the approximation error for different values of $s\in (0,1)$. The rate $\norm{u-u_h}{H_0^s(-1,1)}\sim\sqrt{h}$ is in accordance with the theoretical result.}\label{fig:Sec2_errorFL}
\end{SCfigure}

Despite this fact, we can see in Figure \ref{fig:Sec2_errorFL} that, for all $s\in (0,1)$ the $H^s$ error of our FE approximations decreases with $h$ at a rate $\norm{u-u_h}{H_0^s(-1,1)}\sim\sqrt{h}$, which is the expected one according to Theorem \ref{thm:Sec2_errorHs}. This is a numerical evidence of the accuracy of the proposed FE scheme when discretizing the elliptic problem \eqref{eq:Sec2_DP}.

\section{Interior controllability properties of the fractional heat equation}\label{sec:3}

In this section, we discuss the interior controllability properties of the fractional heat equation. That is, for any $y_0\in L^2(\Omega)$, $\Omega\subset\mathbb{R}^d$ being a bounded and $C^{1,1}$ domain, and $\omega\subset\Omega$ a nonempty and open subset, we are going to consider the following control system
\begin{align}\label{eq:Sec3_HeatControl}
	\begin{cases}
		y_t + \fl{s}{y} = u\chi_\omega &\mbox{ in }\, \Omega\times(0,T),
		\\
		y\equiv 0 &\mbox{ in } \Omega^c\times(0,T),
		\\
		y(\cdot,0) = y_0 & \mbox{ in }\Omega.
\end{cases}
\end{align}

Here $u$ is the control function belonging to some functional space to be specified later. We will start by recalling in Subsection \ref{subs:Sec3_HeatControlResults} the theoretical controllability results which are presently available in the literature. Secondly, in Subsection \ref{subs:Sec3_PenalizedHUM}, we will give an abridged presentation of the well-known \textit{penalized Hilbert Uniqueness Method} that we shall employ for the computation of our numerical controls. Finally, Subsection \ref{subs:Sec3_numerics} is devoted to the presentation and discussion of our numerical experiments.

\subsection{Review of theoretical controllability results}\label{subs:Sec3_HeatControlResults}

We summarize the theoretical controllability results which can be currently found in the literature for the fractional heat equation \eqref{eq:Sec3_HeatControl}.

\begin{theorem}\label{thm:Sec3_controlUnconstrained}
For the fractional heat equation \eqref{eq:Sec3_HeatControl}, the following results hold.
\begin{itemize}
	\item[1.] \textbf{Approximate controllability}. Let $\omega\subset\Omega$ be any nonempty and open subset of $\Omega$ and $s\in(0,1)$. For any $T>0$, $y_0,y_T\in L^2(\Omega)$ and $\varepsilon>0$, there exists a control function $u\in L^2(\omega\times(0,T))$ such that the unique solution $y$ of \eqref{eq:Sec3_HeatControl} satisfies $\norm{y(\cdot,T)-y_T}{L^2(\Omega)}\leq\varepsilon$.
	
	\item[2.] \textbf{$L^2$ null controllability in 1-D}. Let $\Omega=(-1,1)$ and $\omega\subset(-1,1)$ be any nonempty and open subset. For any $T>0$ and $y_0\in L^2(-1,1)$, there exists a control function $u\in L^2(\omega\times(0,T))$ such that the unique solution $y$ of \eqref{eq:Sec3_HeatControl} satisfies $y(\cdot,T)=0$ a.e. in $(-1,1)$, if and only if $s\in(1/2,1)$.
	
	\item[3.] \textbf{$L^\infty$ null controllability in 1-D}. Let $\Omega=(-1,1)$ and $\omega\subset(-1,1)$ be any nonempty and open subset. For any $T>0$ and $y_0\in L^2(-1,1)$, there exists a control function $u\in L^\infty(\omega\times(0,T))$ such that the unique solution $y$ of \eqref{eq:Sec3_HeatControl} satisfies $y(\cdot,T)=0$ a.e. in $(-1,1)$, if and only if $s\in(1/2,1)$.
	
	\item[4.] \textbf{$L^2$ null controllability in multi-D}. Let $\omega\subset\Omega$ be a neighborhood of $\partial\Omega$. For any $T>0$ and $y_0\in L^2(\Omega)$, there exists a control function $u\in L^2(\omega\times(0,T))$ such that the unique solution $y$ of \eqref{eq:Sec3_HeatControl} satisfies $y(\cdot,T)=0$ a.e. in $\Omega$, if and only if $s\in(1/2,1)$.	
\end{itemize}
\end{theorem}

The approximate controllability property has been established in \cite[Theorem 1.2]{biccari2019controllability}. As for any analytic semi-group (as it is the one generated by $-\fl{s}{}_D$, see Appendix \ref{appendixA} or \cite[Lemma 2.5]{biccari2017local}), this is a consequence of the unique continuation of the fractional Laplacian proved in \cite{fall2014unique}, which yields to the unique continuation of the corresponding parabolic problem. The null controllability in space dimension $d=1$ with $L^2(\omega\times(0,T))$ controls has been proved in \cite[Theorem 1.1]{biccari2019controllability}, while the same result with $L^\infty(\omega\times(0,T))$ controls has been obtained in \cite[Theorem 2.2]{biccari2020positive}. Here we shall highlight the restriction $s\in (1/2,1)$. At this regard, we remark that the exponent $s$ determines the diffusion strength of the fractional heat semi-group. The limitation $s\in (1/2,1)$ then means that, when dealing with heat-like processes involving the fractional Laplacian, a minimal amount of diffusivity is required to obtain positive null-controllability results. On a mathematical perspective, this arises from the application of moments techniques and parabolic Ingham inequalities in our proofs, and from the lack of an asymptotic gap between the eigenvalues of the one-dimensional Dirichlet fractional Laplacian when $s\leq 1/2$. Finally, the null controllability for $d\geq 2$ has been obtained in \cite[Theorem 1.1]{biccari2021null} through a combination of partial observability results for the fractional wave equation involving the fractional Laplacian, transmutation techniques, and the Lebeau-Robbiano strategy. In this multi-dimensional case, apart from the limitation $s\in (1/2,1)$, the employment of transmutation techniques imposes the further restriction that $\omega$ has to be a neighborhood of the boundary $\partial\Omega$.

In addition to the above results, the controllability property under state/control constraints has been investigated for the one-dimensional fractional heat equation in \cite{biccari2020positive}. In particular, the following has been obtained.

\begin{theorem}[1-D constrained controllability, {\cite[Theorem 2.1]{biccari2020positive}}]\label{thm:Sec3_controlConstraints}
Let $s\in(1/2,1)$, $y_0\in L^2(-1,1)$ and let $\widehat{y}$ be a positive trajectory, i.e., a solution of \eqref{eq:Sec3_HeatControl} with initial datum $0<\widehat{y}_0\in L^2(-1,1)$ and right-hand side $\widehat{u}\in L^\infty(\omega\times(0,T))$. Assume that there exists $\varpi> 0$ such that $\widehat{u}\geq\varpi$ a.e. in $\omega\times(0,T)$. Then, the following assertions hold.
\begin{itemize}
	\item[1.] There exists a minimal strictly positive controllability time $T_{\rm min}>0$ such that, for all $T> T_{\rm min}$, we can find a non-negative control $u\in L^\infty(\omega\times(0,T))$ whose corresponding solution $y$ of \eqref{eq:Sec3_HeatControl} satisfies $y(\cdot,T) = \widehat{y}(\cdot,T)$ a.e. in $(-1,1)$. Moreover, if $y_0\geq 0$, then $y(x,t)\geq 0$ for every $(x,t)\in (-1,1)\times (0,T)$. 	
	
	\item[2.] For $T = T_{\rm min}$, there exists a non-negative control $u\in\mathcal M(\omega\times(0,T_{\rm min}))$, the space of Radon measures on $\omega\times(0,T_{\rm min})$, such that the corresponding solution $y$ of \eqref{eq:Sec3_HeatControl} satisfies $y(\cdot,T) = \widehat{y}(\cdot,T)$ a.e. in $(-1,1)$.
\end{itemize}
\end{theorem}

For completeness, we shall mention that the null controllability properties of \textit{spectral fractional heat-like equations} have been also analyzed in \cite{micu2006controllability,miller2006controllability} in the case where the integral fractional Laplace operator $\fl{s}{}$ is replaced by the spectral Dirichlet fractional Laplacian mentioned in Section \ref{sec:2}. Recall that this operator is different from the integral one considered in this work. Nevertheless, when addressing the null control of fractional diffusion, it shows an analogous behavior as $\fl{s}{}$. In particular, also for the spectral operator null controllability requires a minimum amount of diffusivity, $s=1/2$ being the critical case, in which the property fails. Consequently, also the spectral fractional heat equation turns out to be null controllable if and only if $1/2<s<1$. This spectral fractional Laplacian is easier to handle and well understood since its eigenfunctions are the same as those of the Laplacian with the Dirichlet boundary condition and the associated eigenvalues are the $s$-power of the eigenvalues of the Laplacian with the Dirichlet boundary condition. Thus, for the spectral fractional Laplacian, the existing observability results for the classical heat equation and its eigenfunctions can be reused to show that null controllability is preserved in the range $1/2<s<1$. By the contrary, for the integral fractional Laplacian, the eigenfunctions depend genuinely on $s$ and the same techniques inherited from the Laplacian do not apply in general. At this regard, let us stress that for the Laplace operator, many fundamental properties are based on Carleman inequalities that have not been properly developed for the fractional setting.

\subsection{The penalized Hilbert Uniqueness Method}\label{subs:Sec3_PenalizedHUM}

We devote this section to an abridged description of the \textit{penalized Hilbert Uniqueness Method} (HUM) that we shall employ for computing the controls for the fractional heat equation. Here we will mostly refer to the works \cite{boyer2013penalised,glowinski1994exact,glowinski2008exact}. 

Let $\mathcal H$ be a Hilbert space and $\mathcal A: \mathcal D(\mathcal A)\subset \mathcal H\to \mathcal H$ be an unbounded operator such that $-\mathcal A$ generates an analytic semi-group. Let $\mathcal U$ be another Hilbert space and $\mathcal B: \mathcal U\to \mathcal D(\mathcal A)^\star$ be a bounded operator. Let $T>0$ be given and, for any $y_0\in \mathcal H$ and $u\in L^2(0,T;\mathcal U)$, let us consider the Cauchy problem
\begin{align}\label{eq:Sec3_abstractPB}
	y_t + \mathcal Ay = \mathcal Bu \;\mbox{ in } [0,T], \quad y(0) = y_0.
\end{align}

The penalized HUM approach consists in finding the control of minimal $L^2(0,T;\mathcal U)$ norm for \eqref{eq:Sec3_abstractPB} by means of the following optimization problem: 
\begin{align}\label{eq:Sec3_HUMprimal}
	& u_\beta = \min_{u\in L^2(0,T;\mathcal U)} F_\beta (u)
	\\
	& F_\beta(u):=\frac{1}{2}\int_0^T\norm{u(t)}{\mathcal U}^2\,dt+\frac{1}{2\beta}\norm{y(T)}{\mathcal H}^2.\notag 
\end{align}

Notice that, if \eqref{eq:Sec3_abstractPB} is controllable (either null, exactly or approximately), then for any $\beta > 0$, the functional $F_\beta$ is strictly convex, continuous and coercive. Hence, it has a unique minimizer $u_\beta\in L^2(0,T;\mathcal U)$. Moreover, we have the following result.

\begin{theorem}[{\cite[Theorem 1.7]{boyer2013penalised}}]\label{thm:Sec3_HUMthm}
The following controllability properties hold.	
\begin{itemize}
	\item[1.] Problem \eqref{eq:Sec3_abstractPB} is approximately controllable at time $T>0$ from $y_0\in \mathcal H$ if and only if $\norm{y_\beta (T)}{\mathcal H}\rightarrow 0$ as $\beta\to 0$, where $y_\beta$ denotes the solution corresponding to $u_\beta$.
	\item[2.] Problem \eqref{eq:Sec3_abstractPB} is null-controllable at time $T>0$ from $y_0\in \mathcal H$ if and only if
	\begin{align}\label{eq:Sec3_bound}
		\mathcal E_{y_0}:=2\sup_{\beta>0}\left(\inf_{u\in L^2(0,T;\mathcal U)}F_\beta(u)\right)<+\infty.
	\end{align}
	In this case, we have 
	\begin{subequations}
		\begin{align}
			&\norm{u_\beta}{L^2(0,T;\mathcal U)}\leq \sqrt{\mathcal E_{y_0}}\label{eq:Sec3_functBehavior1}
			\\
			&\norm{y_\beta(T)}{\mathcal H}\leq \sqrt{\mathcal E_{y_0}\beta}.\label{eq:Sec3_functBehavior2}
		\end{align}	
	\end{subequations}
	Moreover, as $\beta\to 0$, $u_\beta\to \bar{u}$ strongly in $L^2(0,T;\mathcal U)$, $\bar{u}$ being the optimal control obtained from the functional \eqref{eq:Sec3_HUMprimal} without the second penalization term.
\end{itemize}
\end{theorem}

According to Theorem \ref{thm:Sec3_HUMthm}, there is then an essential difference between approximate and null controllability in this penalization context. In both cases, the solution at time $T$ of \eqref{eq:Sec3_abstractPB} corresponding to $u_\beta$ converges to zero in $\mathcal H$ as $\beta\to 0$. Nevertheless, for null-controllability, this convergence has a precise rate $\sqrt{\beta}$ which, as we will see in Section \ref{subs:Sec3_numerics} and had already been observed in \cite{boyer2013penalised}, is typically violated when only approximate controllability holds. Furthermore, when the problem is null controllable, we also have that the control cost $\norm{u_\beta}{L^2(0,T;\mathcal U)}$ remains uniformly bounded which, together with \eqref{eq:Sec3_functBehavior2}, yields the convergence of $u_\beta$ to the solution $\bar{u}$ of the non-penalized problem.

Let us recall that $-\fl{s}{}_D$ generates an analytic semi-group (see \cite{biccari2017local,claus2020realization}). Hence, Theorem \ref{thm:Sec3_HUMthm} applies to \eqref{eq:Sec3_HeatControl} when selecting $\mathcal H = L^2(\Omega)$ and $\mathcal U = L^2(\omega)$. 

Furthermore, to solve the optimization problem \eqref{eq:Sec3_HUMprimal}, we will apply duality and the Fenchel-Rockafellar theory (see e.g. \cite[Chapters VI to VII]{ekeland1999convex}) to build an equivalent dual minimization problem defined on the space $L^2(\Omega)$. Although this duality argument is nowadays classical, for completeness we shall recall below its main steps.

\vspace{0.1cm}
\noindent\textbf{Step 1.} Since the problem \eqref{eq:Sec3_HeatControl} is linear, we can write its solution as $y = \xi + z$, with
\begin{align}\label{eq:Sec3_heatXi}
	\begin{cases}
		\xi_t + \fl{s}{\xi} = 0 &\mbox{ in }\, \Omega\times(0,T),
		\\
		\xi\equiv 0 &\mbox{ in } \Omega^c\times(0,T),
		\\
		\xi(\cdot,0) = y_0 & \mbox{ in }\Omega
	\end{cases}
\end{align} 
and
\begin{align}\label{eq:Sec3_heatZ}
	\begin{cases}
		z_t + \fl{s}{z} = u\chi_\omega &\mbox{ in }\, \Omega\times(0,T),
		\\
		z\equiv 0 &\mbox{ in } \Omega^c\times(0,T),
		\\
		z(\cdot,0) = 0 & \mbox{ in }\Omega.
	\end{cases}
\end{align} 
Moreover, in what follows, for the solution of \eqref{eq:Sec3_heatXi} we will use the notation 
\begin{align*}
	\xi(\cdot,t)= e^{-(-\Delta)_D^st}y_0,
\end{align*}
where $(-\Delta)^s_D$ is the realization in $L^2(\Omega)$ of the fractional Laplacian with the zero Dirichlet exterior condition (see \eqref{eq:App1_FLrealization}).

\noindent\textbf{Step 2.} Let $\mathcal L_T: L^2(\omega\times(0,T)) \rightarrow L^2(\Omega)$ be the linear continuous operator defined as $\mathcal L_T(u) = z(\cdot,T)$, with $z$ solution of \eqref{eq:Sec3_heatZ}. Then, the adjoint operator $\mathcal L_T^\star: L^2(\Omega)\rightarrow L^2(\omega\times(0,T))$ is given by $\mathcal L_T^\star(p_T) = p\chi_\omega$ where, for all $p_T\in L^2(\Omega)$, $p$ solves
\begin{align}\label{eq:Sec3_heatP}
	\begin{cases}
		-p_t + \fl{s}{p} = 0 &\mbox{ in }\, \Omega\times(0,T),
		\\
		p\equiv 0 &\mbox{ in } \Omega^c\times(0,T),
		\\
		p(\cdot,T) = p_T & \mbox{ in }\Omega.
	\end{cases}
\end{align} 

\noindent\textbf{Step 3.} With this notation, we have $F_\beta(u) = \widehat F(u) + G_\beta(\mathcal L_T u)$, where
\begin{align*}
	\widehat F(u):=\frac{1}{2}\int_0^T\norm{u(t)}{L^2(\omega)}^2\,dt \quad\mbox{ and }\quad G_\beta(\mathcal L_T u) :=\frac{1}{2\beta}\norm{\mathcal L_T u + e^{-(-\Delta)_D^sT}y_0}{L^2(\Omega)}^2.
\end{align*}
Since both $\widehat F$ and $G_\beta$ are convex functionals, Fenchel-Rockafellar theory (see \cite[Proposition 1.5]{boyer2013penalised}) yields that
\begin{align}\label{eq:Sec3_FR}
	u_\beta = p_\beta \chi_\omega,
\end{align}
with $p_\beta$ solution of \eqref{eq:Sec3_heatP} corresponding to the initial datum
\begin{align*}
	p_{T,\beta} = \min_{p_T\in L^2(\Omega)} J(p_T),
\end{align*}
and $J(p_T):= \widehat F^\star (\mathcal L_T^\star p_T) + G_\beta^\star(-p_T)$,	$\widehat F^\star$ and $G_\beta^\star$ being the convex conjugates
\begin{equation}\label{eq:Sec3_conjugates}
	\begin{array}{ll}
		\displaystyle\widehat F^\star(u) = \sup_{v\in L^2(\omega\times(0,T))} \left\{\langle u,v\rangle_{L^2(\omega\times(0,T))} - \widehat F(v)\right\}, &  u\in L^2(\omega\times(0,T))
		\\[15pt]
		\displaystyle G_\beta^\star(-p_T) = \sup_{q_T\in L^2(\Omega)} \Big\{-\langle p_T,q_T\rangle_{L^2(\Omega)} - G_\beta(q_T)\Big\}, & q_T\in L^2(\Omega).
	\end{array} 
\end{equation} 

\noindent\textbf{Step 4.} It can be readily checked using \eqref{eq:Sec3_conjugates} that 
\begin{align*}
	&\widehat F^\star(\mathcal L_T^\star p_T)= \frac 12\int_0^T \norm{p(t)}{L^2(\omega)}^2\,dt
	\\
	&G_\beta^\star(-p_T) = \langle p_T,e^{-(-\Delta)_D^sT}y_0\rangle_{L^2(\Omega)} + \frac \beta2 \norm{p_T}{L^2(\Omega)}^2.
\end{align*}
Collecting everything, we then obtain that $J_\beta(p_T)$ is given by  
\begin{align}\label{eq:Sec3_Jfunct}
	J_\beta(p_T)=\frac{1}{2}\int_0^T\norm{p(t)}{L^2(\omega)}^2\,dt + \frac \beta2 \norm{p_T}{L^2(\Omega)}^2 + \langle p_T,e^{-(-\Delta)_D^sT}y_0\rangle_{L^2(\Omega)}.
\end{align}

To compute the numerical control, let us introduce the fully-discrete version of \eqref{eq:Sec3_HeatControl}. Given a uniform mesh $\mathfrak M$ of size $h$ on $\Omega$ and any integer $M>0$, we set $\delta t=T/M$ and we approximate \eqref{eq:Sec3_HeatControl} through an implicit Euler method:
\begin{align}\label{eq:Sec3_HeatControlEuler}
	\begin{cases}
		\displaystyle\mathcal M_h \frac{y_h^{m+1}-y_h^m}{\delta t}+\mathcal A_h y_h^{m+1} = \mathcal B_hu_h^{m+1}, \quad \mbox{ for all } m\in \left\{1,\ldots,M-1\right\}	
		\\
		y_h^1=y_{0,h}, 
\end{cases}
\end{align}
where $y_{0,h}\in \mathbb R^N$ is the projection of the initial datum $y_0\in L^2(\Omega)$ on the mesh $\mathfrak M$, $\mathcal A_h$ and $\mathcal M_h$ are the stiffness and mass matrices given in Section \ref{sec:2}, while the matrix $\mathcal B_h$ has entries (see \cite[Section 1.4.2]{boyer2011uniform})
\begin{align*}
	b_{i,j} = \int_\omega \phi_i(x)\phi_j(x)\,dx, \quad i,j = 1,\ldots, N.
\end{align*}

In \eqref{eq:Sec3_HeatControlEuler}, $u_h=(u_h^m)_{m=1}^M\in\mathbb{R}^{N\times M}$ is a fully-discrete control function, whose cost is given by the discrete $L^2(\Omega\times(0,T))$-norm defined by
\begin{align}\label{eq:Sec3_costDiscr}
	\norm{u_h}{L^2_{h,\delta t}}:=\left(\sum_{m=1}^M\delta t |u^m_h|^2_{L^2_{h,\mathcal M_h}}\right)^{1/2},
\end{align}
and where $|\cdot |_{L^2_{h,\mathcal M_h}}$ is the norm associated with the $L^2$-inner product on $\mathfrak M$ and the mass matrix $\mathcal M_h$:
\begin{align*}
	&\mbox{for all } v = (v_i)_{i=1}^N\in\mathbb{R}^N \quad\mbox{ and }\quad w = (w_i)_{i=1}^N\in\mathbb{R}^N
	\\
	&\langle v,w\rangle_{L^2_{h,\mathcal M_h}} = \langle \mathcal M_h v,w\rangle_{L^2_h} = h \sum_{i=1}^N (\mathcal M_h v)_i w_i \quad\longrightarrow\quad |v|^2_{L^2_{h,\mathcal M_h}} = \langle v,v\rangle_{L^2_{h,\mathcal M_h}}.
\end{align*}

With the above notation, given some penalization parameter $\beta>0$ we can introduce the fully-discrete primal and dual functionals 
\begin{align}\label{eq:Sec3_functDiscr}
	&F_{\beta,h}(u_h)=\frac 12\norm{u_h}{L^2_{h,\delta t}}^2+\frac{1}{2\beta}|y_h^M|_{L^2_{h,\mathcal M_h}}^2 \notag 
	\\
	&J_{\beta,h}(p_h^M) =\,\frac 12\norm{\mathcal B_h p_h}{L^2_{h,\delta t}}^2+\frac{\beta}{2}|p_h^M|^2_{L^2_{h,\mathcal M_h}} +\left\langle p_h^M,e^{\mathcal A_hT}y_{0,h}\right\rangle_{L^2_{h,\mathcal M_h}}, 
\end{align}
with $p_h =(p_h^n)_{n=1}^M\mathbb{R}^{N\times M}$ solution to the adjoint system
\begin{align}\label{eq:Sec3_HeatDualEuler}
	\begin{cases}
		\displaystyle\mathcal M_h \frac{p_h^m-p_h^{m+1}}{\delta t}+\mathcal A_h p_h^m=0, & \mbox{ for all } m\in\{1,\ldots,M-1\}
		\\
		p_h^M = p_{T,h},
	\end{cases}
\end{align}
where $p_{T,h}\in \mathbb R^N$ is the projection of $p_T\in L^2(\Omega)$ on the mesh $\mathfrak M$.

\subsection{Numerical experiments}\label{subs:Sec3_numerics}

To address the actual computation of the fully-discrete controls for \eqref{eq:Sec3_HeatControlEuler}, we apply an optimization algorithm to the dual functional $J_{\beta,h}(p_h^M)$. This functional being quadratic and coercive, the conjugate gradient (CG) method is a natural choice. At this regard, we recall (see \cite{glowinski1994exact,glowinski2008exact}) that the implementation of the CG algorithm requires the gradient of $J_{\beta,h}(p_h^M)$, which is given by $\nabla J_{\beta,h}(p_h^M) = y_h^M + \beta p_h^M$, where $y_h^M$ is the solution at time $T$ of \eqref{eq:Sec3_HeatControlEuler} corresponding to the initial datum $y_{0,h}$ and the control $u_h = \mathcal B_h p_h$. Hence, to compute $\nabla J_{\beta,h}(p_h^M)$ requires to solve two parabolic equations, one forward and the other backward in time. 

In what follows, we will consider the 1-D fractional heat equation on the interval $(-1,1)$ and the 2-D one on the unit ball $B_0(1)$. In both cases, we will discuss the control properties of \eqref{eq:Sec3_HeatControlEuler} from the viewpoint of Theorem \ref{thm:Sec3_HUMthm}. In the context of the fully-discrete problem \eqref{eq:Sec3_HeatControlEuler} this may be a delicate issue. As a matter of fact, it is well known that, in general, we cannot expect for a given bounded family of initial data that the fully-discrete controls are uniformly bounded when $h$, $\delta t$ and $\beta$ tend to zero independently. Instead, we expect to obtain uniform bounds by taking $\beta=\phi(h)$ that tends to zero in connection with the mesh size not too fast (see \cite{boyer2013penalised}) and a time step $\delta t$ verifying some weak condition of the kind $\delta t\leq \zeta(h)$ where $\zeta$ tends to zero logarithmically when $h\to 0$ (see \cite{boyer2011uniform}). In particular, it is crucial to choose properly the penalization parameter $\beta$. Following the discussion in \cite{boyer2013penalised}, a reasonable practical rule is to take $\beta=\phi(h)\sim h^{2p}$, where $p$ is the order of accuracy in space of the numerical method employed for the discretization of the fractional Laplacian. We refer to \cite[Section 4.2]{boyer2013penalised} for some heuristic comments on the choice of this penalization parameter. Nevertheless, we have to stress that, as far as the authors know, the effectiveness of this choice has only been demonstrated numerically and is still not supported by rigorous mathematical results. 

To select the correct value for $p$, let us recall that the solution $y$ to \eqref{eq:Sec3_HeatControl} with $y_0\in L^2(\Omega)$ and $u\in L^2(\omega\times(0,T))$ belongs to $L^2(0,T;H_0^s(\Omega))\cap C([0,T];L^2(\Omega))$. In particular, we have that $y(\cdot,T)\in L^2(\Omega)$. Therefore, we shall choose the value of $p$ as the convergence rate in the $L^2$-norm for the approximation of the elliptic problem \eqref{eq:Sec2_DP}. By virtue of Theorem \ref{thm:Sec2_errorL2}, the appropriate value of $p$ that we shall employ is thus
\begin{align}\label{eq:Sec3_convergenceHeat}
	p = \begin{cases}
		2s, & \textrm{ for }s<\frac 12
		\\
		1-2\delta, & \textrm{ for }s\geq \frac 12
	\end{cases} \quad\longrightarrow\quad \beta = h^{2p} = \begin{cases}
		h^{4s}, & \textrm{ for }s<\frac 12
		\\
		h^{2-4\delta}, & \textrm{ for }s\geq \frac 12,
	\end{cases}
\end{align}
with $\delta>0$ arbitrary small. With this choice of $\beta$, by observing the behavior of the norm of the control, the optimal energy \eqref{eq:Sec3_bound} and the norm of the solution at time $T$, we will obtain numerical evidences for the null and approximate controllability of \eqref{eq:Sec3_HeatControl}, which are in accordance with the theoretical results we recalled in Section \ref{subs:Sec3_HeatControlResults}. Nevertheless, we shall mention that these numerical evidences are, at present time, not supported by a rigorous mathematical analysis, since the controllability of fully-discrete fractional heat equations is still an open question. We refer to Section \ref{sec:6} for more details. 

\subsubsection{1-D simulations}

\begin{SCfigure}[1.1][h]
	\centering
	\includegraphics[scale=0.25]{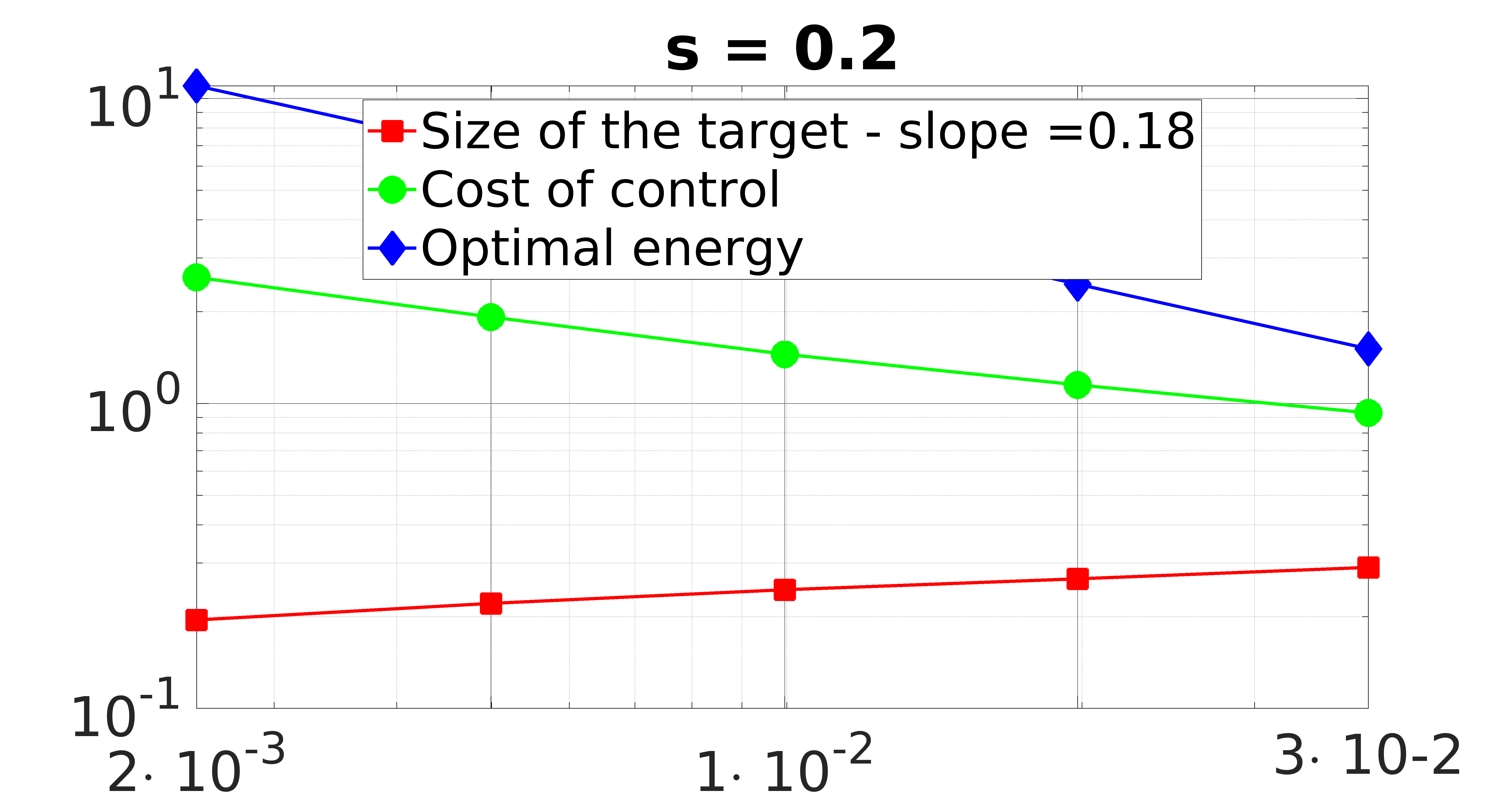}
	\caption{Behavior with respect to the mesh size $h$ of the cost of control, optimal energy and size of the solution to \eqref{eq:Sec3_HeatControlEuler} in space dimension $d=1$ at time $T$ when $s=0.2$.}\label{fig:Sec3_convergence02}
\end{SCfigure}

Let us present here our numerical simulations in one space dimension. To this end, we introduce a uniform $N$-points mesh discretizing the space domain $(-1,1)$. The time interval $(0,T)$ is discretized with a uniform partition as well, this time composed by $M$ points, on which we will implement an implicit Euler method. We begin by considering $s=0.2$, for which we know that \eqref{eq:Sec3_HeatControl} is only approximately controllable. We set $\omega=(-0.3,0.8)$, $T=0.3$ and $y_0(x) = \sin(\pi x)$. We then run the CG algorithm to compute the optimal control $\hat{u}_h$ on several uniform meshes with decreasing mesh-size $h\to 0$, and display in Figure \ref{fig:Sec3_convergence02} the following quantities of interest in dependence of $h$:
\begin{itemize}
	\item The cost of control $\|\hat{u}_h\|_{L^2_{h,\delta t}}$ (see \eqref{eq:Sec3_costDiscr}).
	\item The optimal energy $F_{\beta,h}(\hat{u}_h)$.
	\item The size $|\hat{y}^M|_{L^2_{h,\mathcal M_h}}$ of the corresponding solution at time $T$. 
\end{itemize}

We observe that the discrete $L^2$ norm of $\hat{y}^M$ tends to zero as $h\to 0$, confirming computationally the approximate controllability of \eqref{eq:Sec3_HeatControl}. However, we also see that the cost of the control and the optimal energy increase as $h\to 0$. Hence, numerical evidence indicates that the null controllability is not fulfilled, in agreement with Theorem \ref{thm:Sec3_HUMthm}.

\begin{SCfigure}[1.1][h]
	\centering
	\includegraphics[scale=0.25]{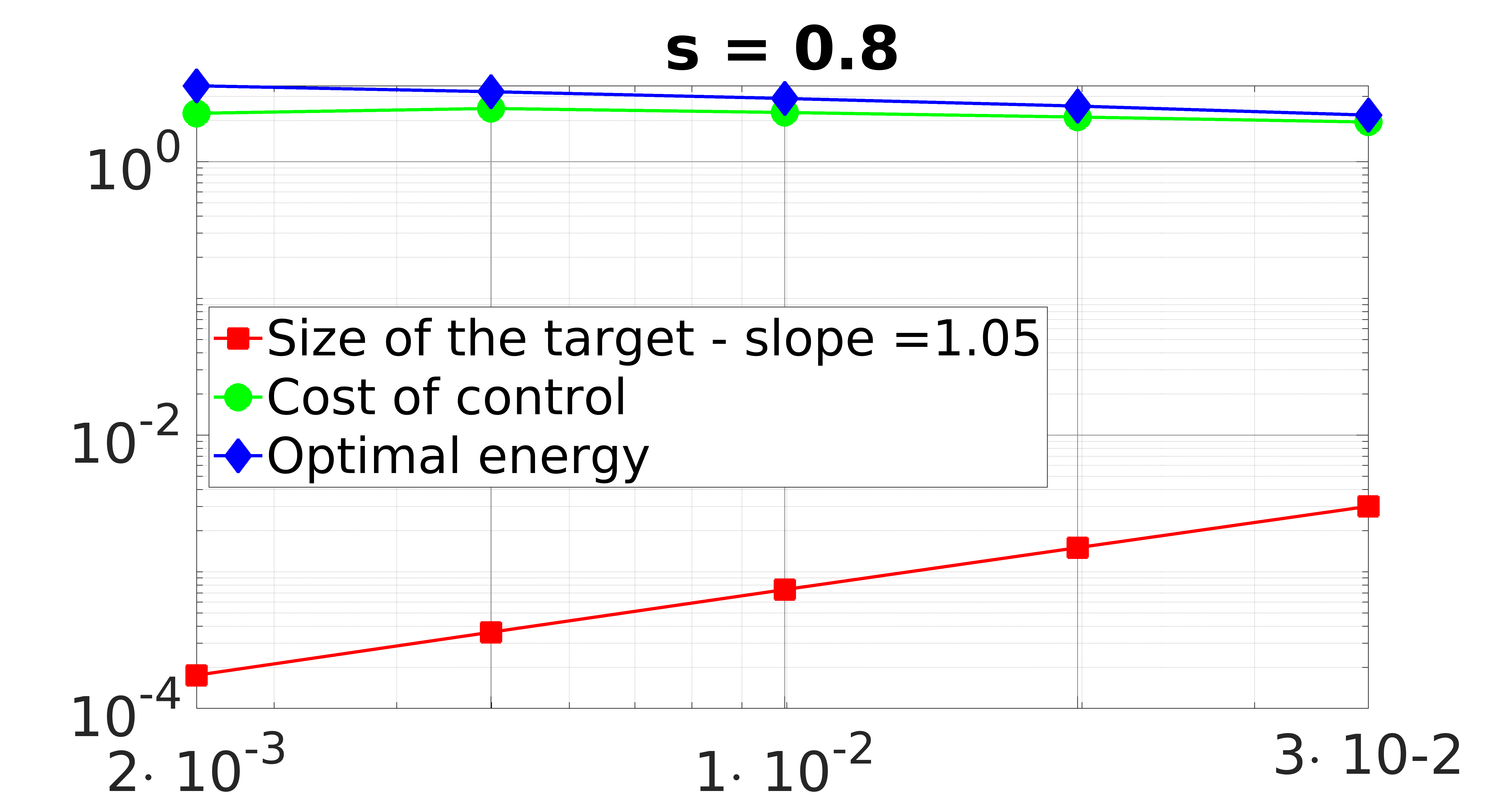}
	\caption{Behavior with respect to the mesh size $h$ of the cost of control, optimal energy and size of the solution to \eqref{eq:Sec3_HeatControlEuler} in space dimension $d=1$ at time $T$ when $s=0.8$.}\label{fig:Sec3_convergence08}
\end{SCfigure}

Let us now take $s=0.8$. In this case, as we can see in Figure \ref{fig:Sec3_convergence08}, the situation changes. Indeed, we can observe that this time the control cost and the optimal energy remain bounded as $h\to 0$. Furthermore, we also see that 
\begin{align*}
	|\hat{y}^M|_{L^2_{h,\mathcal M_h}}\,\sim\, h = \sqrt{\beta},
\end{align*}
which is the expected convergence rate given in \eqref{eq:Sec3_convergenceHeat} for the discrete $L^2$ norm of $y(\cdot,T)$. According to Theorem \ref{thm:Sec3_HUMthm}, all these facts confirm that, for $s=0.8$, \eqref{eq:Sec3_HeatControl} is indeed null controllable as we already proved in Theorem \ref{thm:Sec3_controlUnconstrained}. 

\begin{figure}[ht]
	\centering
	\includegraphics[scale=0.25]{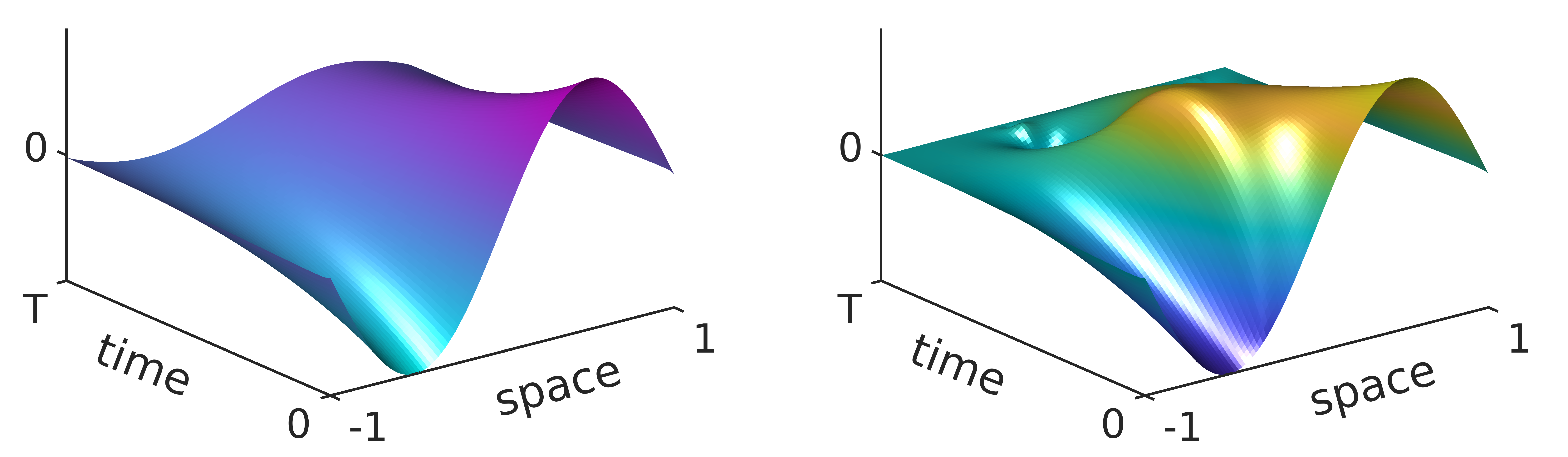}
	\caption{Free (left) and controlled (right) solution of \eqref{eq:Sec3_HeatControl} with $s=0.8$ at time $T=0.3$.}\label{fig:Sec3_solHeat}
\end{figure}

This positive controllability result is also appreciated in Figure \ref{fig:Sec3_solHeat}, where we illustrate the time evolution of the uncontrolled solution as well as the controlled one. We can clearly see that the uncontrolled solution diffuses under the action of the fractional heat semi-group, but does not reach zero at time $T$. On the other hand, the introduction of a control modifies the dynamical behavior of $y$ in such a way that we achieve $y(\cdot,T)=0$. 

\begin{SCfigure}[1][h]
	\centering
	\includegraphics[scale=0.25]{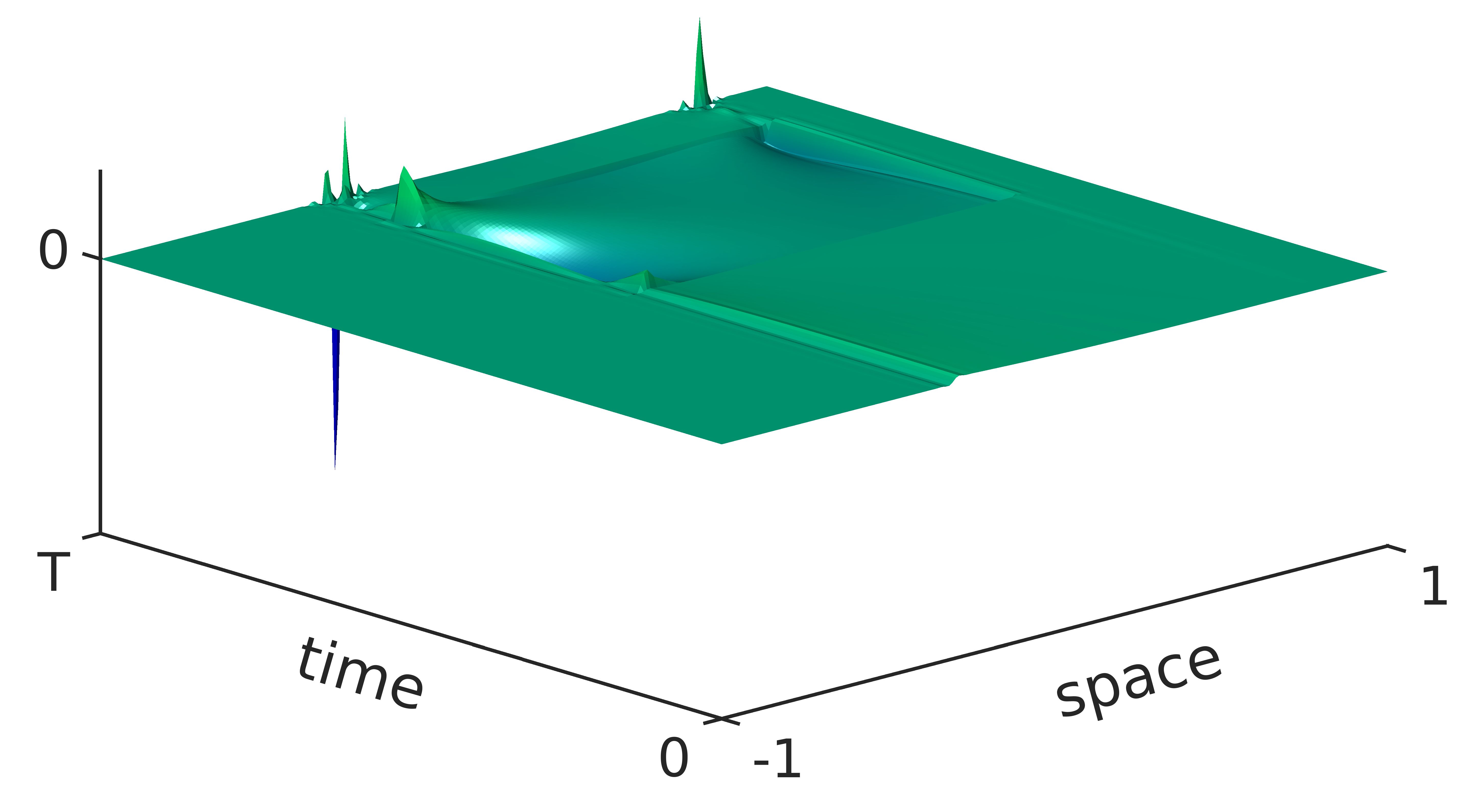}
	\caption{Behavior of the control for the fractional heat equation \eqref{eq:Sec3_HeatControl} with $s=0.8$.}\label{fig:Sec3_controlHeat}
\end{SCfigure}

Finally, in Figure \ref{fig:Sec3_controlHeat}, we see the evolution of the control function which, as we can observe, is almost inactive for a large part of the time interval, and then experiences large oscillations in the proximity of the final time. This fact, related with the characterization of the control as restrictions of solutions of the adjoint system, is in accordance with the \textit{lazy} behavior (observed by Glowinski and Lions in \cite{glowinski1994exact}) of controls for the local heat equation which, at the very beginning, leave the solution evolve under the dissipative effect of the heat semi-group and, only when approaching the final controllability time, inject energy into the system in order to match the desired configuration. 

For completeness, we shall mention that it is by now well-known (see \cite{porretta2013long}) that when the optimal control is chosen differently (in particular, including a \textit{track term} in the cost functional) then the controls are better behaved and in particular exhibit the \textit{turnpike} phenomenon. Nevertheless, this issue is beyond the interest of the present chapter and we will not discuss it further.

\subsubsection{The 2-D case}

Let us now consider the fractional heat equation \eqref{eq:Sec3_HeatControl} in space dimension $d=2$. In this case, we know that null controllability holds in any time $T>0$, acting from a neighborhood of the boundary of the domain $\Omega$, once again if and only if $s\in (1/2,1)$. On the other hand, if $s\in (0,1/2]$, we only have approximate controllability. 

\begin{SCfigure}[1.1][h]
	\centering
	\includegraphics[scale=0.25]{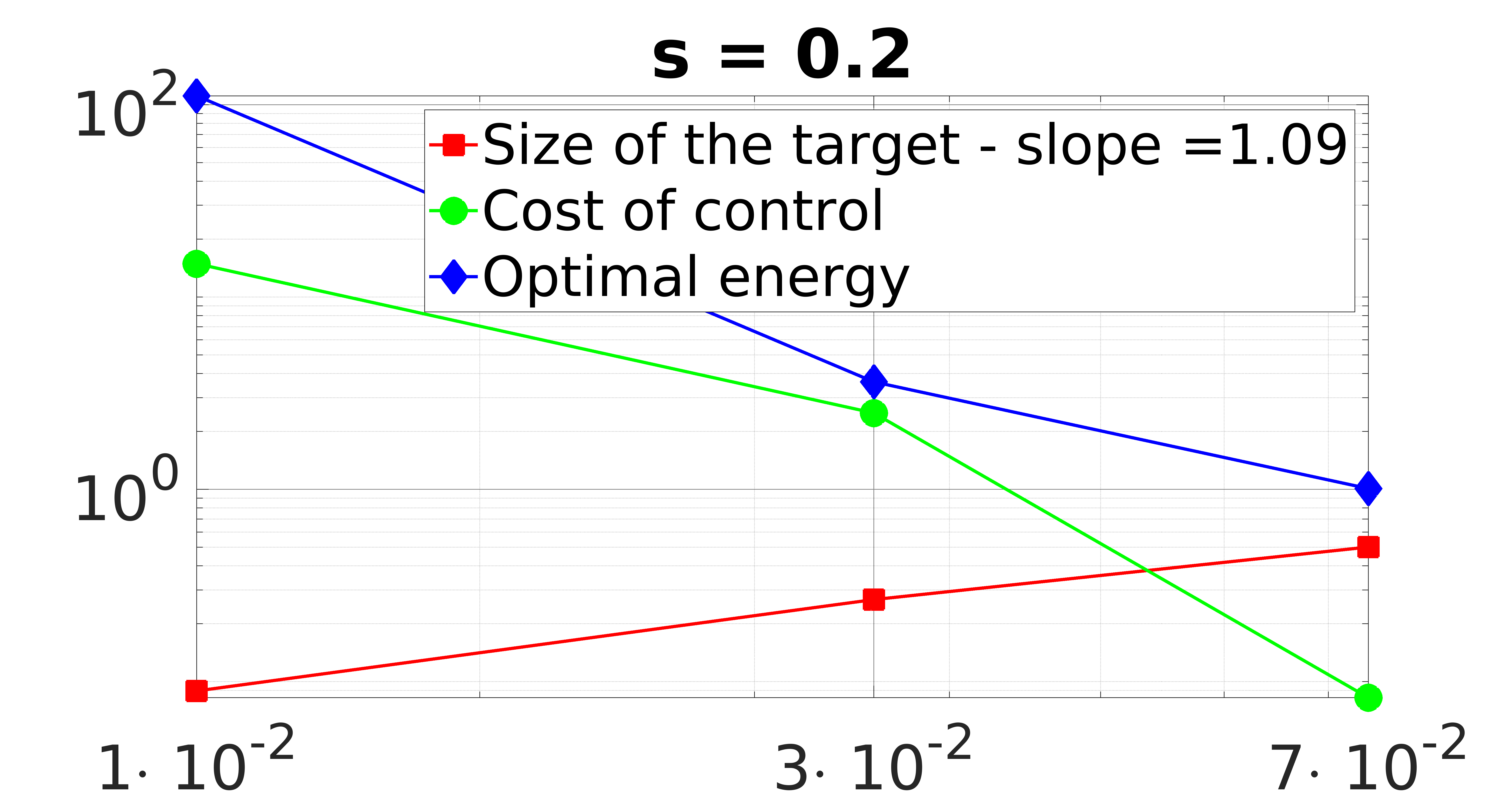}
	\caption{Behavior with respect to the mesh size $h$ of the cost of control, optimal energy and size of the solution to \eqref{eq:Sec3_HeatControlEuler} in space dimension $d=2$ at time $T$ when $s=0.2$.}\label{fig:Sec3_convergence2D02}
\end{SCfigure}

For our simulations, we chose the domain $\Omega=B_0(1)$ as the unit ball centered at zero, while the control region is the ring $\omega = B_0(1)\setminus B_0(4/5)$. Moreover, we considered $y_0(x,y) = \sin(\pi x)\sin(\pi y)$ and the time horizon $T = 0.3$. The control for \eqref{eq:Sec3_HeatControl} has been computed once again by applying the penalized HUM of Section \ref{subs:Sec3_PenalizedHUM}. Nevertheless, now we are considering a two dimensional dynamics described by a full stiffness matrix, whose associated optimal control problem may be computationally very heavy. To alleviate this computational effort and efficiently minimize the discrete functionals \eqref{eq:Sec3_functDiscr}, we have relied on the expert interior-point optimization routine \texttt{IpOpt} (\cite{wachter2006implementation}), with the help of \texttt{CasADi} open-source tool for algorithmic differentiation (\cite{andersson2019casadi}). We start by considering $s=0.2$ and displaying in Figure \ref{fig:Sec3_convergence2D02} the computed values as $h\to 0$ of the cost of control, the optimal energy and the discrete $L^2$ norm of $y(\cdot,T)$. We can see that:
\begin{itemize}
	\item The discrete $L^2$ norm of the solution at time $T$ decreases with $h$, thus confirming the approximate controllability of \eqref{eq:Sec3_HeatControl} according to the first part of Theorem \ref{thm:Sec3_HUMthm}. 
	\item The control cost and the optimal energy both increase as $h\to 0$, thus violating \eqref{eq:Sec3_bound}. Hence, Theorem \ref{thm:Sec3_HUMthm} yields the failure of the null controllability property.
\end{itemize}

Finally, let us consider the case $s=0.8$. This time, we can appreciate how the control cost and the optimal energy remain bounded as $h\to 0$. Besides, the discrete $L^2$ norm of $y(\cdot,T)$ decreases with rate $h$, which is the expected one according to \eqref{eq:Sec3_functBehavior2}. Hence, Theorem \ref{thm:Sec3_HUMthm} confirms once again the theoretical controllability properties of Theorem \ref{thm:Sec3_controlUnconstrained}, assuring that the fractional heat equation in 2-D is null-controllable at time $T$.

\begin{SCfigure}[1.1][h]
	\centering
	\includegraphics[scale=0.25]{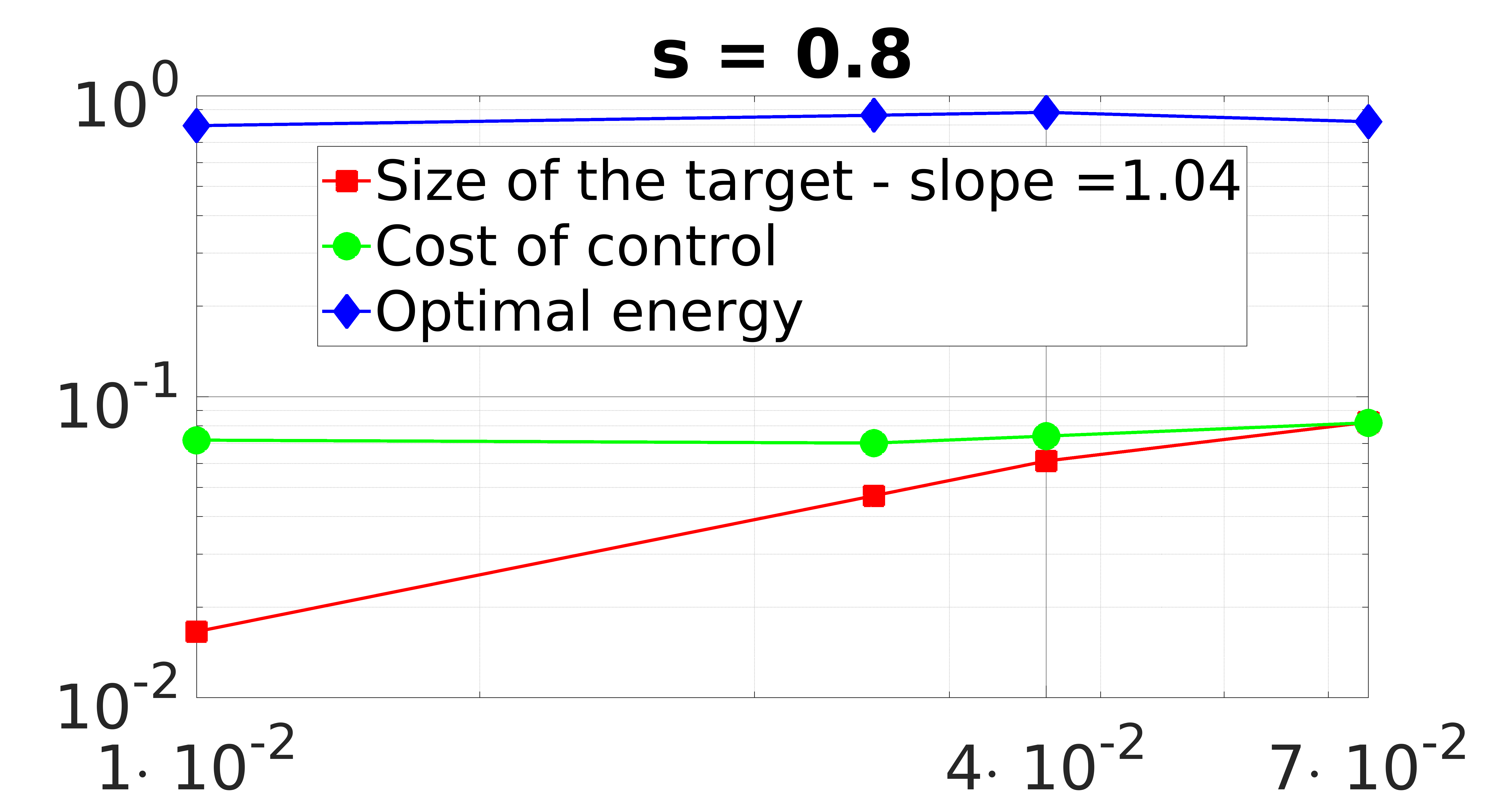}
	\caption{Behavior with respect to the mesh size $h$ of the cost of control, optimal energy and size of the solution to \eqref{eq:Sec3_HeatControlEuler} in space dimension $d=2$ at time $T$ when $s=0.8$.}\label{fig:Sec3_convergence2D08}
\end{SCfigure}

\subsubsection{The constrained controllability case}

In this section, we present some numerical evidences of the constrained controllability properties we obtained for \eqref{eq:Sec3_HeatControl} in Theorem \ref{thm:Sec3_controlConstraints}. In particular, we focus on the controllability to trajectories of \eqref{eq:Sec3_HeatControl} by means of non-negative controls. To this end, we choose the initial datum $y_0(x) = \sin(\pi x)$ and we set as a target $\widehat{y}(\cdot,T)$ the solution at time $T$ of \eqref{eq:Sec3_HeatControl} with initial datum $\widehat{y}_0(x) = 0.5\cos(\pi x/2)$ and right-hand side $\widehat{u}\equiv 0.02$ a.e. in $(-1,1)\times (0,T)$. Moreover, we choose $\omega=(-0.3,0.5)\subset (-1,1)$ as the control region and set $s=0.8$. Hence, Theorem \ref{thm:Sec3_controlConstraints} ensures the existence of a strictly positive controllability time $T_{\rm min}$ and a non-negative control function $u\in L^\infty(\omega\times(0,T))$ such that if $T\geq T_{\rm min}$ the solution of \eqref{eq:Sec3_HeatControl} is controllable to the trajectory $\widehat y(\cdot,T)$.

At this regard, we stress that we do not have analytical bounds for this minimal controllability time. As a matter of fact, the techniques developed in \cite{loheac2017minimal,pighin2018controllability} to obtain these bounds for the local heat equation are not immediately extendable to our fractional context. We refer to \cite[Section 4.4]{biccari2020positive} for a detailed discussion about this specific issue. Notwithstanding that, in what follows, we will provide numerical estimates of $T_{\rm min}$ by solving a suitable constrained optimization problem, which will also give us the minimal-time control $u_{\rm min}$ in the form of a linear combination of Dirac deltas, in accordance with Theorem \ref{thm:Sec3_controlConstraints}. Secondly, we will consider the controllability problem in time $T> T_{\rm min}$. As our simulations will show, also in that case the target trajectory is matched at time $T$, although the control loses its atomic nature. Finally, we will consider a short time horizon $T< T_{\rm min}$, in which the controllability to trajectories is not achieved.

\subsubsection{Numerical approximation of the minimal controllability time}

We start by estimating numerically the minimal controllability time $T_{\rm min}$. To this end, we employ \texttt{IpOpt} and \texttt{CasADi} to solve the following constrained optimization problem:
\begin{align}\label{eq:Sec3_Topt_constrT}
	\textrm{minimize}\;T
\end{align}
subject to the constraints
\begin{align}\label{eq:Sec3_Topt_constr}
	& T>0, \notag
	\\
	&(y_h^m)_{m=1}^M \mbox{ solves the fully-discrete dynamics }\eqref{eq:Sec3_HeatControlEuler}
	\\
	&(y_h^m)_{m=1}^M, (u_h^m)_{m=1}^M\geq 0. \notag
\end{align}

\begin{figure}[h]
	\centering
	\begin{minipage}{0.4\textwidth}
		\includegraphics[scale=0.23]{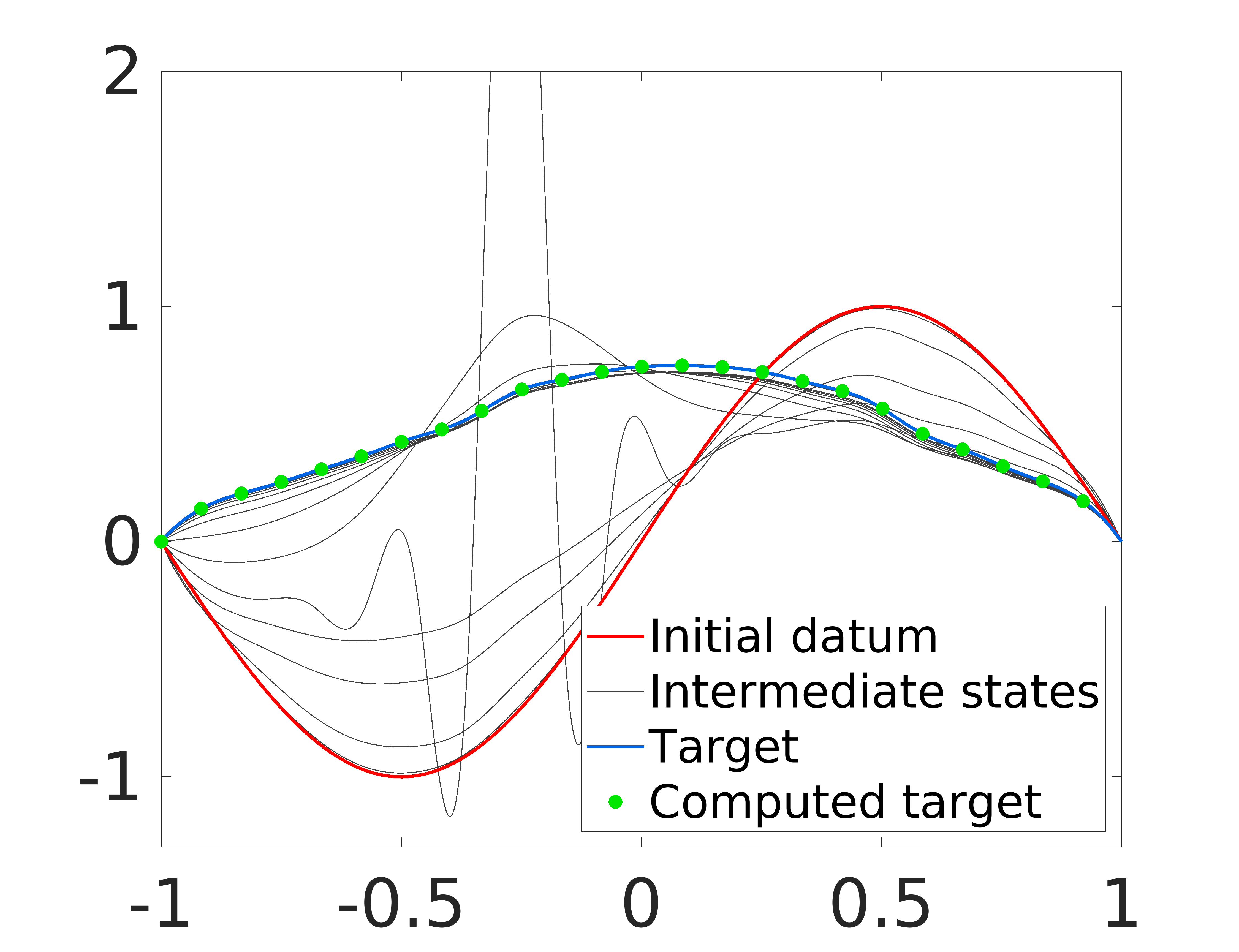}
	\end{minipage}
	\begin{minipage}{0.4\textwidth}
		\includegraphics[scale=0.23]{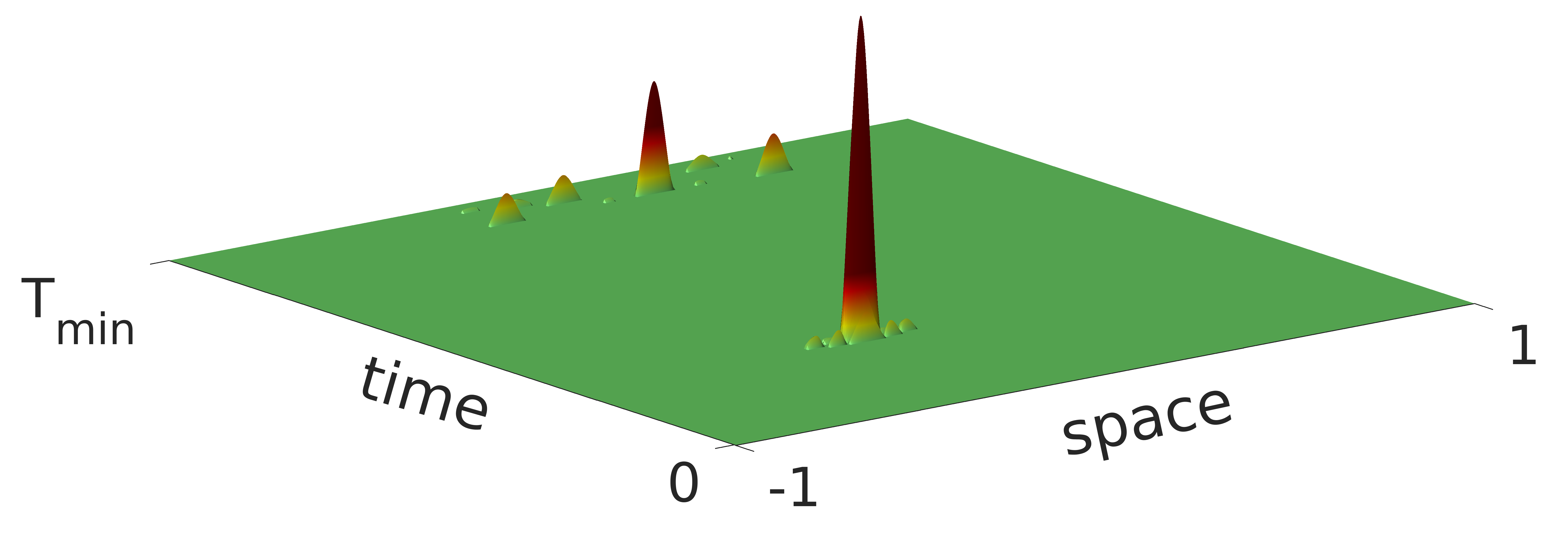}
	\end{minipage}
	\caption{Evolution in the time horizon $T=T_{\rm min}$ of the solution of \eqref{eq:Sec3_HeatControl} with $s=0.8$ (left) and of the control $u_{\rm min}$ (right). The blue curve is the target we want to reach while the green bullets indicate the final state of the computed numerical solution.}\label{fig:Sec3_constrainedTmin}
\end{figure}

By solving \eqref{eq:Sec3_Topt_constrT}-\eqref{eq:Sec3_Topt_constr}, we obtain $T_{\rm min}\simeq 0.68$ and the control $u_{\rm min}$. Figure \ref{fig:Sec3_constrainedTmin} shows the solution to \eqref{eq:Sec3_HeatControl} which, under the action of this minimal-time control, is steered from the initial datum $y_0$ to the desired target. Moreover, we can also see that the minimal-time control $u_{\rm min}$ is localized in certain specific points of the domain and time instants. This is in accordance with Theorem \ref{thm:Sec3_controlConstraints} which states that the minimal-time control is a Radon measure (in particular, a linear combination of Dirac masses).

\subsubsection{Constrained controllability in large time}

\begin{figure}[h]
	\centering
	\begin{minipage}{0.4\textwidth}
		\includegraphics[scale=0.23]{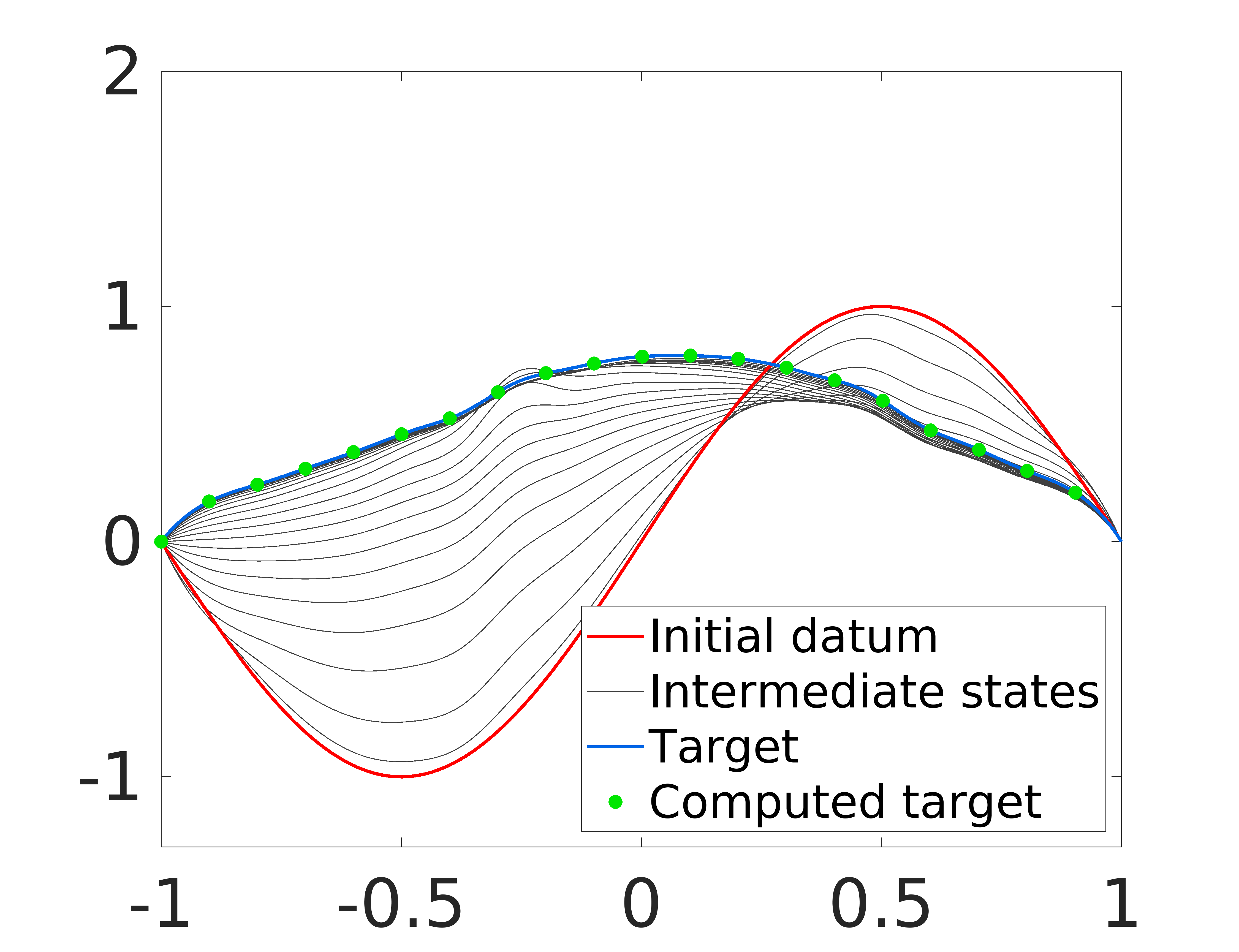}
	\end{minipage}
	\begin{minipage}{0.4\textwidth}
		\includegraphics[scale=0.23]{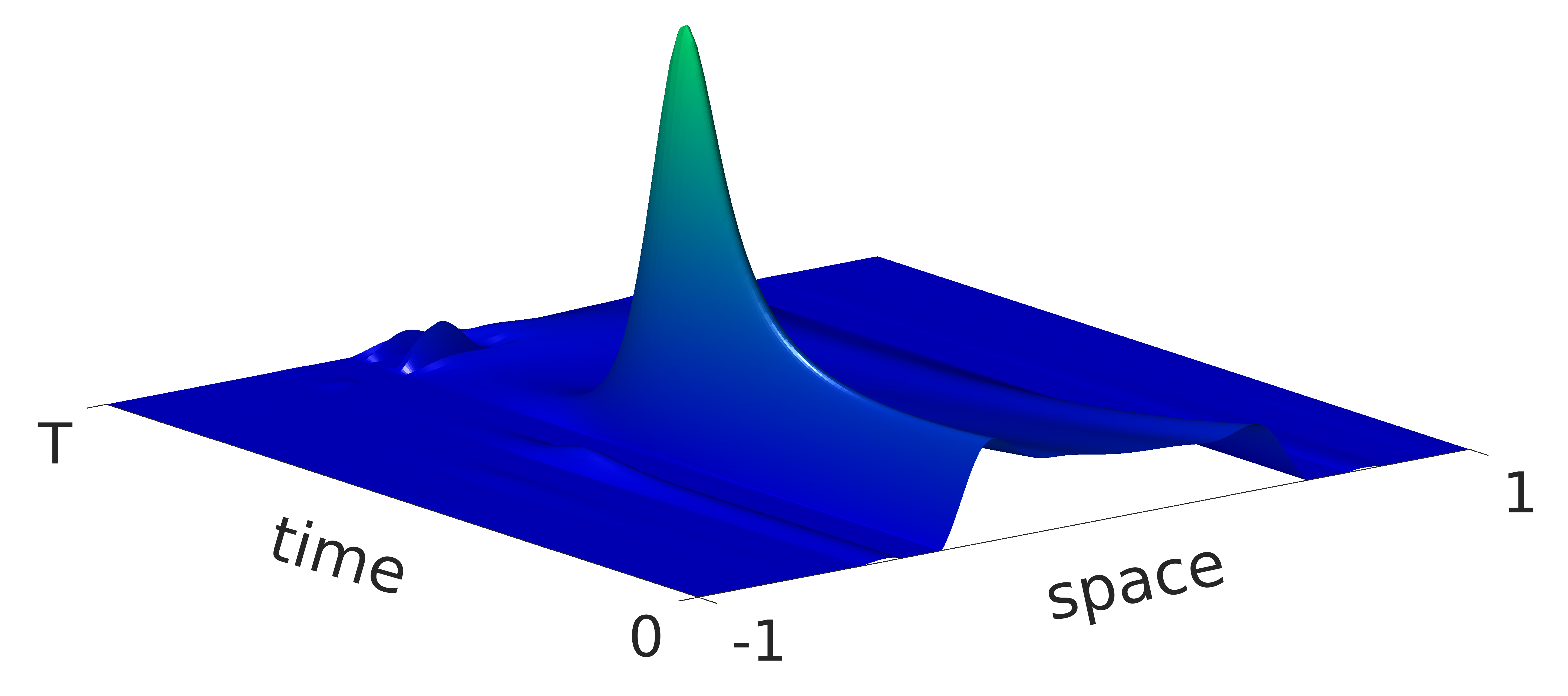}
	\end{minipage}
	\caption{Evolution in the time horizon $T=1>T_{\rm min}$ of the solution of \eqref{eq:Sec3_HeatControl} with $s=0.8$ (left) and of the control $u$ (right), under the constraint $u\geq 0$. The atomic nature of the control is lost.}\label{fig:Sec3_constrainedTlarge}
\end{figure}

The atomic behavior of the control is lost when extending the time horizon beyond $T_{\rm min}$. In Figure \ref{fig:Sec3_constrainedTlarge}, we show the evolution of the solution of the fractional heat equation \eqref{eq:Sec3_HeatControl} from the initial datum $y_0$ to the target $\widehat{y}(\cdot,T)$ in the time horizon $T=1$, and the corresponding optimal control that we computed minimizing in \texttt{IpOpt} the cost functional
\begin{align}\label{eq:Sec3_FunctLinfty}
	F_{\beta,\infty}(u)=\frac 12\norm{u}{L^\infty(\omega\times(0,T))} + \frac{1}{2\beta}\norm{y(\cdot,T)-\widehat y(\cdot,T)}{L^2(-1,1)}^2
\end{align}
with the constraints \eqref{eq:Sec3_Topt_constr}. As we can observe, in accordance with our theoretical results, the equation is still controllable in time $T$. Nevertheless, the control has lost its atomic nature and its action is now more distributed in $\omega$. 

\subsubsection{Lack of controllability in short time}

To conclude this section, let us now consider the controllability problem in a short time horizon $T<T_{\rm{min}}$ (more specifically, $T=0.25$). Once again, we have employed \texttt{IpOpt} combined with \texttt{CasADi} to minimize the functional \eqref{eq:Sec3_FunctLinfty} giving us the optimal control. 

\begin{figure}[h]
	\centering
	\begin{minipage}{0.4\textwidth}
		\includegraphics[scale=0.23]{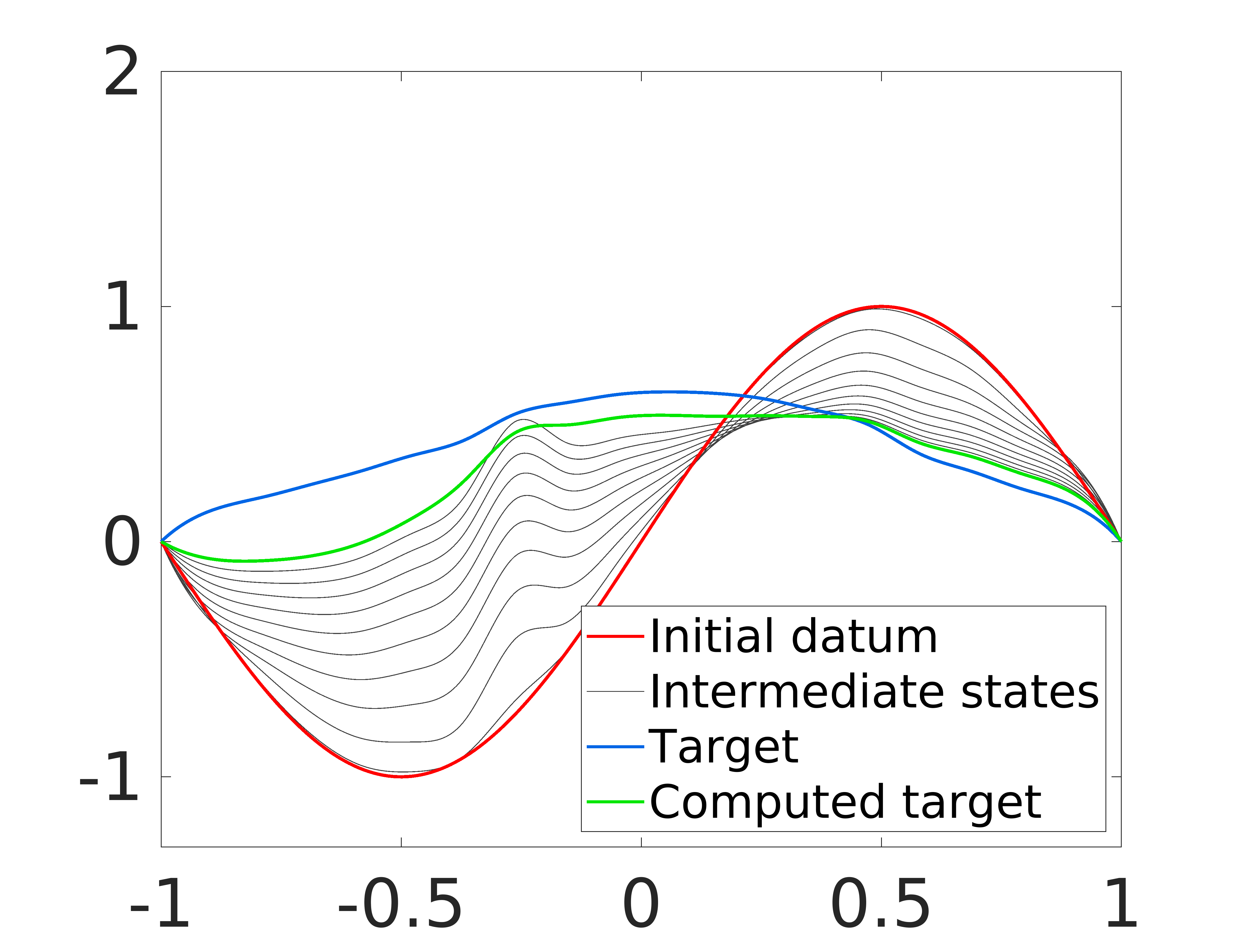}
	\end{minipage}
	\begin{minipage}{0.4\textwidth}
		\includegraphics[scale=0.23]{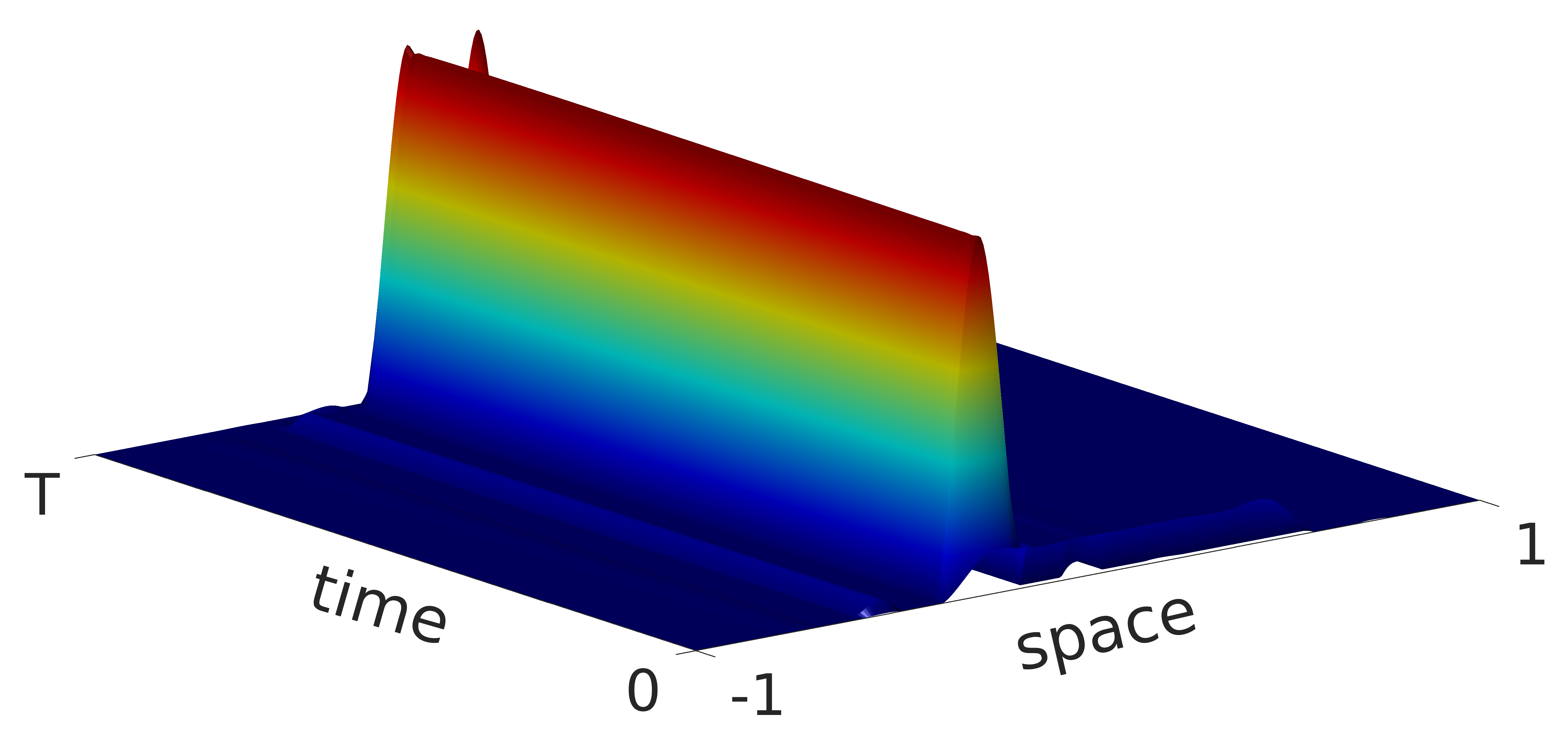}
	\end{minipage}
	\caption{Evolution in the time horizon $T=0.25<T_{\rm min}$ of the solution of \eqref{eq:Sec3_HeatControl} with $s=0.8$ (left) and of the control $u$ (right), under the constraint $u\geq 0$. The equation is not controllable.}\label{fig:Sec3_constrainedTsmall}
\end{figure}

Our simulations displayed in Figure \ref{fig:Sec3_constrainedTsmall} show that the solution of \eqref{eq:Sec3_HeatControl} fails to be controlled. In particular, the numerical solution at time $T$ does not match with the trajectory in blue we want to reach. This is in accordance with the lack of constrained controllability in short time as proved in Theorem \ref{thm:Sec3_controlConstraints}.

\section{Exterior controllability properties of the fractional heat equation}\label{sec:4}

In this section, we discuss the exterior controllability properties of the one-dimensional fractional heat equation. That is, for any $y_0\in L^2(-1,1)$ and $\mathcal O$ bounded, nonempty and open subset of $(-1,1)^c$, we are going to consider the following control system
\begin{align}\label{eq:Sec4_HeatControlExt}
	\begin{cases}
		y_t + \fl{s}{y} = 0 &\mbox{ in }\, (-1,1)\times(0,T),
		\\
		y = g\chi_{\mathcal O} &\mbox{ in } (-1,1)^c\times(0,T),
		\\
		y(\cdot,0) = y_0 & \mbox{ in }(-1,1).
	\end{cases}
\end{align}

Here $g$ is the control function belonging to some functional space to be specified later. 

The concept of exterior controllability for fractional models has been introduced in the literature only recently. First, in \cite{antil2019external,antil2020external}, the authors analyzed several elliptic and parabolic exterior optimal control problems. On the other hand, exterior controllability problems have been addressed in \cite{antil2020controllability,warma2019approximate,warma2018null}. At this regard, we shall recall that, as it has been shown in \cite{warma2019approximate}, a boundary control (that is, a control $g$ localized in a subset of the boundary) does not make sense in the presence of a fractional Laplacian. This is because of the non-locality of the operator and the fact that fractional models with standard boundary conditions (Dirichlet, Neumann or Robin) are ill-posed. For problems involving the fractional Laplacian the correct notion of a boundary controllability is actually the exterior one, requiring the control function to be localized outside the domain where the PDE is satisfied, as in \eqref{eq:Sec4_HeatControlExt}. Let us mention that exterior control problems also appear in many realistic applications, such as for instance:
\begin{itemize}
	\item[1.] Magnetic drug delivery: the drug with ferromagnetic particles is injected in the body and an external magnetic field is used to steer it to a desired location. 
	\item[2.] Acoustic testing: the aerospace structures are subjected to the sound from the loudspeakers.
\end{itemize}

We refer to \cite{antil2019external,antil2020external} and the references therein for a further discussion and the derivation of the exterior control.

In this section, we will give a broad panorama of the exterior controllability problem \eqref{eq:Sec4_HeatControlExt}. We will start by recalling in Subsection \ref{subs:Sec4_ExtControlResults} the theoretical controllability results presently available in the literature, and commenting important aspects on the numerical approximation of exterior controls. Secondly, in Subsection \ref{subs:Sec4_numerics}, we will present our numerical experiments.

\subsection{Review of theoretical controllability results}\label{subs:Sec4_ExtControlResults}

We summarize the theoretical controllability results which can be currently found in the literature for the fractional heat equation \eqref{eq:Sec4_HeatControlExt}.

\begin{theorem}\label{thm:Sec4_controlUnconstrained}
For the fractional heat equation \eqref{eq:Sec4_HeatControlExt}, the following results hold.
\begin{itemize}
	\item[1.] \textbf{Approximate controllability}. Let $\mathcal O\subset(-1,1)^c$ be any nonempty and open subset of $(-1,1)^c$ and $s\in(0,1)$. For any $T>0$, $y_0,y_T\in L^2(-1,1)$ and $\varepsilon>0$, there exists a control function $g\in \mathcal D((-1,1)^c\times(0,T))$ such that the unique solution $y$ of \eqref{eq:Sec4_HeatControlExt} satisfies $\norm{y(\cdot,T)-y_T}{L^2(-1,1)}\leq\varepsilon$.
	
	\item[2.] \textbf{Null controllability}. Let $\mathcal O\subset(-1,1)^c$ be any nonempty and open subset of $(-1,1)^c$. For any $T>0$ and $y_0\in L^2(-1,1)$, there exists a control function $g\in L^2(0,T;H^s((-1,1)^c))$ such that the unique solution $y$ of \eqref{eq:Sec4_HeatControlExt} satisfies $y(\cdot,T)=0$ a.e. in $(-1,1)$, if and only if $s\in(1/2,1)$.
		
	\item[3.] \textbf{$L^\infty$ null controllability}. Let $\mathcal O\subset(-1,1)^c$ be any nonempty and open subset of $(-1,1)^c$. For any $T>0$ and $y_0\in L^2(-1,1)$, there exists a control function $g\in L^\infty(\mathcal O\times(0,T))$ such that the unique solution $y$ of \eqref{eq:Sec4_HeatControlExt} satisfies $y(\cdot,T)=0$ a.e. in $(-1,1)$, if and only if $s\in(1/2,1)$.
\end{itemize}
\end{theorem}

For completeness, we shall mention that in \cite[Theorem 6]{warma2019approximate} the approximate controllability of \eqref{eq:Sec4_HeatControlExt} has actually been obtained in any space dimension $d\geq 1$. In Theorem \ref{thm:Sec4_controlUnconstrained}, we have stated this result in the one-dimensional context to remain consistent with the presentation of this Section. 

\subsubsection{Remarks on the theoretical controllability results}\label{subs:Sec4_remarks}

This section is devoted to some remarks and comments on the theoretical controllability results of Theorem \ref{thm:Sec4_controlUnconstrained}. These remarks will be at the basis of the methodology we will adopt in our forthcoming numerical simulations. 

First of all, for any function $\zeta\in H^s(\mathbb{R})$, let us denote with $\mathcal N_s\zeta$ the non-local normal derivative defined as (see \eqref{eq:App2_nonlocalDer})
\begin{align*}
	\mathcal N_s\zeta(x):= C_s\int_{-1}^1 \frac{\zeta(x)-\zeta(z)}{|x-z|^{1+2s}}\,dz,\quad x\in(-1,1)^c
\end{align*}

Moreover, for any $p_T\in L^2(-1,1)$, let $p\in L^2(0,T;H_0^s(\Omega))$ be the solution of the adjoint equation 
\begin{align}\label{eq:Sec4_adjoint}
	\begin{cases}
		-p_t + \fl{s}{p} = 0 &\mbox{ in }\, (-1,1)\times(0,T),
		\\
		p = 0 &\mbox{ in } (-1,1)^c\times(0,T),
		\\
		p(\cdot,T) = p_T & \mbox{ in }(-1,1).
	\end{cases}
\end{align}

Then multiplying \eqref{eq:Sec4_HeatControlExt} by $p$ and integrating over $(-1,1)\times(0,T)$ using the integration by parts formula given in Proposition \ref{prop:App2_prop}, it is easy to see that
\begin{align}\label{eq:Sec4_identity}
	y(\cdot,T) = 0 \mbox{ a.e. in } (-1,1) \mbox{ if and only if }\quad \int_{-1}^1y_0p(\cdot,0)\,dx = \int_0^T\int_{\mathcal O} g\mathcal N_s p\,dx.
\end{align}

Furthermore, the characterization \eqref{eq:Sec4_identity} yields that \eqref{eq:Sec4_HeatControlExt} is null controllable at time $T>0$ with $g\in L^2(0,T;H^s(\mathcal O))$ if and only if the following observability inequality for the adjoint system \eqref{eq:Sec4_adjoint} holds (see \cite[Lemma 2]{warma2018null}):
\begin{align*}
	\norm{p(\cdot,0)}{L^2(-1,1)}^2\leq C\int_0^T\norm{\mathcal N_s p(t)}{L^2(\mathcal O)}^2\,dt.
\end{align*}

This observability result has been proved in \cite[Theorem 1]{warma2018null} and follows by employing spectral techniques and parabolic Ingham inequalities, taking into account the Fourier decomposition for the solution of \eqref{eq:Sec4_HeatControlExt} given in Theorem \ref{thm:App2_existenceY}.

Finally, in view of the above considerations, we can see that the exterior control for \eqref{eq:Sec4_HeatControlExt} can be obtained from the following optimal control problem:
\begin{align}\label{eq:Sec4_OCP}
	&g_\beta = \min_{g\in L^2((0,T);H^s(\mathcal O))} F_\beta^{ext}(g)
	\\
	&F_\beta^{ext}(g) := \frac 12 \int_0^T \norm{g}{H^s(\mathcal O)}^2\,dt + \frac{1}{2\beta}\norm{y(\cdot,T)}{L^2(-1,1)}^2.\notag
\end{align}

Notice that the optimization process \eqref{eq:Sec4_OCP} has to be solved under the constraints given by the dynamics \eqref{eq:Sec4_HeatControlExt}. As we commented in Section \ref{subs:Sec2_FLext}, there are several different possibility to discretize this dynamics. Here, we will use the approach presented in \cite{antil2019external,antil2020external}, consisting in approximating \eqref{eq:Sec4_HeatControlExt} through the following exterior Robin problem 
\begin{align}\label{eq:Sec4_ExtHeatRobin}
	\begin{cases}
		y^n_t + \fl{s}{y^n} = 0 & \mbox{ in } (-1,1)\times(0,T),
		\\
		\mathcal{N}_s y^n+n \kappa y^n = n\kappa g &\mbox{ in } (-1,1)^c\times (0,T), 
		\\
		y^n(\cdot,0) = y_0 & \mbox{ in } (-1,1),
	\end{cases}
\end{align}
where $n\in\NN$ is a fixed natural number and $\kappa\in L^1((-1,1)^c)\cap L^\infty((-1,1)^c)$ is a non-negative function. 

The fact that \eqref{eq:Sec4_ExtHeatRobin} indeed approximates \eqref{eq:Sec4_HeatControlExt} has been proved in \cite{antil2020external}. As a matter of fact, we have the following result.

\begin{theorem}[{\cite[Theorem 5.3]{antil2020external}}]\label{thm:Sec4_Robin}
Let $y_0\in L^2(-1,1)$, $g\in L^2((0,T);H^s(-1,1)^c)$, $\kappa\in L^1((-1,1)^c)\cap L^\infty((-1,1)^c)$ non-negative and 
\begin{align*}
	y^n\in L^2((0,T);H^s_\kappa(-1,1))\cap H^1((0,T);H^{-s}_\kappa(-1,1))
\end{align*} 
be the weak solution of \eqref{eq:Sec4_ExtHeatRobin} according to Definition \ref{def:App2_weakSolRobin}. Let $y\in L^2((0,T);H^s(\RR))$ be the weak solution of \eqref{eq:Sec4_HeatControlExt}. There is a constant $C>0$, independent of $n$, such that
\begin{align}\label{approximation}
	\norm{y-y^n}{L^2(\RR\times(0,T))}\leq \frac Cn \norm{y}{L^2((0,T);H^s(\RR))}.
\end{align}
In particular, $y^n$ converges strongly to $y$ in $L^2((-1,1)\times(0,T))$ as $n\to+\infty$. 
\end{theorem}

We refer to Appendix \ref{appendixB} for the definition of the Hilbert space $H^s_\kappa(-1,1)$ and of the weak solutions to \eqref{eq:Sec4_ExtHeatRobin}. At this regard, we shall stress that, at the numerical level, passing to the limit rigorously as $n\to+\infty$ in \eqref{eq:Sec4_ExtHeatRobin} would require a better understanding of the problems mentioned in Section \ref{subs:Sec2_FLext} about the FE treatment of the exterior problem in the elliptic case. We shall comment more on this aspect in Section \ref{sec:6}.

Moreover, also in the case of the Robin problem \eqref{eq:Sec4_ExtHeatRobin} controllability can be characterized through a dual argument. To this end, for any $p_T^n\in L^2(-1,1)$ we denote with $p^n$ the solution of the following adjoint problem with Robin exterior conditions 
\begin{align}\label{eq:Sec4_AdjointRobin}
	\begin{cases}
		-p^n_t + \fl{s}{p^n} = 0 & \mbox{ in } (-1,1)\times(0,T),
		\\
		\mathcal{N}_s p^n+n \kappa p^n = 0 &\mbox{ in } (-1,1)^c\times (0,T), 
		\\
		p^n(\cdot,T) = p_T^n & \mbox{ in } (-1,1).
	\end{cases}
\end{align}
It can be readily checked that 
\begin{align*}
	y^n(\cdot,T) = 0 \mbox{ a.e. in } (-1,1) \mbox{ if and only if }\quad \int_{-1}^1y_0p^n(\cdot,0)\,dx + n\int_0^T\int_{\mathcal O} p^n\kappa g\,dx = 0.
\end{align*}

In view of that, the controllability of \eqref{eq:Sec4_ExtHeatRobin} $g$ replaced with $g\chi_{\mathcal O}$, where the control region $\mathcal O\subset (-1,1)^c$ is an arbitrary nonempty open set, is equivalent to the following observability inequality for \eqref{eq:Sec4_AdjointRobin}
\begin{align*}
	\norm{p^n(\cdot,0)}{L^2(-1,1)}^2\leq C\int_0^T\norm{p^n(t)}{L^2(\mathcal O)}^2\,dt,
\end{align*}
and the exterior control can be obtained from the following optimal control problem:
\begin{align}\label{eq:Sec4_OCProbin}
	&g_\beta = \min_{g\in L^2(\mathcal O\times (0,T))} G_\beta^{ext}(g)
	\\
	&G_\beta^{ext}(g) := \frac 12 \int_0^T \norm{g}{L^2(\mathcal O)}^2\,dt + \frac{1}{2\beta}\norm{y^n(\cdot,T)}{L^2(-1,1)}^2.\notag
\end{align}

To conclude this section, let us stress that the considerations which led to \eqref{eq:Sec4_OCProbin} are only formal, and currently not supported by rigorous mathematical results. We will comment more on this issue in Section \ref{sec:6}. 

\subsection{Numerical experiments}\label{subs:Sec4_numerics}

Let us present our numerical simulations for the exterior control problem. At this regard, we recall that, as discussed in the previous section and supported by Theorem \ref{thm:Sec4_Robin}, instead of the Dirichlet problem \eqref{eq:Sec4_HeatControlExt} we will consider the Robin problem \eqref{eq:Sec4_ExtHeatRobin}. The computation of the exterior control will be done through the optimization process \eqref{eq:Sec4_OCProbin}.

For completeness, we stress that replacing our original dynamics \eqref{eq:Sec4_HeatControlExt} with \eqref{eq:Sec4_ExtHeatRobin} introduces an approximation error associated with the parameter $n\in\NN$. Nevertheless, we shall also notice that this approximation error is of the order of $n^{-1}$ in the $L^2$-norm, and can be kept small (in comparison with the error introduced by the space/time discretization of the dynamics) by selecting $n$ large enough. In particular, we will take $n=10^9$.

To approximate \eqref{eq:Sec4_ExtHeatRobin}, we consider the interval $\mathcal I = (-2,2)\supset (-1,1)$ and assume that the control function $g$ is supported in $\mathcal O\subset ((-2,2)\setminus (-1,1))$. Notice that, in this case, the regularity required in Theorem \ref{thm:Sec4_Robin} for the function $\kappa$, namely $\kappa\in L^1((-2,2)\setminus(-1,1))\cap L^\infty((-2,2)\setminus(-1,1))$, simply reduces to $\kappa\in L^\infty((-2,2)\setminus(-1,1))$ and is fulfilled by considering $\kappa$ to be constant. For simplicity, we take $\kappa=1$. In other words, we approximate the following problem:
\begin{align}\label{eq:Sec4_ExtHeatRobinApprox}
	\begin{cases}
		y^n_t + \fl{s}{y^n} = 0 & \mbox{ in } (-1,1)\times(0,T),
		\\
		\mathcal{N}_sy^n + ny^n = ng\chi_{\mathcal{O}\times(0,T)} & \mbox{ in } ((-2,2)\setminus(-1,1))\times (0,T), 
		\\
		y^n(\cdot,0) = y_0 & \mbox{ in } (-1,1).
	\end{cases}
\end{align}

To discretize \eqref{eq:Sec4_ExtHeatRobinApprox} in space, we introduce a uniform $N$-points mesh on $(-2,2)$ with mesh-size $h$ and we use a globally continuous piece-wise linear FE scheme based on the following variational formulation (see Definition \ref{def:App2_weakSolRobin}): for all $v\in H^s_\kappa(-1,1)$, 
\begin{align}\label{eq:Sec4_RobinSol}
	&n\int_{\mathcal O} \kappa gv\,dxdt = \int_{-1}^1 y_t^n v\,dxdt + \mathcal F(y^n,v) + n\int_{(-1,1)^c}\kappa y^n v \,dxdt.
\end{align}

Notice that in \eqref{eq:Sec4_RobinSol} the non-homogeneous datum $g$ enters in the integral on the left-hand side and not in the bilinear form.

For the time discretization, instead, we apply implicit Euler on a uniform $M$-points grid discretizing $[0,T]$. 

To discuss the numerical controllability problem, in what follows, we adopt again the methodology of Section \ref{sec:3} and use Theorem \ref{thm:Sec3_HUMthm} to find numerical evidences of the control results of Theorem \ref{thm:Sec4_controlUnconstrained}. We do this by minimizing the functional \eqref{eq:Sec4_OCProbin} on several uniform meshes with decreasing mesh-size $h\to 0$ to compute the optimal control $g_\beta$, and by analyzing the behavior with respect to $\beta(h)$ (chosen as in \eqref{eq:Sec3_convergenceHeat}) of: 
\begin{itemize}
	\item The cost of control $\norm{g_\beta}{L^2((0,T);H^s(\mathcal O))}$.
	\item The optimal energy $G_\beta^{ext}(g_\beta)$.
	\item The $L^2$-norm of the corresponding solution at time $T$. 
\end{itemize}

\begin{SCfigure}[1.1][h]
	\centering
	\includegraphics[scale=0.25]{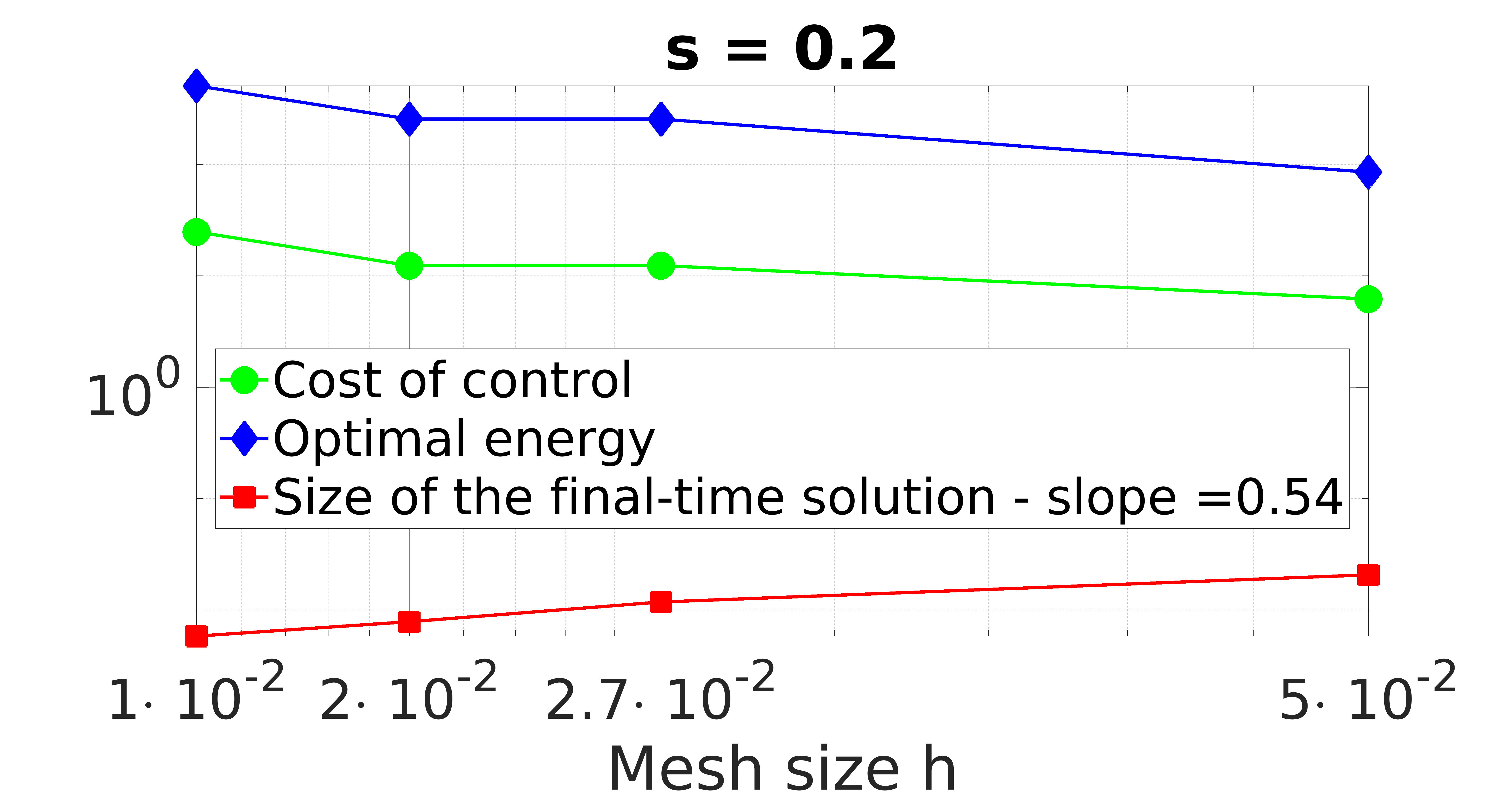}
	\caption{Behavior with respect to the mesh size $h$ of the cost of control, optimal energy and size of the solution to \eqref{eq:Sec4_ExtHeatRobinApprox} at time $T$ when $s=0.2$.}\label{fig:Sec4_convergence02}
\end{SCfigure}

We begin by considering $s=0.2$, for which we know that \eqref{eq:Sec4_HeatControlExt} is only approximately controllable. We set $\mathcal O=(1.7,1.9)$, $T=0.4$ and $y_0(x) = \cos(\pi x/2)$. We then employ \texttt{IpOpt} and \texttt{CasADi} to solve \eqref{eq:Sec4_OCProbin}. The results of these numerical experiments are displayed in Figure \ref{fig:Sec4_convergence02}.

We observe that the $L^2$ norm of the final state tends to zero as $h\to 0$, confirming computationally the approximate controllability of \eqref{eq:Sec4_HeatControlExt}. Notwithstanding that, we can also see that the cost of the control and the optimal energy slightly increase as $h\to 0$. Therefore, according to Theorem \ref{thm:Sec3_HUMthm}, the null controllability of \eqref{eq:Sec4_HeatControlExt} is not fulfilled.

\begin{SCfigure}[1.1][h]
	\centering
	\includegraphics[scale=0.25]{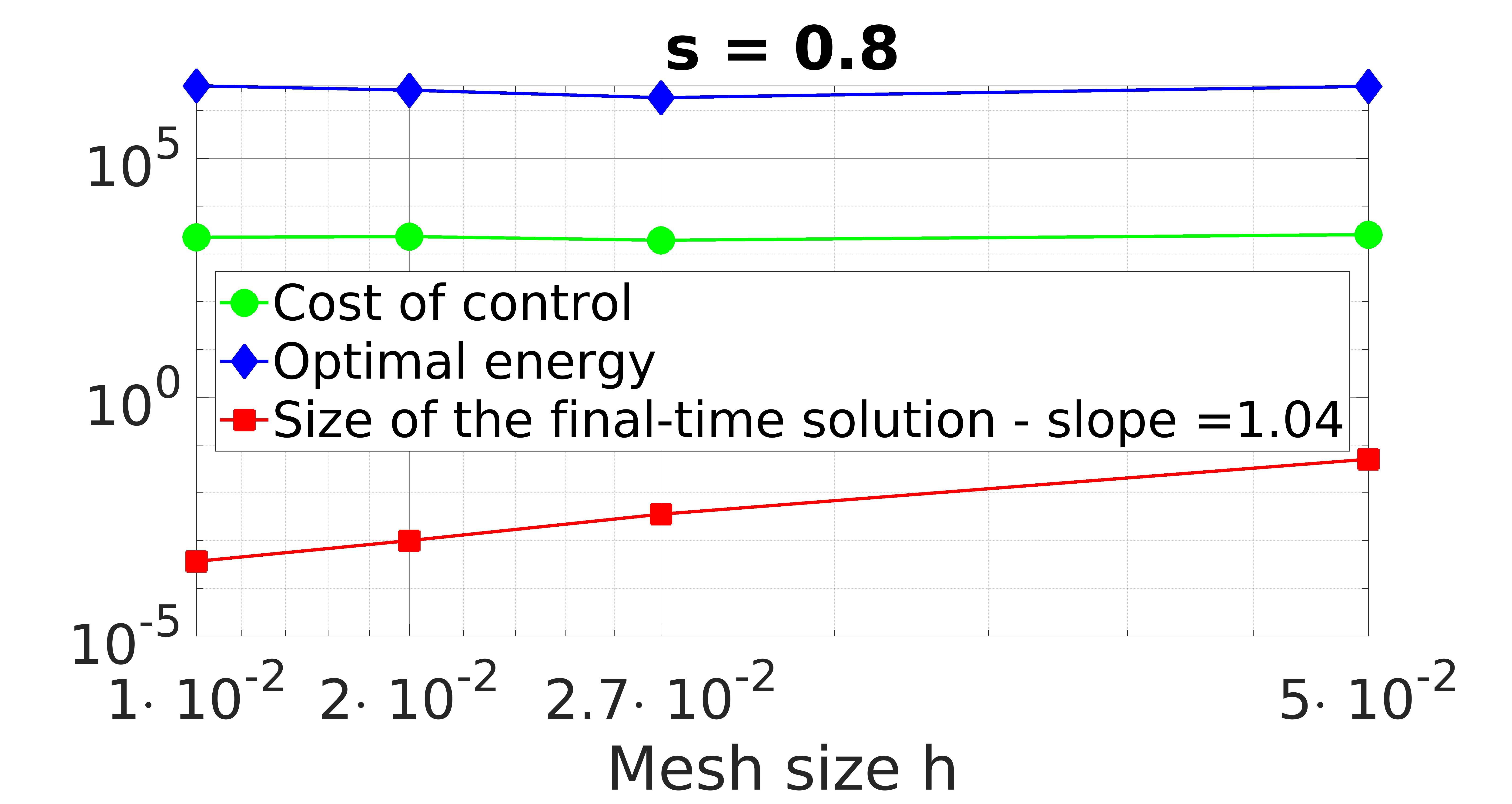}
	\caption{Behavior with respect to the mesh size $h$ of the cost of control, optimal energy and size of the solution to \eqref{eq:Sec4_ExtHeatRobinApprox} at time $T$ when $s=0.8$.}\label{fig:Sec4_convergence08}
\end{SCfigure}

Let us now take $s=0.8$. In this case, as we can see in Figure \eqref{fig:Sec4_convergence08}, the situation changes. Indeed, we can observe that this time the control cost and the optimal energy remain bounded as $h\to 0$. Furthermore, we also see that the $L^2$ norm of the solution at time $T$ decreases with the expected rate $h=\sqrt{\beta}$ (see \eqref{eq:Sec3_convergenceHeat}). According to Theorem \ref{thm:Sec3_HUMthm}, all these facts confirm that, for $s=0.8$, \eqref{eq:Sec4_HeatControlExt} is indeed null controllable. 

\begin{figure}[ht]
	\centering
	\includegraphics[scale=0.25]{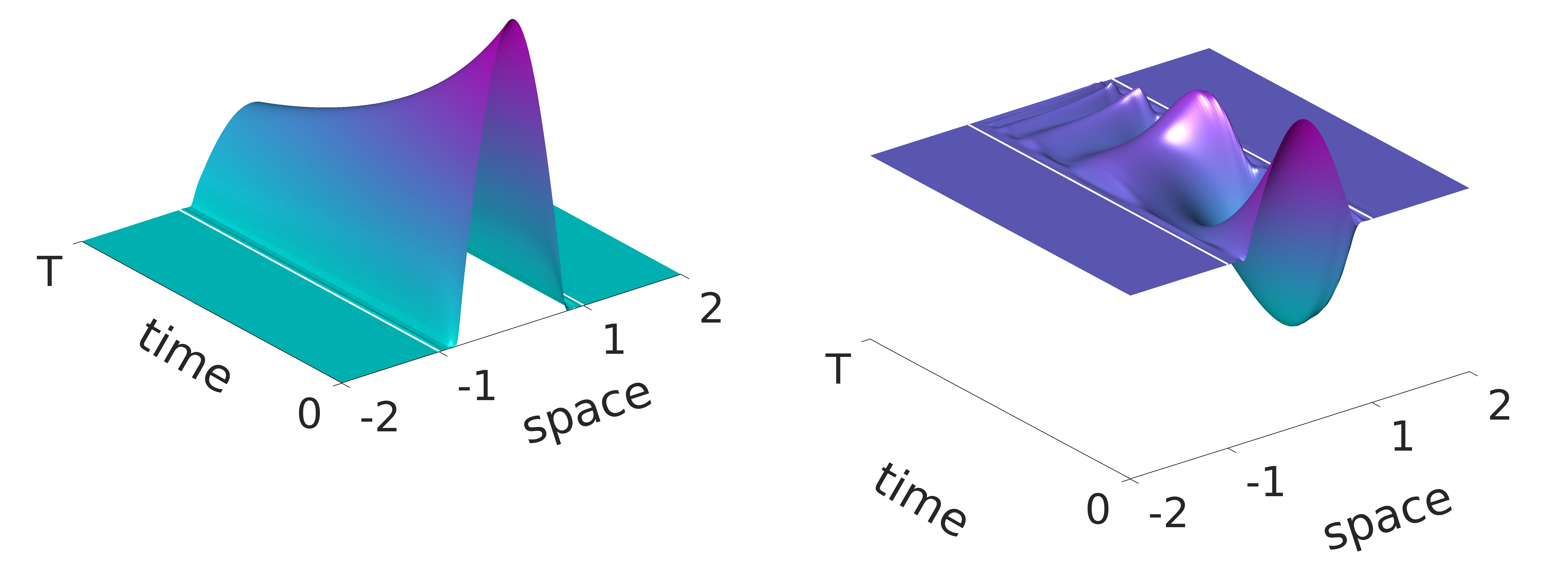}
	\caption{Free (left) and controlled (right) solution of \eqref{eq:Sec4_HeatControlExt} with $s=0.8$ at time $T=0.4$.}\label{fig:Sec4_solHeatExt}
\end{figure}

This positive controllability result is also appreciated in Figure \ref{fig:Sec4_solHeatExt}, where we can clearly see that the uncontrolled solution diffuses under the action of the fractional heat semi-group, but does not reach zero at time $T$. On the other hand, the introduction of the computed optimal control modifies the dynamical behavior of $y$ in such a way that we achieve $y(\cdot,T)=0$. 

\begin{SCfigure}[1.1][ht]
	\centering
	\includegraphics[scale=0.2]{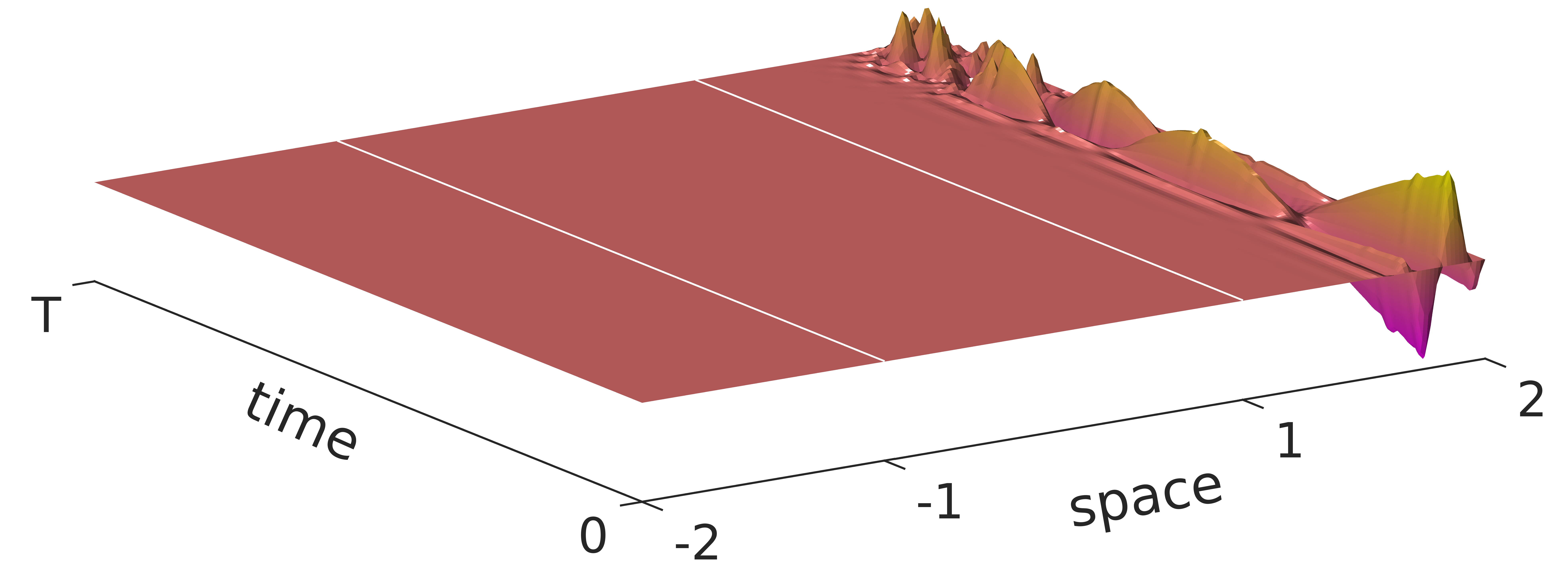}
	\caption{Exterior control obtained via the minimization of the functional \eqref{eq:Sec4_OCProbin}.}\label{fig:Sec4_controlHeatExt}
\end{SCfigure}

Finally, in Figure \ref{fig:Sec4_controlHeatExt}, we can see the time evolution of the control function, which is only acting outside the interval  $(-1,1)$ where the dynamics evolves. Specifically, the control is active only on the right part of the exterior domain, since we have chosen $\mathcal O = (1.7,1.9)$ as the control region.

\begin{remark}[Constrained controllability problem]
For completeness, we shall mention that the study of constrained controllability has been recently extended in \cite{antil2020controllability} to the exterior control problem \eqref{eq:Sec4_HeatControlExt}. In particular, in the same spirit of the interior control problem discussed in Section \ref{sec:3}, the following results have been established in \cite[Theorems 2.3, 2.4 and 2.5]{antil2020controllability}.
\begin{itemize}
	\item[1.] For all $s\in(1/2,1)$, $y_0\in L^2(-1,1)$, and any positive trajectory $\widehat{y}$, there exists a minimal strictly positive controllability time $T_{\rm min}>0$ such that, for all $T> T_{\rm min}$, we can find a non-negative control $g\in L^\infty(\mathcal O\times(0,T))$ whose corresponding solution $y$ of \eqref{eq:Sec4_HeatControlExt} satisfies $y(\cdot,T) = \widehat{y}(\cdot,T)$ a.e. in $(-1,1)$. Moreover, if $y_0\geq 0$, then $y(x,t)\geq 0$ for every $(x,t)\in (-1,1)\times (0,T)$. 	
		
	\item[2.] For $T = T_{\rm min}$, the above constrained controllability result holds with controls $g\in\mathcal M(\mathcal O\times(0,T_{\rm min}))$, the space of Radon measures on $\mathcal O\times(0,T_{\rm min})$.
\end{itemize}

Furthermore, the numerical computation of exterior non-negative controls has also been addressed in \cite[Section 5]{antil2020controllability}, by combining the techniques developed for the interior control problem \eqref{eq:Sec3_HeatControl} (and described in Section \ref{sec:3}) with the methodology we presented in Section \ref{subs:Sec4_remarks} for the numerical resolution of \eqref{eq:Sec4_HeatControlExt}.

An exhaustive discussion of the constrained exterior controllability properties of \eqref{eq:Sec4_HeatControlExt} (both at the continuous and discrete level) is omitted here for the sake of brevity, and can be found in \cite{antil2020controllability}.
\end{remark}
\section{Simultaneous control of parameter-dependent fractional heat equations}\label{sec:5}

The concept of simultaneous control arises naturally in many contexts including parameter-dependent models, transmission problems, PDE on graphs, synchronization or multi-agent systems with applications, for instance, to robot coordination. The aim is to design a unique control policy which is independent of the model's changes and robust in a variegated spectrum of different realizations.

In this section, we discuss the simultaneous control problem in the context of parameter-dependent fractional heat equations. 

One of the key issues when designing a control strategy for a dynamical system is the efficient computation of the control. This becomes an even more critical aspect in simultaneous control, where the corresponding optimization problem typically depends on a large amount of input data. For this reason, many analytical and computational techniques have been developed in the past years in order to speed up the simulation of parameterized control problems. Among others, we can mention Proper Orthogonal Decomposition, other more general Reduced Basis approaches, or the so-called \textit{greedy methodology} (\cite{hernandez2019greedy,lazar2020control,lazar2016greedy}). 

Here, we propose an alternative approach to simultaneous control, based on the employment of stochastic optimization techniques. We will start by presenting in Subsection \ref{subs:Sec5_formulation} the problem we are going to address. Secondly, in Subsection \ref{subs:Sec5_algorithms}, we will give an overview of several deterministic and stochastic optimization algorithms adapted to the computation of simultaneous controls, and we will discuss their convergence properties and computational cost. Finally, in Subsection \ref{subs:Sec5_numerics} we will present some numerical experiments to compare the efficiency of deterministic and stochastic optimization in the framework of simultaneous control.

\subsection{Problem formulation}\label{subs:Sec5_formulation}

Let $y_0\in L^2(-1,1)$, $\mathcal K = \{s_1,s_2,\ldots,s_{|\mathcal K|}\}\subset (1/2,1)$ a finite set of cardinality $|\mathcal K |$, and $s_\ell\in\mathcal K$ a random parameter following a uniform probability law $\mu$, that is, 
\begin{align*}
	\mu(\mathcal K) = 1 \quad\quad\text{ and }\quad\quad \mu(s_\ell) = |\mathcal K|^{-1}, \;\text{ for all } \ell = 1,\ldots,|\mathcal K|.
\end{align*}
We consider the parameter-dependent fractional heat equation 
\begin{align}\label{eq:Sec5_heatParam}
	\begin{cases}
		y_{{s_\ell},t} + \fl{{s_\ell}}{y_{s_\ell}} = u\chi_\omega &\mbox{ in }\, (-1,1)\times(0,T),
		\\
		y_{s_\ell} = 0 &\mbox{ in } (-1,1)^c\times(0,T),
		\\
		y_{s_\ell}(\cdot,0) = y_0 & \mbox{ in }(-1,1).
	\end{cases}
\end{align}

We are interested in the simultaneous control of \eqref{eq:Sec5_heatParam}, that is, we want to find a unique \textit{parameter-independent} control $u$ such that, at time $T>0$, the solution $y_{s_\ell}$ satisfies 
\begin{equation}\label{eq:Sec5_simulControlCond}
	y_{s_\ell}(\cdot,T)=0, \quad s_\ell\in \mathcal{K} \;\;\;\;\mu-\text{a. e.}
\end{equation}

In what follows, to simplify the notation, we will drop the sub-index $\ell$ and we will simply denote by $s$ any parameter in the set $\mathcal K$. The corresponding solution of \eqref{eq:Sec5_heatParam} will then be denoted by $y_s$. 

Simultaneous control problems are typically very difficult to be tackled. In particular, to determine whether a parameter-dependent system is simultaneous controllable is not a trivial task. 

For linear finite-dimensional (ODE) models, the question has been addressed and solved in \cite{loheac2016averaged}. Nevertheless, to the best of our knowledge, in the infinite-dimensional PDE setting the problem is still poorly understood. In particular, there are no results of simultaneous control in the sense of \eqref{eq:Sec5_simulControlCond}. We will present a more complete discussion on this point in Section \ref{sec:6}. 

In view of this fact, in what follows, we will relax the problem to simply consider the numerical implementation of simultaneous controls for \eqref{eq:Sec5_heatParam} via the minimization of a suitable cost function, leaving open the rigorous mathematical justification of the simulation evidences we will obtain. At this regard, we shall remark that, as pointed out in \cite[Remark 1.1-5]{loheac2016averaged}, the simultaneous control property \eqref{eq:Sec5_simulControlCond} is guaranteed by the fact that 
\begin{align*}
	\mathbb{E}\left[\norm{y_s(\cdot,T)}{L^2(-1,1)}^2\right] = \frac{1}{|\mathcal K|}\sum_{s\in\mathcal K}\norm{y_s(\cdot,T)}{L^2(-1,1)}^2 = 0, 
\end{align*}
where $\mathbb{E}[\cdot]$ denotes the expectation and we took into account the uniform probability distribution of the parameters. Hence, it is natural to address the simultaneous control of \eqref{eq:Sec5_heatParam} by solving the optimization problem
\begin{equation}\label{eq:Sec5_controlproblem}
	\begin{array}{l}
		\displaystyle u^\ast=\min_{u\in L^2(\omega\times(0,T))} F_s(u)
		\\[10pt]
		\displaystyle F_s(u):= \frac 12\int_0^T\norm{u(t)}{L^2(\omega)}^2\,dt + \frac{1}{2\beta} \mathbb{E}\left[\norm{y_s(\cdot,T)}{L^2(-1,1)}^2\right].
	\end{array} 
\end{equation}
Notice that, since $\mathbb{E}[\cdot]$ is convex, the functional $F_s$ is convex as well.

A classical way to address the minimization problem \eqref{eq:Sec5_controlproblem} would be to use the Gradient Descent (GD) or the Conjugate Gradient (CG) algorithm. Nevertheless, when applied to parameter-dependent problems, these approaches have a main drawback. Indeed, their implementation requires, in each iteration, to solve the state equation \eqref{eq:Sec5_heatParam} and the corresponding adjoint equation for all parameter values. This may rapidly increase the computational cost, especially when the dimension of $\mathcal K$ is large. 

To bypass this issue, a possible approach is to employ a stochastic algorithm to reduce the number of gradient calculations and, consequently, the total computational complexity. Here we will consider the well-known Stochastic Gradient Descent (SGD) algorithm (see \cite{robbins1951stochastic}). Our main goal in this section will be to analyze to which extent this stochastic approach may be successfully applied for solving \eqref{eq:Sec5_controlproblem}. 

\subsection{GD and SGD approaches}\label{subs:Sec5_algorithms}

In this section, we give an overall description of GD, CG and SGD for the resolution of \eqref{eq:Sec5_controlproblem}. In particular, we recall the main results about the convergence of these methods and their computational complexity. 

\subsubsection{The GD approach}

Let us start with the well-known GD procedure, consisting in finding the minimizer $\widehat{u}$ in \eqref{eq:Sec5_controlproblem} as the limit $k\to +\infty$ of the following iterative process
\begin{align}\label{eq:Sec5_GDalgo}
	u^{k+1}=u^k-\eta_k \nabla F_s(u^k),
\end{align}
where $\eta_{k}>0$ is called the \textit{step-size} or \textit{learning rate}. The selection of a correct learning rate is crucial for the algorithm performances. As a matter of fact, if $\eta_k$ is not properly chosen, \eqref{eq:Sec5_GDalgo} may actually not converge to the minimum of $F_s$. See e.g. \cite[Section 9.2]{boyd2004convex} or \cite[Section 2.1.5]{nesterov2004introductory} for more details. 

For what concerns the gradient $\nabla F_s$, this can be easily computed by means of a standard adjoint methodology. In our case, we can readily check that
\begin{align*}
	\nabla F_s(u) = u - \frac 1\beta\mathbb{E}[p_s\chi_\omega] = u - \frac{1}{\beta|\mathcal K|}\sum_{s\in\mathcal K}p_s\chi_\omega,
\end{align*}
where $p_s$ is the solution of the backward problem
\begin{align}\label{eq:Sec5_systemp}
	\begin{cases}
		-p_{s,t} + \fl{s}{p_s} = 0 &\mbox{ in }\, (-1,1)\times(0,T),
		\\
		p_s = 0 &\mbox{ in } (-1,1)^c\times(0,T),
		\\
		p_s(\cdot,T) = -y_s(\cdot,T) & \mbox{ in }(-1,1).
	\end{cases}
\end{align}
Consequently, the GD scheme to solve the optimization problem \eqref{eq:Sec5_controlproblem} becomes
\begin{equation}\label{eq:Sec5_GDscheme}
	\textbf{GD:} \quad u^{k+1} = u^k-\eta_k \left(u^k - \frac{1}{\beta|\mathcal K|}\sum_{s\in\mathcal K}p_s^k\chi_\omega\right).
\end{equation}
Hence, applying \eqref{eq:Sec5_GDscheme} for minimizing $F_s(u)$ requires to solve at each iteration $|\mathcal K|$ times the system 
\begin{align}\label{eq:Sec5_coupled-system}
	\begin{cases}
		y_{s,t} + \fl{s}{y_s} = u\chi_\omega &\mbox{ in }\, (-1,1)\times(0,T),
		\\
		-p_{s,t} + \fl{s}{p_s} = 0 &\mbox{ in }\, (-1,1)\times(0,T),
		\\
		y_s = p_s = 0 &\mbox{ in } (-1,1)^c\times(0,T),
		\\
		y_s(\cdot,0) = y_0,\; p_s(\cdot,T) = -y_s(\cdot,T) & \mbox{ in }(-1,1).
	\end{cases}
\end{align}

Concerning now the convergence rate, since the functional $F_s$ is convex, it is known (see, e.g., \cite[Theorem 2.1.15]{nesterov2004introductory} or \cite[Theorem 3.3]{nocedal2006numerical}) that if we take $\eta_k$ constant small enough we have that  
\begin{align}\label{eq:Sec5_GDconvergence}
	\norm{u^k-u^\ast}{L^2(\omega\times(0,T))}^2 \leq \norm{u^0-u^\ast}{L^2(\omega\times(0,T))}^2 e^{-2\mathcal C_{GD}k},
\end{align}
where the positive constant $\mathcal C_{GD}$ is given by 
\begin{align}\label{eq:Sec5_constGD}
	\mathcal C_{GD} = \ln\left(\frac{\rho+1}{\rho-1}\right),
\end{align}
$\rho$ denoting the problem's conditioning number. 

As for the computational effort, let us denote with $\mathfrak{C}$ the cost of solving \eqref{eq:Sec5_coupled-system} once. Then, the per-iteration cost of GD will be $\mathfrak{C}|\mathcal K|$, due to the need of solving \eqref{eq:Sec5_coupled-system} for all $s\in\mathcal K$ to compute $\nabla F_s(u^k)$. Combining this with \eqref{eq:Sec5_GDconvergence}, we then conclude that, for computing the control $\widehat{u}$ up to some given tolerance $\varepsilon>0$, i.e. 
\begin{align*}
	\norm{u^k-u^\ast}{L^2(\omega\times(0,T))}^2 <\varepsilon,
\end{align*}
the cost of the GD algorithm will be
\begin{align}\label{eq:Sec5_GDwork}
	cost_{GD} = \mathcal O\left(\frac{\mathfrak{C}|\mathcal K|\ln(\varepsilon^{-1})}{\mathcal C_{GD}}\right).
\end{align}

\subsubsection{The CG approach}\label{CGsec}

Let us now describe the CG approach and comment its convergence properties. 

CG is an efficient algorithm to solve linear systems (see \cite{saad2003iterative}). To apply it for minimizing $F_s$, the starting point is to notice that the gradient $\nabla F_s$ can be written as
\begin{align}\label{eq:Sec5_gradCG}
	\nabla F_s(u) = \underbrace{\Big(\,I + \Lambda\,\Big)}_{\mathbb{A}}u + \underbrace{\frac 1\beta\mathbb{E}[q_s\chi_\omega]}_{-b} = \mathbb{A}u-b,
\end{align}
where
\begin{itemize}
	\item[1.] The operator $\Lambda$ is defined as $\Lambda u = \mathbb E[p_s\chi_\omega]$ with $p_s$ computed solving  
	\begin{align*}
		\begin{cases}
			z_{s,t} + \fl{s}{z_s} = u\chi_\omega &\mbox{ in }\, (-1,1)\times(0,T),
			\\
			-p_{s,t} + \fl{s}{p_s} = 0 &\mbox{ in }\, (-1,1)\times(0,T),
			\\
			z_s = p_s = 0 &\mbox{ in } (-1,1)^c\times(0,T),
			\\
			z_s(\cdot,0) = 0,\;p_s(\cdot,T) = z_s(\cdot,T) & \mbox{ in }(-1,1).
		\end{cases}
	\end{align*}
	\item[2.] The function $q_s$ is obtained solving the system 
	\begin{align*}
		\begin{cases}
			\zeta_{s,t} + \fl{s}{\zeta_s} = 0 &\mbox{ in }\, (-1,1)\times(0,T),
			\\
			-q_{s,t} + \fl{s}{q_s} = 0 &\mbox{ in }\, (-1,1)\times(0,T),
			\\
			\zeta_s = q_s = 0 &\mbox{ in } (-1,1)^c\times(0,T),
			\\
			\zeta_s(\cdot,0) = y_0,\;q_s(\cdot,T) = \zeta_s(\cdot,T) & \mbox{ in }(-1,1).
		\end{cases}
	\end{align*}
\end{itemize}

Since, clearly, the minimizer $\widehat{u}$ of $F_s$ has to satisfy $\nabla F_s(\widehat{u}\,) = 0$, we see from \eqref{eq:Sec5_gradCG} that computing $\widehat{u}$ is equivalent to solve the linear system $\mathbb{A} u = b$, for which purpose we can use the CG algorithm. Nevertheless, we immediately see from the above discussion that this requires to solve at each iteration the state and adjoint equations $|\mathcal K|$ times. Hence, as for GD before, if $|\mathcal K|$ is large to employ CG to compute the simultaneous control may become a very demanding task. Concerning now the convergence rate, we know from \cite[Theorem 6.29, Equation 6.128]{saad2003iterative} that 
\begin{align}\label{eq:Sec5_CGconvergence}
	\norm{u^k-\widehat{u}\,}{L^2(\omega\times(0,T)}^2 \leq 4\norm{u^0-\widehat{u}\,}{L^2(\omega\times(0,T))}^2 e^{-2\mathcal C_{CG}k},
\end{align}
where the positive constant $\mathcal C_{CG}$ is given by 
\begin{align}\label{eq:Sec5_constCG}
	\mathcal C_{CG} = \ln\left(\frac{\sqrt{\rho}+1}{\sqrt{\rho}-1}\right).
\end{align}
From \eqref{eq:Sec5_CGconvergence} we get that, for achieving $\varepsilon$-optimality, the cost of the CG algorithm will be
\begin{align}\label{eq:Sec5_CGwork}
	cost_{CG} = \mathcal O \left(\frac{\mathfrak{C}|\mathcal K|\ln(\varepsilon^{-1})}{\mathcal C_{CG}}\right).
\end{align}

Finally, let us stress that the convergence properties of CG are known to be better than the GD ones. This is due to two main reasons. First, since by definition of conditioning number we have $\rho>1$, the constant $\mathcal C_{CG}$ in \eqref{eq:Sec5_constCG} is larger than $\mathcal C_{GD}$ given in \eqref{eq:Sec5_constGD}. Hence, even if both GD and CG algorithms converge exponentially, this convergence will actually be faster for CG. Moreover, CG is known to enjoy the so-called \textit{finite termination} property (see, e.g., \cite[Remark 2.4]{glowinski2015variational}). This means that, when solving a $N$-dimensional problem, the algorithm will converge in at most $N$ iterations. Practical implementations of CG may partially lose this finite termination property due to round-off errors. Nevertheless, this iterative method still provides monotonically improving approximations to the exact solution, which usually reach the required tolerance after a small (compared to $N$) number of iterations. See \cite[Section 6.11.3]{saad2003iterative} for more details.

\subsubsection{The SGD approach}

Let us now describe the SGD algorithm. The main difference with respect to GD and CG is that, in this iterative scheme, we do not employ all the components of $\nabla F_s(u)$. Instead, we pick a parameter $s_k$ i.i.d. from $\mathcal K$ (here the sub-index $k$ refers to the $k$-th iteration of the algorithm) with respect to the uniform probability $\mu$ and we use the corresponding gradient as descent direction. Hence, the SGD recursion process is given by
\begin{align*}
	u^{k+1} = u^{k}-\eta_k \nabla F_{s_k}(u^k),
\end{align*}
where $(\eta_k)_{k\geq 1}$ is a deterministic sequence of positive scalars which we still refer to as the \textit{learning rates sequence}. Moreover, in view of the computations we already presented for GD, the descent direction $\nabla F_{s_k}$ is given by
\begin{align*}
	\nabla F_{s_k}(u^k) = u^k - \frac 1\beta p_{s_k}^k\chi_\omega,
\end{align*}
with $p_{s_k}^k$ solution of the adjoint equation \eqref{eq:Sec5_systemp}. Hence, the complete SGD scheme to solve the optimization problem \eqref{eq:Sec5_controlproblem} is given by
\begin{align}\label{eq:Sec5_SGDscheme}
	\textbf{SGD:} \quad u^{k+1} = u^k-\eta_k \left(u^k - \frac 1\beta p_{s_k}^k\chi_\omega\right).
\end{align}

We then see that applying \eqref{eq:Sec5_SGDscheme} for minimizing the functional $F_s(u)$ requires, at each iteration $k$, only one resolution of the system \eqref{eq:Sec5_coupled-system}. Because of that, each iteration of this stochastic approach is very cheap.

Concerning now the convergence properties of SGD, some preliminary observations have to be made:
\begin{itemize}
	\item[1.] First of all, in the SGD method the iterate sequence $(u^k)_{k\geq 1}$ is a stochastic process whose behavior is determined by the random sequence $(s_k)_{k\geq 1}\subset\mathcal K$. In particular, this implies that the convergence properties of the algorithm have to be defined in terms of stochastic quantities, namely
	\begin{align*}
		\mathbb{E}\Big[\norm{u^k-u^\ast\,}{L^2(\omega\times(0,T))}^2\Big]
	\end{align*}  
	(see, e.g., \cite{bach2011non,bottou2018optimization}), or in the context of \textit{almost sure convergence} (\cite{bottou1998online}).
	
	\item[2.] As for the deterministic case, the choice of a good learning rate is crucial for the performances of the algorithm. In the stochastic framework, this is a delicate issue. At this regard, let us stress that choosing a constant learning rate is not a viable option for SGD. Indeed, if SGD is run with a fixed step-size $\eta_k = \overline{\eta}$, even if $\overline{\eta}$ is small we may not reach convergence (see \cite[Theorem 4.6]{bottou2018optimization}). This is essentially due to the noise introduced by the randomness of the stochastic process defined by \eqref{eq:Sec5_SGDscheme}, as it has been exhaustively discussed for instance in \cite[Section 4.2]{bottou2018optimization}. Because of this noise, the convergence of SGD is guaranteed as long one is able to maintain the second moment $\mathbb{E}\left[\norm{\nabla F_{s_k}(u^k)}{}^2\right]$ bounded above by a deterministic quantity, namely
	\begin{align*}
		\mathbb{E}\left[\norm{\nabla F_{s_k}(u^k)}{}^2\right]\leq \sigma^2, \quad 0<\sigma\in\mathbb{R}.
	\end{align*}
	This may be achieved by properly reducing the value of $\eta_k$ at each iteration. See, for instance, \cite{bach2011non,bottou2018optimization,robbins1951stochastic} for a complete discussion on this issue. 
	
	\item[3.] If the learning rate is properly chosen, by means of martingale techniques (see \cite[Section 4.5]{bottou1998online}) we can show that SGD converges almost surely
	\begin{align*}
		u^k \overset{a.s}{\longrightarrow} \widehat{u}, \quad \textrm{ as } k\to +\infty.
	\end{align*}
	In practice, this means that for solving \eqref{eq:Sec5_controlproblem} it is enough to run the SGD algorithm only once and we will have probability one of converging to the minimum $u^\ast$. 
\end{itemize} 

As for the convergence rate of SGD, if we choose the steps-size $\eta_k$ reducing as $\eta_k = k^{-\alpha}$, $\alpha\in (0,1)$, it has been proved in \cite[Theorem 1]{bach2011non} that 
\begin{align}\label{eq:Sec5_SGDconvergence}
	\mathbb{E}\left[\norm{u^k-\widehat{u}\,}{L^2(\omega\times(0,T))}^2\right] = \mathcal O\left(\sigma^2k^{-\alpha}\right).
\end{align}

We then immediately understand the importance of keeping $\sigma$ (which, we recall, is associated with the noise of the stochastic process) small, in order to have better convergence behavior. Moreover, confronting \eqref{eq:Sec5_SGDconvergence} with \eqref{eq:Sec5_GDconvergence} and \eqref{eq:Sec5_CGconvergence}, we see that SGD converges slower than GD and CG to the minimum $u^\ast$. Nevertheless, it is crucial to recall that, at each iteration $k$, SGD only requires to approximate the dynamics \eqref{eq:Sec5_coupled-system} once instead of $|\mathcal K|$ times, as it was for the deterministic algorithms. Hence we can conclude that for computing the control $\widehat{u}$ up to some given tolerance $\varepsilon>0$, i.e. 
\begin{align*}
	\mathbb{E}\left[\norm{u^k-u^\ast}{L^2(\omega\times(0,T))}^2\right] <\varepsilon,
\end{align*}
the cost of the SGD algorithm with step-size $\eta_k = k^{-\alpha}$, $\alpha\in (0,1)$, will be
\begin{align}\label{eq:Sec5_SGDwork}
	cost_{SGD} = \mathcal O\left(\mathfrak{C}\sigma^2\varepsilon^{-1}\right).
\end{align}

Of course, \eqref{eq:Sec5_SGDwork} is larger than \eqref{eq:Sec5_GDwork} and \eqref{eq:Sec5_CGwork} for $|\mathcal K|$ small, while the comparison favors SGD when $|\mathcal K|$ is large. This suggests that, when $|\mathcal K|$ is small, GD and CG are expected to perform better than SGD due to the lower amount of iterations they require to converge. On the other hand, in case of a large parameter set for which a single iteration of GD and CG becomes very costly, a stochastic approach will be more efficient than a deterministic one for solving \eqref{eq:Sec5_controlproblem}. Our numerical experiments in Section \ref{subs:Sec5_numerics} will confirm this behavior. 

As a final consideration, let us mention that the scheme \eqref{eq:Sec5_SGDscheme} is the most basic version of a SGD algorithm, as it was introduced in the original paper \cite{robbins1951stochastic}. In recent years, more sophisticated SGD-based algorithms have been proposed, to reduce the noise in the stochastic process and therefore achieve better convergence and stability properties. One of the most widely used nowadays is the so-called \textit{Adam} scheme (see \cite{kingma2015adam}), which combines the momentum approach proposed in \cite{polyak1964some} with a proper reduction of the step-size. In more detail, the Adam scheme consists in the following steps:
\begin{align}\label{eq:Sec5_Adam}
	&u^{k+1} = u^k - \frac{\eta}{\sqrt{\tilde v^k}+\delta}\widetilde m^k \notag 
	\\
	&\widetilde m^k = \frac{m^k}{1-\gamma_1}, \quad m^k = \gamma_1 m^{k-1} + (1-\gamma_1)\nabla F_{s_k}(u^k) 
	\\
	&\widetilde v^{\,k} = \frac{v^k}{1-\gamma_2}, \quad v^k = \gamma_1 v^{k-1} + (1-\gamma_1)|\nabla F_{s_k}(u^k)\,|^2, \notag 
\end{align} 
where $\eta$, $\gamma_1$ and $\gamma_2$ are suitable chosen parameters and $0<\delta\ll 1$ is introduced to avoid division by zero. This is the scheme that we will use for the simulations in Section \ref{subs:Sec5_numerics}.

\subsection{Practical considerations on the implementation of GD and CG}\label{subs:Sec5_considerations}

In \eqref{eq:Sec5_GDwork} and \eqref{eq:Sec5_CGwork}, we gave convergence rates for GD and CG in terms of the explicit constants $\mathcal C_{GD}$ and $\mathcal C_{CG}$ (see \eqref{eq:Sec5_constGD} and \eqref{eq:Sec5_constCG}), which depend on the conditioning $\rho$ of the problem we are solving. Nevertheless, if we analyze the behavior of these constants with respect to $\rho$, since by definition $\rho>1$, we immediately notice that both $\mathcal C_{GD}$ and $\mathcal C_{CG}$ are positive decreasing functions of $\rho$ and they converge to zero as $\rho\to +\infty$. This implies that a bad conditioning in a minimization problem affects the actual convergence of GD and CG, which may deteriorate and violate \eqref{eq:Sec5_GDconvergence} and \eqref{eq:Sec5_CGconvergence}. 

This is a well-known computational limitation of gradient optimization methods. In particular, the GD algorithm is very sensitive to the problem conditioning and, if $\rho$ is large, the convergence properties may deteriorate up to a linear rate. An illustrative example of this phenomenon is provided in \cite{meza2010steepest}. 

The situation is less critical for CG, because of the considerations we presented in the previous section. In particular, CG is less sensible to the conditioning of the problem since the constant $\mathcal C_{CG}$ depends on $\sqrt{\rho}$ instead of $\rho$ (see \eqref{eq:Sec5_constCG}). Furthermore, recall that CG enjoys the finite termination property, which helps in achieving convergence in a relatively small number of iterations.

In addition to that, let us mention that for a control problem $\rho$ is typically very large (see \cite[Remark 4.2]{boyer2013penalised} for some explicit estimates in the CG context). As a consequence of this bad conditioning, and according to the discussion above, we will see in our numerical simulations that the convergence rate of GD is significantly reduced with respect to the expected one given by \eqref{eq:Sec5_GDconvergence}, thus making this algorithm not very efficient in the context of our simultaneous control problem.

As a final remark, we stress that the considerations of this subsection concerning the conditioning of gradient methodologies apply also to SGD. Nevertheless, in this case, other key factors such as for instance the noise $\sigma$ come into play, thus making more delicate the precise quantification of the problem's conditioning and of its influence on the overall optimization problem.

\subsection{Numerical experiments}\label{subs:Sec5_numerics}

This section is devoted to some numerical experiments. The main goal is to confirm our previous discussion by comparing GD, CG and SGD for the simultaneous control of the linear parameter-dependent model \eqref{eq:Sec5_heatParam}. 

We have chosen the initial state $y_0=\sin(\pi x)$, the time horizon $T=0.4$ and the control region $\omega=(-0.5,0.8)$. The parameter set $\mathcal K$ is a $|\mathcal K|$-points discretization of the interval $(0.6,0.9)$, in which we know that each fractional heat equation \eqref{eq:Sec5_heatParam} is null-controllable. The tolerance is set to $\varepsilon = 10^{-4}$. To test the efficiency of each algorithm, we have performed simulations for increasing values of $|\mathcal K|$. For SGD, we have implemented the Adam scheme \eqref{eq:Sec5_Adam} with $\eta=10^{-3}$, $\gamma_1=0.9$, $\gamma_2 = 0.999$ and $\delta = 10^{-8}$.

\begin{SCfigure}[1][h]
	\centering
	\includegraphics[scale=0.25]{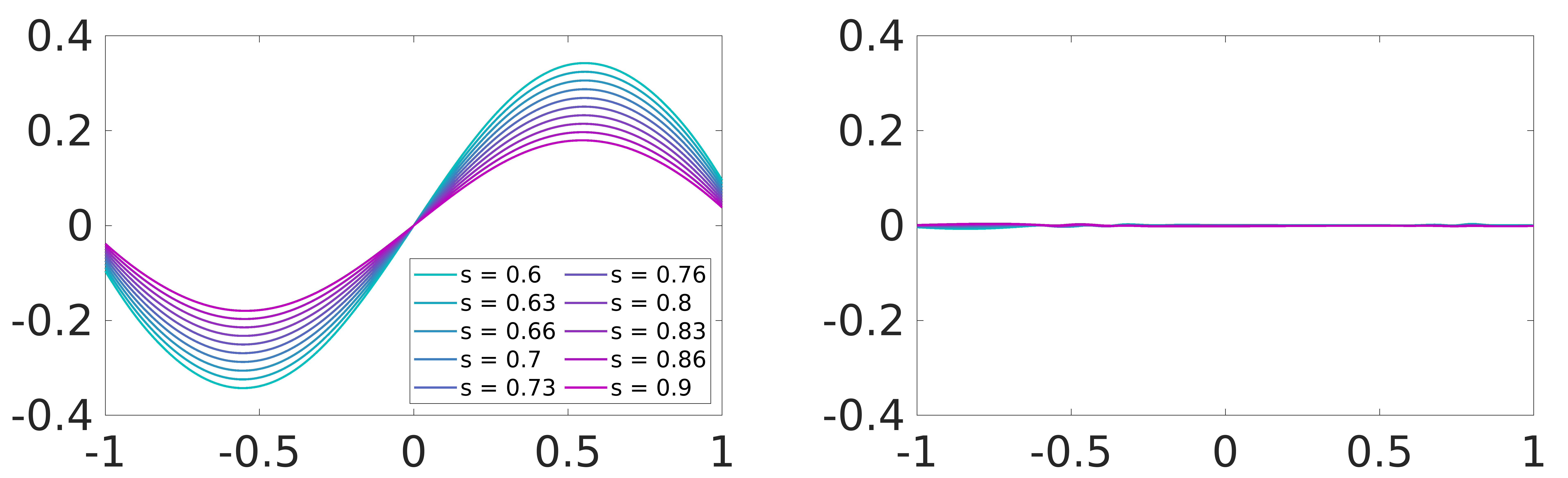}
	\caption{Final state $y_s(\cdot,T)$ of the free (left) and controlled (right) dynamics associated to \eqref{eq:Sec5_heatParam} for different values of $s\in (0.6,0.9)$.}\label{fig:Sec5_stateEvol}
\end{SCfigure}

Before comparing the performances of the three algorithms, let us show that the optimization problem \eqref{eq:Sec5_controlproblem} indeed provides an effective simultaneous control for \eqref{eq:Sec5_heatParam}. To this end, in Figure \ref{fig:Sec5_stateEvol}, we display the final state $y_s(\cdot,T)$ for the free (that is, when $u\equiv 0$) and controlled dynamics associated to \eqref{eq:Sec5_heatParam}. In order to increase the visibility of our plots, we consider here only the case of $|\mathcal K|=10$ parameters in our system. Moreover, the simulations displayed in Figure \ref{fig:Sec5_stateEvol} have been performed with the CG algorithm, although the other two approaches that we described provide the same result.

We see how, while the free dynamics simply dissipates under the action of the heat semi-group without reaching zero at the final time, introducing a control allows us to steer all the realizations of \eqref{eq:Sec5_heatParam} to zero at time $T$. For completeness, we shall stress that in the experiment displayed in Figure \ref{fig:Sec5_stateEvol} the initial datum $y_0$ is the same for all the realization of the parameter $s$. In principle, parameter-dependent initial condition could also be considered, and the sensitivity of the computed results with respect to these initial conditions could be analyzed. This kind of questions, although interesting, are beyond the interest of the present work and will not be discussed here.

\begin{SCfigure}[1.3][h]
	\centering
	\includegraphics[scale=0.28]{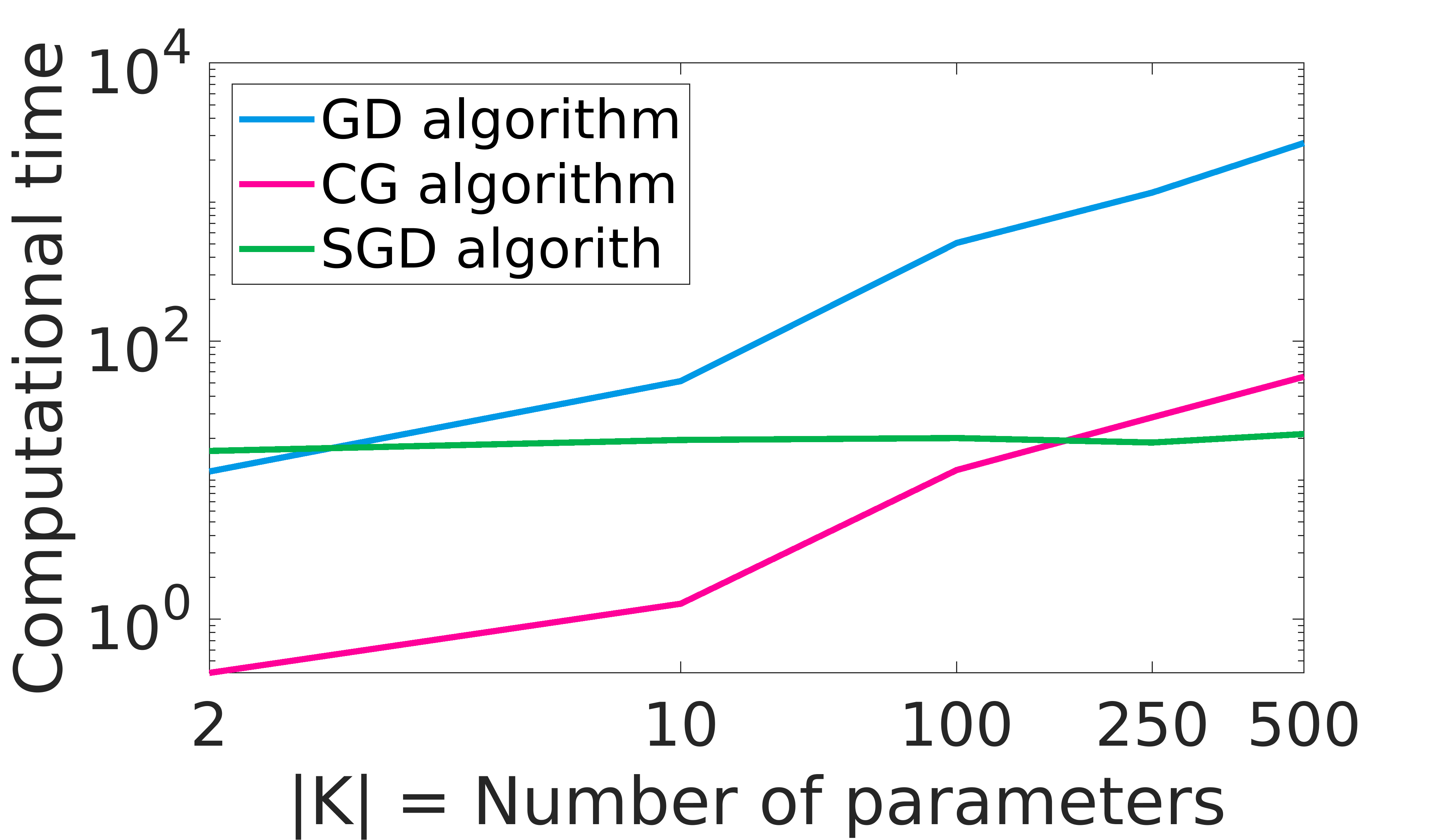}
	\caption{Computational time in logarithmic scale to converge to $\varepsilon = 10^{-4}$ for GD, CG and SGD applied to \eqref{eq:Sec5_controlproblem} with increasing values of $|\mathcal K|$.}\label{fig:Sec5_timeFig} 
\end{SCfigure}

Let us now analyze and discuss the behavior of GD, CG and SGD with respect to the amount of parameters included in our model. To this end, we have run simulations for increasing values of the cardinality of $\mathcal K$, namely $|\mathcal K| = 2, 10, 100, 250, 500$. The results of our numerical experiments are collected in Table \ref{tab:Sec5_timeTable} and displayed in Figure \ref{fig:Sec5_timeFig}.

\begin{table}[h]
	\centering 
	\begin{small}
	\begin{tabular}{|c|c|c|c|c|c|c|c|c|c|c|}
		\hline & \multicolumn{2}{|c|}{\textbf{GD}} & \multicolumn{2}{|c|}{\textbf{CG}} & \multicolumn{2}{|c|}{\textbf{SGD}} 
		\\
		\hline $|\mathcal K|$ & Iter. & Time (sec) & Iter. & Time (sec) & Iter. & Time (sec)
		\\
		\hline $2$ & $1424$ & $11.5$ & $35$ & $0.4$ & $3563$ & $16.2$
		\\
		\hline $10$ & $1363$ & $51.5$ & $25$ & $1.3$ & $3987$ & $19.4$
		\\
		\hline $100$ & $1343$ & $507.6$ & $25$ & $11.8$ & $4110$ & $20.1$
		\\
		\hline $250$ & $1341$ & $1169.3$ & $25$ & $28.3$ & $4019$ & $18.6$
		\\
		\hline $500$ & $1341$ & $2659.9$ & $25$ & $55.39$ & $4485$ & $21.4$
		\\
		\hline
	\end{tabular}
	\end{small}\caption{Number of iterations and computational time to converge to $\varepsilon = 10^{-4}$ for GD, CG and SGD applied to \eqref{eq:Sec5_controlproblem} with increasing values of $|\mathcal K|$.}\label{tab:Sec5_timeTable}
\end{table}

We notice that the number of iterations for GD and CG remains stable with respect to $|\mathcal K|$, while it grows for SGD. We believe that this is due to the noise introduced by the stochastic approach, although we do not have a rigorous mathematical justification of this phenomenon.
  
Nevertheless, despite this fact, we can see that our simulations confirm the behavior we expected from the discussion in Section \ref{subs:Sec5_algorithms}. In particular, we can make the following observations: 
\begin{itemize}
	\item[1.] GD is the worst performing algorithm. For $|\mathcal K|=2$, its computational time is comparable with SGD. Nevertheless, already for $|\mathcal K|=10$, its computational effort becomes considerable. This is a consequence of the bad conditioning of control problems, as we commented in Section \ref{subs:Sec5_considerations}. 
	\item[2.] CG is the algorithm requiring the lowest number of iterations to converge. On the one hand, this confirms that CG is less sensible to the conditioning of the problem. On the other hand, this implies that CG is the best approach when dealing with a low and moderate amount of parameters, since when $|\mathcal K|$ is not too large the algorithm is capable to compensate the per-iteration cost with the very limited amount of iterations it requires to achieve $\varepsilon$-optimality. 
	\item[3.] SGD is the algorithm requiring the highest numbers of iterations to converge. This is in line with the fact that, for this algorithm, only a linear or sub-linear convergence rate is expected (see \eqref{eq:Sec5_SGDconvergence}), while for GD and CG this rate is expected to be exponential (see \eqref{eq:Sec5_GDconvergence} and \eqref{eq:Sec5_CGconvergence}). Notwithstanding that, SGD appears to be insensitive to the cardinality of $\mathcal K$. This is not surprising if we consider that, no matter how many parameters enter in our control problem, with SGD each iteration of the optimization process always requires only one resolution of the coupled system \eqref{eq:Sec5_coupled-system}. Due to this low per-iteration cost, despite of the slower convergence rate, when $|\mathcal K|$ is large the total computational time for SGD is actually smaller than GD and CG. As a matter of fact, we can see in Figure \ref{fig:Sec5_timeFig} that for a parameter set of cardinality around $|\mathcal K|=150$ SGD starts outperforming CG. 
\end{itemize}

All these considerations are aligned with our previous discussion and corroborate the fact that, when $|\mathcal K|$ is large, the SGD approach is preferable to the GD and CG ones to address the simultaneous control of \eqref{eq:Sec5_heatParam}. 
\section{Conclusion and open problems}\label{sec:6}

In this work, we have addressed several aspects of the numerical control of non-local diffusive models involving the fractional Laplacian. Our discussion embraced a wide spectrum of situations: we considered both interior and exterior control problems, possibly under positivity constraints, as well as the relevant issue of simultaneous control. After having recalled the existing results in this fractional diffusive setting, we have shown how to compute the numerical controls by combining FE approximations of the fractional Laplacian with efficient optimization tools. Our simulations, on the one hand, have provided a numerical representation and validation of many known control properties for fractional diffusive processes. On the other hand, they have highlighted some analogies and differences between local and non-local heat-like models when facing the approximation of numerical controls. Nevertheless, many key questions related to our work remain currently unaddressed. We present below a short collection of them, which may be of interest for future investigation.

\begin{itemize}
	\item[1.] \textbf{Control of fully-discrete fractional problems.} Despite of the numerical evidences we presented in Sections \ref{sec:3} and \ref{sec:4} for the controllability properties of fully-discrete fractional heat equations, the complete analysis and proofs of these results are open. As a matter of fact, the controllability of fully-discrete problems has been considered only in a few research works and many important issues still remain unanswered. In \cite{ervedoza2010observability}, the authors proved in a quite general framework that any parabolic equation is null-controllable after time discretization by applying an adequate filtering of high frequencies, under the condition that the space semi-discrete approximation schemes are uniformly observable with respect to the mesh size parameters. We remark that, in the context of the discrete Laplace operator, the uniform controllability of semi-discretized solutions has been considered in \cite{lopez1998some,zuazua2006control} by means of spectral techniques and a careful analysis of the eigenvalues. Notice, however, that this analysis is not extendable to the discrete fractional Laplacian, due to the lack of an explicit knowledge of the spectrum of this operator. On the other hand, in \cite{boyer2011uniform}, the controllability of fully-discrete finite-difference heat equations has been studied by adopting a Lebeau-Robbiano strategy, based on discrete spectral inequalities. Finally, in \cite{hernandez2020carleman}, similar results have been obtained by means of a Carleman approach. To the best of the authors' knowledge, there are no results currently available apart from the aforementioned ones. In particular, the study of control properties for the fully-discrete approximation of \eqref{eq:Sec3_HeatControl} is still completely open.
	
	\item[2.] \textbf{Exterior control of the fractional heat equation.} In Section \ref{sec:4}, we considered the exterior control problem for the fractional heat equation and discussed its theoretical and numerical aspects. Nevertheless, several important issues are still only partially understood. On the one hand, from the theoretical viewpoint, the only existing controllability results for \eqref{eq:Sec4_HeatControlExt} have been obtained in space dimension $d=1$, while the multi-dimensional problem is completely open. At this regard, it would be interesting to analyze whether the techniques we recently developed in \cite{biccari2021null} for the interior control problem in multi-D are applicable also in the exterior control setting or, instead, new methodologies need to be derived. On the other hand, also the numerical control of \eqref{eq:Sec4_HeatControlExt} requires further investigation. In Section \ref{subs:Sec4_remarks}, we have proposed to compute the exterior controls passing through the Robin problem \eqref{eq:Sec4_ExtHeatRobin} since, as we discussed throughout this chapter, the numerical approximation of the exterior Dirichlet problem is still not completely clarified. Our simulations have shown the efficacy of this approach. Nevertheless, to completely justify this procedure, we shall rigorously analyze the control properties of the fully-discrete approximation of \eqref{eq:Sec4_ExtHeatRobin}, in the same spirit of the discussion in point 1 above. Finally, also the Fenchel-Rockafellar duality has not been developed yet int the context of the exterior control problem, and it would be considered in a future work. 
	
	\item[3.] \textbf{Simultaneous control of fractional PDE.} In Section \ref{sec:5}, we addressed the simultaneous control of parameter-dependent fractional heat equations and we proposed the use of stochastic optimization techniques for an efficient computation of numerical controls. Our discussion focused mainly on the implementation details, and the existence of a simultaneous control for the problem we considered has been addressed only at a numerical level. The reason behind this choice is that, at present time, the theory for simultaneous control of PDE is not sufficiently developed to tackle fractional models. As a matter of fact, existing techniques available in the local setting do not seem to be applicable in the non-local one. Just to give an example, in \cite{zuazua2016stable}, we have proposed a methodology for the stable observation of additive superpositions of heat and wave equations which is applicable to the simultaneous control of such models. This methodology is rather systematic and easy to apply. It consists in observing only one component of the system and considering the others as unknown perturbations. One then composes all but one (the observable one) PDE operators so to reduce the problem to the consideration of a single equation and applying its known observability properties. Notwithstanding that, this approach is hardly extendable in the context of the fractional heat equation we have considered, since the commutation of different powers of the fractional Laplacian leads to non-local lower order perturbations which cannot be handled in terms of observability. It would then be of interest to provide a solid mathematical background to simultaneous control for non-local and fractional models, supporting the numerical evidences we displayed in Section \ref{sec:5}. In addition to that, it would be worth to consider the extension of Fenchel-Rockafellar duality in this simultaneous control framework. 
	
	\item[4.] \textbf{Wave-type equations.} Our discussion has focused only on the parabolic setting. Nevertheless, in recent years, several results have been obtained on the controllability properties of hyperbolic (wave) and dispersive (Schr\"odinger) models involving the fractional Laplacian. The interested reader may refer, for instance, to \cite{biccari2020internal,biccari2020null,warma2020analysis}. As for the numerical approximation of wave-type models, this issue is known to be quite delicate. As a matter of fact, already in the local framework, numerical high-frequency solutions of the wave equation can exhibit pathological behaviors such as lack of propagation in space or the so-called \textit{rodeo effect}, i.e. waves that are trapped by the numerical grid in closed loops (see, e. g., \cite{biccari2020propagation,zuazua2005propagation}). When dealing with control and inversion problems, these behaviors then yield to the necessity of filtering high-frequency numerical components, to cope with the loss of uniform observability properties through numerical discretization. To the authors' knowledge, a rigorous analysis of how these mentioned pathologies transfer to the non-local setting of the fractional Laplacian has yet to be developed, and would be a very interesting question to be considered.
	
	\item[5.] \textbf{Variable-order fractional Laplacian.} In some recent contributions (\cite{antil2019sobolev,bahrouni2018new,rahmoune2021multiplicity,xiang2019multiplicity}), elliptic problems involving a variable-order fractional Laplacian (that is, with $s=s(x):\Omega\to (0,1)$) have been considered, analyzing the existence of (possible several) solutions and some optimal control issues applied to image denoising. It would be of interest to investigate parabolic problems (which,as far as we can tell, are still unaddressed) and associated theoretical and numerical control problems.  
	
	\item[6.] \textbf{Memory-type equations.} Evolution equations involving memory terms appear in several different applications, for modeling
	natural and social phenomena which, apart from their current state, are influenced also by their history. Some classical examples are viscoelasticity or non-Fickian diffusion. In recent years, this class of models has got the attention of the control community (see for instance \cite{biccari2019memory,chaves2017controllability,lu2017null} or \cite{biccari2020null} for models involving the fractional Laplacian). In particular, it has been observed that these models can be cast as coupled  PDE-ODE systems, in which the ODE component introduces non-propagation effects similar to those produced by a low-order ($s\leq 1/2$) fractional Laplacian. This motivated the introduction of a \textit{moving control} strategy to obtain positive controllability results. An illustrative example of this approach can be found in \cite{post}. Nevertheless, as far as the authors know, the controllability of memory-type equation has never been considered under a numerical perspective. In particular, the implementation of the moving control strategy in the hyperbolic setting, which already requires a special treatment of high-frequency solutions (see Point 4 above) would be an important issue to be further investigated. 		
	
	\item[7.] \textbf{Non-linear fractional control} The control of non-linear dynamical systems is most often a very difficult task, since a non-linear dynamics, in contrast with much simpler linear systems, may appear chaotic, unpredictable, or counter-intuitive. In the last decades, different techniques have been develop to control non-linear models, such as the employment of linearization and fixed-point arguments, or the so-called return method (see \cite{coron2007control}). Nevertheless, these techniques are most of the time of difficult application to fractional models, whose complex and pathological behaviors are typically enhanced in a non-linear framework. As for the numerical implementation of the controls, also this is a very difficult question. In some simple situation (see e.g. \cite[Section 4.3.2]{boyer2013penalised}) the numerical control problem for semi-linear models has been tackled by combining the HUM approach we presented in Section \ref{subs:Sec3_PenalizedHUM} with some fixed-point strategy. However, the extension to general non-linear (possibly fractional) problems is, to the best of our knowledge, still missing. All this renders non-linear fractional control an almost unexplored yet very interesting field within control theory.
\end{itemize}

{\appendix
	\section{Fractional order Sobolev spaces and the fractional Laplacian}\label{appendixA}

In this Appendix, we introduce the appropriate functional framework to work with the fractional Laplacian and present some technical results. 

We start by giving a rigorous definition of the fractional Laplace operator. To this end, for any $0<s<1$, we consider the space 
\begin{align*}
	\mathcal L^s(\RR^d):=\left\{u:\RR^d\to\mathbb R \;\mbox{ measurable}:\;\int_{\Omega}\frac{|u(x)|}{(1+|x|)^{d+2s}}\;dx<+\infty\right\}.
\end{align*}
For $u\in \mathcal L^s(\mathbb R^d)$ and $\varepsilon>0$, we let
\begin{align*}
	(-\Delta)_\varepsilon^su(x):=C_{d,s}\int_{\{y\in\mathbb R^d:\;|x-y|\ge \varepsilon\}}\frac{u(x)-u(y)}{|x-y|^{d+2s}}\;dy,\;\;x\in\mathbb R^d,
\end{align*}
where the normalization constant is given by
\begin{align*}
	C_{d,s}:=\frac{s2^{2s}\Gamma\left(\frac{2s+d}{2}\right)}{\pi^{\frac d2}\Gamma(1-s)},
\end{align*}
and $\Gamma$ is the Euler Gamma function. The fractional Laplace operator $(-\Delta)^s$ is then defined for every $u\in \mathcal L^s(\mathbb R^d)$ by the formula
\begin{align}\label{eq:App1_FLO}
	(-\Delta)^su(x):=C_{d,s}\mbox{P.V.}\int_{\mathbb R^d}\frac{u(x)-u(y)}{|x-y|^{d+2s}}\;dy=\lim_{\varepsilon\downarrow 0}(-\Delta)_\varepsilon^su(x),\quad x\in\mathbb R^d,
\end{align}
provided that the limit exists for a.e. $x\in\mathbb R^d$. 

Let us now introduce the appropriate function spaces needed to work with the fractional Laplacian, that is, the fractional order Sobolev spaces. In what follows, we will only provide the definitions and some relevant properties. More complete presentations can be found in several references, including but not limited to \cite{adams2003sobolev,di2012hitchhiker,grisvard2011elliptic,lions1968problemes,warma2015the}.

Let $\Omega\subset\RR^d$ ($d\ge 1$) be a bounded open set. For any $0<s<1$, we define the fractional order Sobolev space 
\begin{align*}
	H^s(\Omega):=\left\{u\in L^2(\Omega):\;\int_{\Omega}\int_{\Omega}\frac{|u(x)-u(y)|^2}{|x-y|^{d+2s}}\;dxdy< +\infty\right\}
\end{align*}
and we endow it with the norm given by
\begin{align*}
	\|u\|_{H^s(\Omega)}:=\left(\int_{\Omega}|u|^2\;dx+\int_{\Omega}\int_{\Omega}\frac{|u(x)-u(y)|^2}{|x-y|^{d+2s}}\;dxdy\right)^{\frac 12}.
\end{align*}
We let
\begin{align*}
	H_0^s(\Omega) &:=\Big\{u\in H^s(\mathbb R^d):\; u=0\;\mbox{ in } \Omega^c\Big\} =\Big\{u\in H^s(\mathbb R^d):\operatorname{supp}[u]\subset\overline{\Omega}\,\Big\}
\end{align*}
and notice that if $0<s\ne 1/2<1$, and $\partial\Omega$ is Lipschitz, then by \cite[Chapter 1]{grisvard2011elliptic}, 
\begin{align*}
	H_0^s(\Omega)=\overline{\mathcal D(\Omega)}^{H^s(\Omega)},
\end{align*}
where $\mathcal D(\Omega)$ denotes the space of all continuous infinitely differentiable functions with compact support in $\Omega$. In that case, we endow $H_0^s(\Omega)$ with the $H^s(\Omega)$-norm. Otherwise one uses the norm of $H^s(\RR^d)$.

We denote by $H^{-s}(\Omega) := (H_0^s(\Omega))^\star$ the dual space of $H_0^s(\Omega)$ with respect to the pivot space $L^2(\Omega)$ so that the following continuous embeddings hold: 
\begin{align*} 
	H_0^s(\Omega)\hookrightarrow L^2(\Omega)\hookrightarrow H^{-s}(\Omega).
\end{align*} 

Besides, let $\langle\cdot,\cdot\rangle_{-s,s}$ be the duality pairing between $H^{-s}(\Omega)$ and $H_0^s(\Omega)$. Finally, we denote by $H^s_{\rm loc}(\Omega)$ the space defined by
\begin{align*}
	H^s_{\rm loc}(\Omega):= \Big\{u\in L^2_{\rm loc}(\Omega) \; : \; uv\in H^s(\Omega) \mbox{ for all } v\in\mathcal D(\Omega)\Big\}.
\end{align*}

Let us now introduce the following Dirichlet problem associated with the fractional Laplace operator:
\begin{equation}\label{eq:App1_DP}
	\begin{cases}
		(-\Delta)^su=f & \mbox{ in }\;\Omega,
		\\
		u=0 & \mbox{ in }\; \Omega^c.
	\end{cases}
\end{equation}
We have the following definition of weak solutions.

\begin{definition}\label{def:App1_weakSol} 
Let $f\in H^{-s}(\Omega)$. A function $u\in H_0^s(\Omega)$ is said to be a weak solution of the Dirichlet problem \eqref{eq:App1_DP} if the equality 
\begin{align*}
	\frac{C_{d,s}}{2}\int_{\mathbb R^d}\int_{\mathbb R^d}\frac{(u(x)-u(y))(v(x)-v(y))}{|x-y|^{d+2s}}\;dxdy=\langle f,v\rangle_{-s,s},
\end{align*}
holds for every $v\in H_0^s(\Omega)$.
\end{definition}

The existence and uniqueness of weak solutions to \eqref{eq:App1_DP} is a direct consequence of the classical Lax-Milgram Theorem. In particular, we have the following result (see e.g. \cite[Proposition 2.1]{biccari2017local} or \cite[Theorem 12]{leonori2015basic}).

\begin{proposition}\label{prop:App1_WPprop}
Let $\Omega\subset \RR^d$ be a bounded open set and $0<s<1$. Then for every $f\in H^{-s}(\Omega)$, the Dirichlet problem \eqref{eq:App1_DP} has a unique weak solution $u\in H_0^s(\Omega)$ in the sense of Definition \ref{def:App1_weakSol}. In addition, there is a constant $C>0$ such that
\begin{align*}
	\norm{u}{H_0^s(\Omega)} \leq C\norm{f}{H^{-s}(\Omega)}.
\end{align*}
\end{proposition}

We have the following maximal local regularity result for weak solutions of \eqref{eq:App1_DP}, which has been proved in \cite[Theorem 1.3]{biccari2017local} (see also \cite{grubb2015fractional}).

\begin{proposition}\label{regProp}
Let $\Omega\subset \RR^d$ be a bounded open set and $0<s<1$. Let $f\in L^2(\Omega)$ and $u\in H_0^s(\Omega)$ the unique weak solution of \eqref{eq:App1_DP} in the sense of Definition \ref{def:App1_weakSol}. Then $u\in H^{2s}_{\rm loc}(\Omega)$.
\end{proposition}

Next, we introduce the realization in $L^2(\Omega)$ of the fractional Laplacian with the zero Dirichlet exterior condition, that is, the operator 
\begin{align}\label{eq:App1_FLrealization}
	& D((-\Delta)_D^s):=\Big\{u\in H_0^s(\Omega):\; ((-\Delta)^su)|_{\Omega}\in L^2(\Omega)\Big\}\notag
	\\
	&(-\Delta)_D^su=((-\Delta)^su)|_{\Omega}\;\mbox{ a.e. in }\;\Omega.
\end{align}

It is well-known (see, e.g., \cite{claus2020realization}) that $-(-\Delta)_D^s$ generates a strongly continuous sub-markovian semi-group on $L^2(\Omega)$. Moreover, it has been shown in \cite{servadei2014spectrum} that $(-\Delta)_D^s$ has a compact resolvent, hence, it has a discrete spectrum which is formed with eigenvalues $(\lambda_j)_{j\in\NN}$ satisfying
\begin{align*}
	0<\lambda_1\le\lambda_2\le\cdots\le\lambda_j\le\cdots\;\mbox{ and }\;\lim_{j\to+\infty}\lambda_j=+\infty.
\end{align*}

We denote by $(\phi_j)_{j\in\NN}$ the normalized eigenfunctions associated with the eigenvalues $(\lambda_j)_{j\in\NN}$, i.e., the solutions of the Dirichlet problem
\begin{equation*}
	\begin{cases}
		(-\Delta)^s_D\phi_j=\lambda_j\phi_j & \mbox{in }\;\Omega,
		\\
		\phi_j=0 & \mbox{in }\; \Omega^c.
	\end{cases}
\end{equation*}
Finally, let us give our notion of solutions to the heat equation associated with the fractional Laplacian.

\begin{definition}\label{def:App1_weakSolHeat}
Let $y_0\in L^2(\Omega)$, $f\in L^2((0,T);H^{-s}(\Omega))$, and consider the following parabolic system
\begin{equation}\label{eq:App1_HeatControl}
	\begin{cases}
		y_t+(-\Delta)^sy= f\;\;&\mbox{ in } \Omega\times (0,T),
		\\
		y\equiv 0&\mbox{ in }\Omega^c\times(0,T),
		\\
		y(\cdot,0)= y_0&\mbox{ in }\Omega.
	\end{cases}
\end{equation}
We say that $y\in C([0,T];L^2(\Omega))\cap L^2((0,T);H_0^s(\Omega))\cap H^1((0,T);H^{-s}(\Omega))$ is a finite energy solution to \eqref{eq:App1_HeatControl} if $y(\cdot,0)=y_0$ a.e. in $\Omega$ and the identity
\begin{align*}
	\int_0^T \langle f,w\rangle_{-s,s}&\,dt = \int_0^T \langle y_t,w\rangle_{-s,s}\,dt 
	\\
	&+ \frac{C_{d,s}}{2}\int_0^T\int_{\RR^d}\int_{\RR^d}\frac{(y(x,t)-y(z,t))(w(x)-w(z))}{|x-z|^{d+2s}}\,dxdzdt\notag
\end{align*}
holds, for every $w\in H_0^s(\Omega)$.
\end{definition}

\noindent We have the following well-posedness result (see \cite[Chapter 10, Theorem 10.9]{brezis2010functional} or \cite[Theorem 26]{leonori2015basic}).

\begin{proposition}
Assume that $f\in L^2((0,T);H^{-s}(\Omega))$. Then, for every initial datum $y_0\in L^2(\Omega)$ the fractional heat equation \eqref{eq:App1_HeatControl} has a unique finite energy solution $y$ given by
\begin{align*}
	y(\cdot,t) = e^{-(-\Delta)^s_D t}y_0 + \int_0^t e^{-(-\Delta)^s_D (t-\tau)}f(\cdot,\tau)\;d\tau,
\end{align*}
where $e^{-(-\Delta)^s_D t}$, $t\ge 0$, is the strongly continuous semi-group on $H^{-s}(\Omega)$ generated by the operator $-(-\Delta)_D^s$.
\end{proposition}

	\section{The fractional Laplace operator with exterior conditions}\label{appendixB}

We present here a general overview of elliptic and parabolic problems associated with the fractional Laplace operator on a bounded domain $\Omega\subset\RR^d$ with non-zero exterior data. We first consider the following elliptic Dirichlet problem:
\begin{equation}\label{eq:App2_DP}
	\begin{cases}
		(-\Delta)^su =f & \mbox{in }\Omega 
		\\
		u=g & \mbox{in } \Omega^c.
	\end{cases}
\end{equation}

Let us define our notion of solutions for \eqref{eq:App2_DP}. To this end, we need first to introduce the non-local normal derivative $\mathcal{N}_s$, defined for all $u \in H^s(\RR^d)$ as
\begin{align}\label{eq:App2_nonlocalDer}
   \mathcal{N}_s u(x) := C_{d,s} \int_\Om \frac{u(x)-u(y)}{|x-y|^{d+2s}}\;dy, \quad x \in \RR^d \setminus \overline{\Omega}. 
\end{align}

Clearly, $\mathcal{N}_s$ is a non-local operator. Moreover, it is well defined on $H^s(\RR^d)$ as the following result shows (see \cite[Lemma A.2]{ghosh2020calderon} for the proof).

\begin{lemma}\label{lem:App2_Nmap}
The non-local normal derivative $\mathcal{N}_s$ maps $H^s(\RR^d)$ continuously into $ H^s_{\rm loc}(\Omega^c)\subset L^2_{\rm loc}(\Omega^c)$. 
\end{lemma}

Even if $\mathcal{N}_s$ is defined on the unbounded domain $\RR^d \setminus \overline{\Omega}$, it is still known as the \textit{normal} derivative. This is due to its similarity with the classical normal derivative as the following result taken from \cite{DiRoVa} shows.

\begin{proposition}\label{prop:App2_prop}
The following assertions hold.
\begin{itemize}
   \item[1.] {\bf Divergence theorem.} Let $u\in C^2(\RR^d)$ vanishing at $\pm \infty$. Then
   \begin{align*} 
        \int_\Omega (-\Delta)^s u\;dx = -\int_{\Omega^c} \mathcal{N}_s u\;dx.
   \end{align*}
      
   \item[2.] {\bf Integration by parts formula.} Let $u\in H^s(\RR^d)$ be such that $(-\Delta)^su \in L^2(\Omega)$ and $\mathcal N_s u\in L^2(\Omega^c)$. Then, for every $v\in H^s(\RR^d)$ we have 
   \begin{align*}
   		\int_\Omega v(-\Delta)^s u\;dx =&\; \frac{C_{d,s}}{2} \int\int_{\RR^{2d}\setminus(\Omega^c)^2} \frac{(u(x)-u(y))(v(x)-v(y))}{|x-y|^{d+2s}} \;dxdy \notag
   		\\
    	&\; - \int_{\Omega^c} v\mathcal{N}_s u\;dx,
   \end{align*}
   where $\RR^{2d}\setminus(\Omega^c)^2 = (\Omega\times\Omega)\cup(\Omega\times\Omega^c)\cup(\Omega^c\times\Omega)$. 
     
   \item[3.] {\bf Limit as $s\uparrow 1^-$.} Let $u,v\in C^2(\RR^d)$ vanishing at $\pm\infty$. Then
   \begin{align*}
   		\lim_{s\uparrow 1^-}\int_{\Omega^c} v \mathcal{N}_s u\;dx = \int_{\partial\Omega} v\frac{\partial u}{\partial\nu}\,d\sigma. 
   \end{align*} 
\end{itemize}
\end{proposition}

\noindent We are now ready to introduce the notion of transposition solutions to the Dirichlet problem \eqref{eq:App2_DP}.

\begin{definition}\label{def:App2_vweak_d}
Let $g \in L^2(\Omega^c)$ and $f \in H^{-s}(\overline\Om)$. A function $u \in L^2(\RR^d)$ is said to be a solution by transposition to \eqref{eq:App2_DP} if the identity
\begin{align}\label{eq:App2_transposition}
	\int_{\Om} u (-\Delta)^s v\;dx = \langle f,v \rangle_{-s,s} - \int_{\Omega^c} g \mathcal{N}_s v\;dx,
\end{align}
holds for every $v \in V := \big\{v \in H^{s}_0(\Om) \;:\; (-\Delta)^s v \in L^2(\Om)\big\}$. 
\end{definition}
 
Moreover, we have the following existence and uniqueness result of solutions by transposition. We refer to \cite[Theorem 3.5]{antil2019external} for the proof.

\begin{theorem}\label{thm:App2_vwdexist}
Let $f \in H^{-s}(\Om)$ and $g \in L^2(\Omega^c)$. Then, \eqref{eq:App2_DP} has a unique solution by transposition $u\in L^2(\RR^d)$, and there is a constant $C>0$ such that
\begin{align*}
	\|u\|_{L^2(\RR^N)} \le C\left(\|f\|_{H^{-s}(\overline\Om)} + \|g\|_{L^2(\Om^c)} \right).
\end{align*}
\end{theorem}
 
\noindent We now move to the parabolic problem associated to the fractional Laplacian with exterior condition.  
\begin{equation}\label{eq:App2_Sd}
	\begin{cases}
    	y_t + (-\Delta)^s y = 0 & \mbox{in } \Omega\times(0,T)
    	\\
        y = g & \mbox{in } \Om^c\times(0,T)
        \\
        y(\cdot,0) = y_0 & \mbox{in } \Om.
	\end{cases}                
\end{equation}
Let $y_0\in L^2(\Omega)$, $g\in L^2((0,T);H^s(\Omega^c))$ and consider the following two systems:
\begin{equation}\label{eq:App2_systXi}
	\begin{cases}
		\xi_t + (-\Delta)^s \xi = 0 & \mbox{in } \Omega\times(0,T)
		\\
		\xi = 0 & \mbox{in } \Om^c\times(0,T)
		\\
		\xi(\cdot,0) = y_0 & \mbox{in } \Om.
	\end{cases}                
\end{equation}
and
\begin{equation}\label{eq:App2_systZeta}
	\begin{cases}
		z_t + (-\Delta)^s z = 0 & \mbox{in } \Omega\times(0,T)
		\\
		z = g & \mbox{in } \Om^c\times(0,T)
		\\
		z(\cdot,0) = 0 & \mbox{in } \Om.
	\end{cases}                
\end{equation}

Then, the solution of \eqref{eq:App2_Sd} is given by $y = \xi+z$. Moreover, noticing that \eqref{eq:App2_systXi} can be cast as a Cauchy problem for the operator $(-\Delta)_D^s$ we introduced in \eqref{eq:App1_FLrealization}, using semi-group theory and the spectral theorem, one has the following result.
\begin{theorem}\label{thm:App2_existenceXi}
Let $(\phi_k)_{k\in\NN}$ be the normalized eigenfunctions of the operator $(-\Delta)_D^s$ associated with the eigenvalues $(\lambda_k)_{k\in\NN}$. For every $y_0\in L^2(\Omega)$, define $y_{0,k}:=\langle y_0,\phi_k\rangle_{L^2(\Omega)}$. Then, there is a unique function 
\begin{align*}
	\xi\in C([0,T];L^2(\Omega))\cap L^2((0,T);H_0^s(\Omega))\cap H^1((0,T);H^{-s}(\Omega))
\end{align*}
satisfying \eqref{eq:App2_systXi} which is given for a.e. $x\in\Omega$ and every $t\in [0,T]$ by
\begin{align*}
	\xi(x,t) = \sum_{j\geq 1} y_{0,k}e^{-\lambda_k t}\phi_k(x).	
\end{align*}
\end{theorem}

We now consider the non-homogeneous exterior problem \eqref{eq:App2_systZeta}, for which we introduce the following notion of weak solution.
\begin{definition}
Let $g\in L^2((0,T);H^s(\Omega^c))$. By a weak solution of \eqref{eq:App2_systZeta} we mean a function $z\in L^2((0,T);H^s(\mathbb R))$ such that $z=g$ a.e. in $\Omega^c\times (0,T)$ and the identity
\begin{align}\label{eq:App2_weakSolHeat}
	\int_0^T \langle -w_t+(-\Delta)^s w,z\rangle_{-s,s}\,dt = \int_\Omega z(x,T)w(x,T)\,dx + \int_0^T\int_{\Omega^c} g\mathcal N_s w\,dxdt
\end{align}
holds for every $w\in C([0,T];L^2(\Omega))\cap L^2((0,T);H_0^s(\Omega))\cap H^1((0,T);H^{-s}(\Omega))$ with $\mathcal N_s w\in L^2((0,T)\times\Omega^c)$.
\end{definition}

\noindent We then have the following existence result (see \cite{warma2019approximate}).
\begin{theorem}\label{thm:App2_existenceZeta}
For every $g \in L^2((0,T);H^s(\Om^c))$ , the system \eqref{eq:App2_systZeta} has a unique weak solution $z\in L^2((0,T);H^s(\mathbb R))$ given by 
\begin{align*}
	z(x,t)=\sum_{k\geq 1}\left(\int_0^t\big(g(\cdot,t-\tau),\mathcal N_s\phi_k\big)_{L^2(\Omega^c)} e^{-\lambda_k\tau}\;d\tau\right)\phi_k(x). 
\end{align*}
\end{theorem}

Hence, from Theorems \ref{thm:App2_existenceXi} and \ref{thm:App2_existenceZeta}, we finally have the following result of existence and series representation of solutions to \eqref{eq:App2_Sd}.

\begin{theorem}\label{thm:App2_existenceY}
For every $y_0\in L^2(\Omega)$ and $g\in L^2((0,T);H^s(\Om^c))$ , the system \eqref{eq:App2_Sd} has a unique weak solution $y\in L^2((0,T)\times\RR^d)$ given by 
\begin{align*}
	y(x,t)=\sum_{k\geq 1} y_{0,k}e^{-\lambda_k t}\phi_k + \sum_{k\geq 1}\left(\int_0^t\big(g(\cdot,t-\tau),\mathcal N_s\phi_k\big)_{L^2(\Omega^c)} e^{-\lambda_k\tau}\;d\tau\right)\phi_k(x). 
\end{align*}
\end{theorem}

Let us conclude this appendix by quickly discussing the following fractional heat equation with exterior Robin conditions
\begin{equation}\label{eq:App2_Robin}
	\begin{cases}
		y_t + (-\Delta)^s y = 0 & \mbox{in } \Omega\times(0,T)
		\\
		\mathcal N_sy + \kappa y = \kappa g & \mbox{in } \Om^c\times(0,T)
		\\
		y(\cdot,0) = y_0 & \mbox{in } \Om
	\end{cases}                
\end{equation}
which plays a fundamental role in the numerical simulations of Section \ref{sec:4}. In \eqref{eq:App2_Robin}, $\kappa \in L^1(\Omega^c)\cap L^\infty(\Omega^c)$ is a non-negative function.

To define the notion of solutions for the Robin problem \eqref{eq:App2_Robin}, we first need to introduce the Sobolev space
\begin{align*}
	H^s_\kappa (\Omega):= \left\{y:\RR^d\to\RR \mbox{ measurable },\; \norm{y}{H^s_\kappa(\Omega)}<+\infty\right\},
\end{align*} 
where
\begin{align*}
	\norm{y}{H^s_\kappa(\Omega)}:= \left( \norm{y}{L^2(\Omega)}^2 + \norm{\kappa^{\frac 12}y}{L^2(\Omega^c)}^2 + \int_{\RR^{2d}}\frac{|y(x)-y(z)|}{|x-z|^{d+2s}}\,dxdz\,\right)^{\frac 12}.
\end{align*}

We know from \cite[Proposition 3.1]{DiRoVa} that $H^s_\kappa(\Omega)$ is a Hilbert space. We denote with $H^{-s}_\kappa(\Omega) = (H^s_\kappa(\Omega))^\star$ its dual with respect to the pivot space $L^2(\Omega)$.

\begin{definition}\label{def:App2_weakSolRobin}
Let $y_0\in L^2(\Omega)$, $g\in L^2(\Omega^c\times (0,T))$ and $\kappa \in L^1(\Omega^c)\cap L^\infty(\Omega^c)$ non-negative. A function $y\in L^2((0,T);H^s_\kappa(\Omega))\cap H^1((0,T);H^{-s}_\kappa(\Omega))$ is said to be a weak solution of \eqref{eq:App2_Robin} if $y(\cdot,0) = y_0$ a.e. in $\Omega$ and the identity
\begin{align}\label{eq:App2_weakRobin}
	\int_0^T \langle y_t,v\rangle_{-s,s}\,dt + \int_0^T\mathcal F(y,v)\,dt + \int_0^T\int_{\Omega^c} \kappa yv\,dxdt = \int_0^T\int_{\Omega^c} \kappa gv\,dxdt,
\end{align}
holds for every $v\in H^s_\kappa (\Omega)$, where $\mathcal F(y,v)$ denotes the bilinear form
\begin{align*}
	\mathcal F(y,v):= \int_{\RR^{2d}\setminus(\Omega^c)^2} \frac{(y(x)-y(z))(v(x)-v(z))}{|x-z|^{d+2s}}\,dxdz.
\end{align*}
\end{definition}

\noindent Finally, we have the following existence result for \eqref{eq:App2_Robin} (see \cite[Theorem 3.11]{antil2020external}).
\begin{theorem}\label{thm:App2_existenceRobin}
Let $\kappa \in L^1(\Omega^c)\cap L^\infty(\Omega^c)$ be non-negative. Then for every $y_0\in L^2(\Omega)$ and $g\in L^2(\Omega^c\times (0,T))$, there exists a unique weak solution $y\in L^2((0,T);H^s_\kappa(\Omega))\cap H^1((0,T);H^{-s}_\kappa(\Omega))$ of \eqref{eq:App2_Robin}.
\end{theorem}

}

\bibliographystyle{abbrv}
\bibliography{biblioHandbook}
\end{document}